%% file: rapport.tex
\documentclass[english,a4paper]{article}
\usepackage[english]{babel}
\usepackage{nada-tr}
\usepackage[latin1]{inputenc}
\usepackage{latexsym}
\usepackage[final]{graphicx}
\usepackage{amsfonts}
\usepackage{amsmath}
\usepackage{amssymb}
\usepackage{subfigure}
\usepackage{multirow}

\title{A Stochastic Phase-Field Model \\
  Computed From
  Coarse-Grained Molecular Dynamics}
\author{Erik von Schwerin\\
  KTH -- Computer Science and Communication \\
  SE-100 44 Stockholm\\
  SWEDEN}
\date{October 8, 2007}
\coverauthor{Erik von Schwerin}
\authoremail{schwerin@csc.kth.se}

\begin{document}

\maketitle

\newcommand{\R}{\mathbb{R}}
\newcommand{\backdefeq}{=\colon\;}
\newcommand{\defeq}{\colon\!\!=}
\newcommand{\dX}{\nabla_{X_j}}
\newcommand{\dx}{\nabla_x}
\newcommand{\dxett}{\frac{\partial}{\partial x_1}}
\newcommand{\ddxett}{\frac{\partial^2}{\partial x_1^2}}
\newcommand{\transpose}[1]{{\ensuremath{#1^\mathrm{T}}}}
\newcommand{\funcx}{\psi} 
\newcommand{\funcxX}[1]{\psi(\cdot;{#1})} 
\newcommand{\I}[1]{{\ensuremath{\mathbf 1_{#1}}}}
\newcommand{\E}{{\ensuremath{\mathrm E}}}
\newcommand{\Var}{{\ensuremath{\mathrm{Var}}}}
\newcommand{\average}[2]
{{\ensuremath{\mathcal {A}_{#2}\left(#1\right)}}}
\newcommand{\variance}[2]
{{\ensuremath{\mathrm{Var}_{#2}\left(#1\right)}}}
\newcommand{\sampleset}{{\ensuremath{\mathcal{S}}}}
\newcommand{\conv}{\xrightarrow}
\newcommand{\Ito}{It\^{o}}
\newcommand{\kb}{{\ensuremath{k_\mathrm{B}}}}
\newcommand{\cv}{{\ensuremath{c_V}}}
\newcommand{\latheat}{{\ensuremath{L}}}
\newcommand{\thermcond}{{\ensuremath{\lambda}}}
\newcommand{\nrp}{{\ensuremath{N}}}
\newcommand{\totpot}{{\ensuremath{U}}}
\newcommand{\toten}{{\ensuremath{E}}}
\newcommand{\pairpot}{{\ensuremath{\Phi}}}
\newcommand{\force}{{\ensuremath{F}}}
\newcommand{\pairforce}{{\ensuremath{f}}}
\newcommand{\divforce}{{\ensuremath{G}}}
\newcommand{\divpairforce}{{\ensuremath{g}}}
\newcommand{\no}{{\mbox{}}}
\newcommand{\mdpos}[2][t]{{\ensuremath{X_{#2}^{#1}}}}
\newcommand{\mdposgen}{{\ensuremath{\mdpos[\no]{\no}}}}
\newcommand{\mdvel}[2][t]{{\ensuremath{v_#2^#1}}}
\newcommand{\mdposd}[2][n]{{\ensuremath{\bar{X}_{#2}^{#1}}}}
\newcommand{\Smo}{Smoluchowski}
\newcommand{\sumall}[2][\nrp]{{\ensuremath{\sum_{#2=1}^#1}}}
\newcommand{\sumneq}[3][\nrp]
        {{\ensuremath{\sum_{#2\neq#3,#2=1}^#1}}}
\newcommand{\rdf}[1][r]{{\ensuremath{g(#1)}}}
\newcommand{\density}{{\ensuremath{\rho}}}
\newcommand{\damping}{{\ensuremath{\tau}}}
\newcommand{\kbT}{{\ensuremath{\gamma}}}
\newcommand{\fenfun}{{\ensuremath{F}}}
\newcommand{\fendty}{{\ensuremath{\tilde{f}}}}
\newcommand{\monoton}{{\ensuremath{g}}}
\newcommand{\dwell}{{\ensuremath{f}}}
\newcommand{\graden}{{\ensuremath{K_1}}}
\newcommand{\barrier}{{\ensuremath{K_0}}}
\newcommand{\mobility}{{\ensuremath{M}}}
\newcommand{\nrw}{{\ensuremath{M}}}
\newcommand{\indepw}{{\ensuremath{\widetilde{W}}}}
\newcommand{\potenp}[1]{{\ensuremath{m_#1}}}
\newcommand{\pfgen}{{\ensuremath{\phi}}}
\newcommand{\pfen}{{\ensuremath{m}}}
\newcommand{\pfent}[1][t]{{\ensuremath{m^#1}}}
\newcommand{\pfcg}{{\ensuremath{m_{\mathrm{cg}}}}}
\newcommand{\pfcgt}[1][t]{{\ensuremath{m_{\mathrm{cg}}^#1}}}
\newcommand{\driftmd}{{\ensuremath{\alpha}}}
\newcommand{\diffumd}{{\ensuremath{\beta}}}
\newcommand{\driftcg}{{\ensuremath{a}}}
\newcommand{\diffucg}{{\ensuremath{b}}}
\newcommand{\driftcgx}{{\ensuremath{\overline{a}}}}
\newcommand{\driftcgxs}{{\ensuremath{\overline{A}}}}
\newcommand{\diffucgx}{{\ensuremath{\overline{b}}}}
\newcommand{\diffucgxmat}{{\ensuremath{B}}}
\newcommand{\diffumat}{{\ensuremath{\overline{B}}}}
\newcommand{\diffumatinst}{{\ensuremath{\tilde{B}}}}
\newcommand{\eigenbasemat}{{\ensuremath{V}}}
\newcommand{\poseigenmat}{{\ensuremath{V_+}}}
\newcommand{\eigenvalmat}{{\ensuremath{\Lambda}}}
\newcommand{\poseigenval}{{\ensuremath{\Lambda_+}}}
\newcommand{\driftone}{{\ensuremath{a_1}}}
\newcommand{\driftzero}{{\ensuremath{a_0}}}
\newcommand{\molli}{{\ensuremath{\eta}}}
\newcommand{\molliscale}{{\ensuremath{\epsilon}}}
\newcommand{\temp}{{\ensuremath{T}}}
\newcommand{\Tmelt}{{\ensuremath{T_\mathrm{M}}}}
\newcommand{\vol}{{\ensuremath{V}}}
\newcommand{\pressure}{{\ensuremath{P}}}
\newcommand{\smooth}{{\ensuremath{y}}}
\newcommand{\tend}{{\ensuremath{\mathcal{T}}}}
\newcommand{\dual}{{\ensuremath{\overline{u}}}}

\begin{abstract}
Results are presented from numerical experiments aiming at the
computation of stochastic phase-field models for phase transformations
by coarse-graining molecular dynamics.
The studied phase transformations occur between a solid crystal and a
liquid.
Nucleation and growth, sometimes dendritic, of crystal grains in a
sub-cooled liquid is determined by diffusion and convection of heat, on 
the macroscopic level, and by interface effects, where the width of
the solid--liquid interface is on an atomic length-scale.
Phase-field methods are widely used in the study of the continuum
level time evolution of the phase transformations; they introduce an
order parameter to distinguish between the phases. 
The dynamics of the order parameter is modelled by an Allen--Cahn
equation and coupled to an energy equation, where the latent heat at
the phase transition enters as a source term.
Stochastic fluctuations are sometimes added in the coupled system of
partial differential equations to introduce nucleation and to get
qualitatively correct behaviour of dendritic side-branching.
In this report the possibility of computing some of the Allen--Cahn
model functions from a microscale model is investigated. The
microscopic model description of the material by stochastic,
Smoluchowski, dynamics is considered given. 
A local average of contributions to the potential energy in the micro
model is used to determine the local phase, and a stochastic
phase-field model is computed by coarse-graining the molecular
dynamics. 
Molecular dynamics simulations on a two phase system at
the melting point are used to compute a double-well
reaction term in the Allen--Cahn equation
and a diffusion matrix
describing the noise in the coarse-grained phase-field.
\end{abstract}

~\thanks{This work was supported by 
  the  Swedish Foundation for Strategic Research grant A3~02:123,
  "Mathematical theory and simulation tools for phase
  transformations is materials".}

\section{Introduction}
\label{sec:Intro}

Phase-field methods are widely used for modelling phase
transformations in materials on the continuum level and exist in many
different versions for different applications.
In this report the considered phase transformation occurs in a single
component system with a solid and a liquid phase.

The phase-field model of solidification studied here
is a coupled system of partial 
differential equations for the temperature, \temp, and a phase-field,
\pfgen, which is an order parameter used to distinguish
between the solid and the liquid subdomains. Two different values,
$\pfgen_s$ and $\pfgen_l$, are equilibrium values of the phase-field in
solid and liquid respectively. The phase-field varies continuously
between the two values and the interface between solid and liquid, at
a time $t$, is defined as a level surface of the phase-field; for
example  
$\{x\in\R^d\,\colon\,\pfgen(x,t) = 0.5(\pfgen_s+\pfgen_l)\}$. 
From a computational point of view the implicit definition of the
phases in the phase-field method, as in the level set
method~\cite{Level_set_Osher,Level_set_Sethian}, is an 
advantage over sharp interface methods, since it avoids the explicit
tracking of the interface.
A local change of the phase-field from $\pfgen_l$ to $\pfgen_s$ in a
subdomain translates into solidification of that region with a
corresponding release of latent heat and the reverse change from
$\pfgen_s$ to $\pfgen_l$ means melting which requires energy. 
The release or absorption of latent heat is modelled as a continuous
function of \pfgen\ so that the energy released when a unit volume
solidifies is $\latheat(\monoton(\pfgen_l)-\monoton(\pfgen_s))$, where
\latheat\ is the latent heat and 
$\monoton(\pfgen)$ is a model function, monotone with 
$\monoton(\pfgen_s) = 0$, $\monoton(\pfgen_l) = 1$, 
$\monoton'(\pfgen_s) = 0$, and $\monoton'(\pfgen_l) = 0$.
Then the energy equation for a unit volume becomes 
a heat equation with a source term
\begin{align*}
  \frac{\partial}{\partial t}
  \left( \cv \temp + \latheat \monoton(\pfgen)\right)
  & = \nabla\cdot\left(\thermcond \nabla \temp\right),
\end{align*}
where \cv\ is the heat capacity at constant volume and
\thermcond\ is the thermal conductivity.
Here, and in the following, the usual notation for differentiation
with respect to the spatial variables is applied, with $\nabla$ and
$\nabla\cdot$ denoting the gradient and the divergence respectively. 
The phase-field, and the related model function \monoton, are
exceptional in the energy equation in the sense that, while all the
other quantities are standard physical quantities on the macroscopic 
level, the phase-field need not be associated with a measurable quantity.
A phenomenological model of the phase change is given by the energy
equation coupled to the Allen-Cahn equation 
\begin{align}
  \label{eq:deterministic_AllenCahn}
  \frac{\partial \pfgen}{\partial t}
  & = 
  \nabla\cdot\left(k_1 \nabla\pfgen\right)
  - 
  k_2\Big(\dwell'(\pfgen)+\monoton'(\pfgen)k_3 (\Tmelt-\temp) \Big) 
\end{align}
for the time evolution of the phase-field; here \Tmelt\ denotes the
melting point, $k_1$, $k_2$, and $k_3$, are positive model parameters
($k_1$ may be an anisotropic matrix introducing directional dependence
on the growth of the solid),
and the model function \dwell\ is a double well potential with minima
at $\pfgen_s$ and $\pfgen_l$. 
Standard examples of the model functions are 
\begin{align*}
  \dwell(\pfgen) & = -\frac{1}{2}\pfgen^2 + \frac{1}{4}\pfgen^4,
  &
  \monoton(\pfgen) & = \frac{15}{16}\left(
      \frac{1}{5}\pfgen^5 - \frac{2}{3}\pfgen^3 + \pfgen\right)
    + \frac{1}{2},
\end{align*}
when $\pfgen_s=-1$ and $\pfgen_l=1$. 
By construction of the model functions, the reaction term in the
Allen-Cahn equation vanishes where $\pfgen=\pfgen_s$ or
$\pfgen=\pfgen_l$ independently of the temperature. Since the
diffusion term is zero for any constant function the two constant
phase-fields $\pfgen\equiv\pfgen_s$ and $\pfgen\equiv\pfgen_l$ are
stationary solutions to the Allen-Cahn equation for all temperatures. 
This means, for example, that nucleation of solid in a region of
subcooled liquid can not occur in a phase-field modelled by the
deterministic Allen-Cahn equation above. The effect of nucleation can
be introduced in the model by adding a noise term in the Allen-Cahn
equation, giving a stochastic partial differential equation.
Simulation of dendrite growth in an subcooled liquid is another
example where the deterministic system is inadequate; its solutions
fail to develop the side branches seen to form in real dendrites as
the tips grow. 
Stochastic phase-field models where noise is added to either one, or
both, of the Allen-Cahn equation and the energy equation are used to
include the effect of side branching; 
see for example~\cite{phasefield}.

\begin{figure}[hbp]
  \centering
  \includegraphics[width=3.35in]{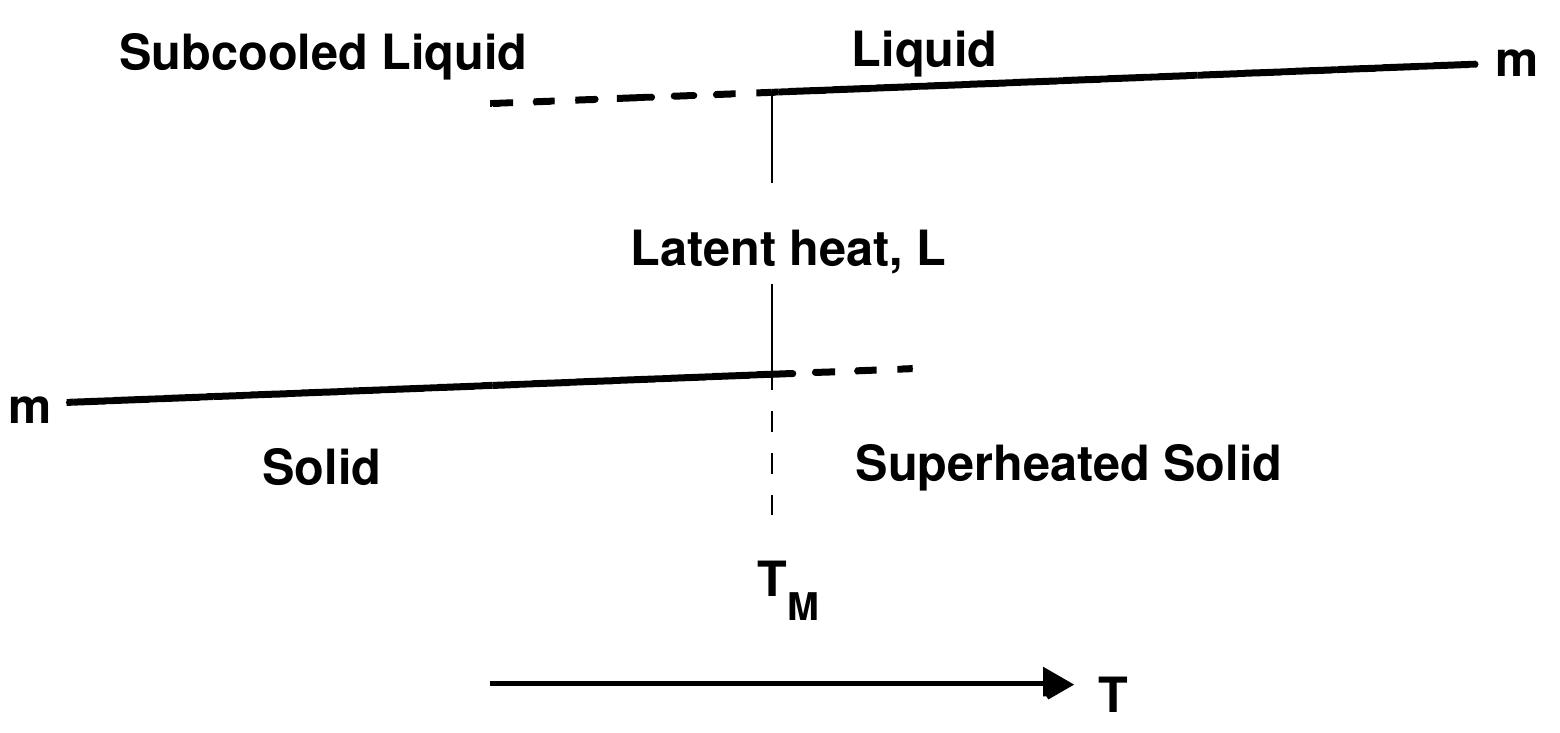}
  \caption{Schematic picture of $\pfen(\temp)$ for a pure liquid (top
    curve) and a pure solid (bottom curve) and the latent heat as the
    jump in \pfen\ at a phase transition.}
  \label{fig:m_of_T}
\end{figure}

The present report contains the results from numerical experiments
on a method presented and analysed in~\cite{Anders} and the rest of
this introduction summarises the ideas from~\cite{Anders} needed here.  
That report takes the stochastic phase-field model
\begin{subequations}
\label{eq:pfsyst}
\begin{align}
  \label{eq:energy}
  \frac{\partial}{\partial t}
  \left( \cv \temp + \latheat \monoton(\pfgen)\right)
  & = \nabla\cdot\left(\thermcond \nabla \temp\right),
  \\
  \label{eq:AllenCahn}
  \frac{\partial \pfgen}{\partial t}
  & = 
  \nabla\cdot\left(k_1 \nabla\pfgen\right)
  - 
  k_2\Big(\dwell'(\pfgen)+\monoton'(\pfgen)k_3 (\Tmelt-\temp) \Big) 
  + \text{noise},
\end{align}
\end{subequations}
as its starting point and asks whether it is possible to obtain 
the model functions and parameters,
$\dwell(\pfgen)$, $\monoton(\pfgen)$, $k_1$, $k_2$, $k_3$, and 
the noise, from computations on a microscale model.
To answer this question the phase-field, $\pfgen$, must be defined 
in terms of quantities computable on the microscale.
The microscopic model used for this purpose 
is a molecular dynamics model of \nrp\ particles 
in a microscopic domain $D$ in $\R^3$ where the motion of the 
particles is given by the Smoluchowski dynamics; see for
example~\cite{Majda_microfluid}. Thus, 
with $\mdpos{\no}\in\R^{3N}$ denoting the positions of all particles
in the system at the time $t$ and $\mdpos{i}\in\R^3$ the position of
particle $i$, the dynamics are given by the \Ito\ stochastic
differential equations
\begin{align}
  \label{eq:Smoluchowski}
  d\mdpos{i} = 
  -\nabla_{\mdpos[\no]{i}} \totpot(\mdpos{\no})\;dt 
  + \sqrt{2\kb \temp}\;dW_i^t,
  && i=1,2,\ldots,\nrp,
\end{align}
where \totpot\ is the total potential energy of the system, 
$\nabla_{\mdpos[\no]{i}}$ denotes the gradient with respect to the
position of particle $i$, \kb\ is the Boltzmann constant, 
and $W_i=\transpose{(W_{i,1},W_{i,2},W_{i,3})}$ are independent three
dimensional Brownian motions, with independent components. 
The macroscopic temperature, \temp, is a constant input parameter in
the microscopic model. 
We may identify the latent heat, in the macroscopic model, 
with the difference in total potential energy per unit volume of the
liquid and the solid at the melting point, in the microscopic model.
The idea is then to let the local contributions to the total potential
energy define the phase variable. 
Since the potential energy decreases with the temperature even in a
single phase system the equilibrium values of such a phase-field,
\pfen, unlike those of \pfgen, depend on the temperature; see
Figure~\ref{fig:m_of_T}. 
Assuming that in pure solid or pure liquid the phase-field, \pfen,
varies slowly, compared to the latent heat release, with the
temperature close to the melting point, the energy equation becomes 
\begin{align*}
  \frac{\partial}{\partial t}
  \left( \cv \temp + \pfen\right)
  & = \nabla\cdot\left(\thermcond \nabla \temp\right),
\end{align*}
where \cv\ and \thermcond\ are approximately the same as
in~\eqref{eq:energy} for $\temp\approx\Tmelt$.

For a model where the total potential energy of the system can be
naturally split into a sum of contributions arising from the
interaction of individual atoms with their environment,  
\begin{align}
  \label{eq:def_totpot}
  \totpot(\mdposgen) & = 
  \sumall{i} \potenp{i}(\mdposgen),
\end{align}
phase-fields can be introduced on the micro level by localised averages
of these contributions; a given configuration \mdposgen\ defines a 
phase-field $\pfen(\,\cdot\,;\mdposgen)\;\colon\;D\rightarrow\R$
through 
\begin{align}
  \label{eq:def_pfen}
  \pfen(x;\mdposgen) 
  & = \sumall{i} 
  \potenp{i}(\mdposgen) \molli(x-\mdpos[\no]{i}),
\end{align}
where the choice of mollifier, \molli, determines the spatial
smoothness of the phase-field. If, for example, the potential energy
is defined entirely by pairwise interactions 
\begin{align*}
  \totpot(\mdposgen) & = 
  \frac{1}{2}\sumall{i} \sumneq{k}{i} 
  \pairpot(\mdpos[\no]{i}-\mdpos[\no]{k}),
\end{align*}
as is common in simple molecular dynamics models, it is natural to let 
\begin{align*}
  \potenp{i}(\mdposgen) & = 
  \frac{1}{2}\sumneq{k}{i} 
  \pairpot(\mdpos[\no]{i}-\mdpos[\no]{k})
\end{align*}
be particle $i$'s contribution to the total potential energy.

With the definition~\eqref{eq:def_pfen} of the potential energy
phase-field, \pfen, and with the microscopic system defined
by~\eqref{eq:Smoluchowski} and~\eqref{eq:def_totpot}, \Ito's formula
gives a stochastic differential equation
\begin{align}
  \label{eq:sde_pfen}
  d\pfen(x;\mdpos{\no}) & = \driftmd(x;\mdpos{\no})\;dt + 
  \sumall{j} \sum_{k=1}^3\diffumd_{j,k}(x;\mdpos{\no})\;dW_{j,k}^t,
\end{align}
for \pfen\ evaluated in a point $x\in D$. The drift,
$\driftmd(x;\cdot)$, and the diffusions, $\diffumd_{j,k}(x;\cdot)$,
are explicitly known functions expressed in terms of the $\pfen_i$:s, 
the mollifier, \molli, and their derivatives up to second order.
While \pfen\ by definition is a continuous field it is still an atomic
scale quantity since it is defined in terms the particle positions
\mdpos{\no}. A macroscopic phase-field, similar to \pfgen\
in~\eqref{eq:pfsyst}, must lose both the dependence on the particle
positions, \mdpos{\no}, and the explicit dependence on the
microscale space variable $x$. To achieve this, a coarse-grained 
approximation $\pfcg(x)$ of $\pfen(x)$ is introduced as a solution of
a stochastic differential equation 
\begin{align}
  \label{eq:sde_pfcg}
  d\pfcgt(x) & = \driftcg(\pfcgt)(x)\;dt + 
  \sumall[\nrw]{j} \diffucg_j(\pfcgt)(x)\;d\indepw_j^t,  
\end{align}
where the independent Wiener processes $\indepw_j^t$, 
$j=1,2,\ldots,\nrw\ll\nrp$, also are independent of the Wiener processes
$W_i$ in the micro model.  
Here the drift and diffusion coefficient functions, $\driftcg(\pfcgt)$
and $\diffucg_j(\pfcgt)$, may depend on more information about the 
coarse-grained phase-field than just the point value; compare the
stochastic Allen-Cahn equation~\eqref{eq:AllenCahn}, where the
diffusion term in the drift contains second derivatives of the
phase-field. 

The choice of the coarse-grained drift and diffusion functions 
proceeds in two steps: first, finding a general form the
coarse-grained equation where the drift and diffusion coefficient
functions, defined as time averaged expected values of the microscopic
drift and diffusions over simulation paths, still depend on the micro
scale space variable, $x$; 
second, expressing the $x$ dependent coarse-grained drift and
diffusion coefficients by drift and diffusion functions depending only
on the phase-field \pfcg, using that \pfcg\ is a smooth monotone
function of $x$ in the interface.

In the first step, a coarse-grained stochastic differential equation
\begin{align*}
  d\pfcgt(x) & = \driftcgx(x)\;dt + 
  \sumall[\nrw]{j} \diffucgx_j(x)\;d\indepw_j^t,  
\end{align*}
is introduced by defining the drift
\begin{subequations}
\label{eq:def_coeffcgx}
\begin{align}
  \label{eq:def_driftcgx}
  \driftcgx(x) & = 
  \frac{1}{\tend}
  \E\biggl[\int_0^\tend \driftmd(x;\mdpos{\no})\,dt
  \;\Bigl\lvert\;
  \mdpos[0]{\no}=\mdpos[\no]{0}
  \biggr],
  && x\in D,
  \intertext{and choosing a diffusion matrix that fulfil}
  \label{eq:def_diffucgx}
  \sumall[\nrw]{j}\diffucgx_j(x)\diffucgx_j(x')
  & = 
  \frac{1}{\tend}
  \E\Biggl[\int_0^\tend
  \sumall{j}\sum_{k=1}^3
  \diffumd_{j,k}(x;\mdpos{\no}) \diffumd_{j,k}(x';\mdpos{\no})
  \,dt
  \;\Bigl\lvert\;
  \mdpos[0]{\no}=\mdpos[\no]{0}
  \Biggr],
  && x,x'\in D,
\end{align}
\end{subequations}
for some fixed, deterministic, initial conditions 
$\mdpos[0]{\no}=\mdpos[\no]{0}$. The initial condition for the
coarse-grained phase-field is
$\pfcgt[0]=\pfen(\cdot;\mdpos[\no]{0})$. 
This particular coarse-graining is motivated by the argument that the
coarse-grained model will be used to compute properties on the form 
$\E\left[\smooth(\pfen(\cdot;\mdpos[\tend]{\no}))\right]$, where 
$\smooth\,\colon\,D\to\R$ is a smooth function and $\tend>0$ is a
fixed final time.
The optimal coarse-grained model is the one that minimises the error
in the expected value; 
using the conditional expected values 
$\dual(\mu,t)=\E[\smooth(\pfcgt[\tend])\,|\,\pfcgt=\mu]$,
this error can be expressed as
\begin{multline*}
  \E\left[\smooth(\pfen(\cdot;\mdpos[\tend]{\no})) \right]
  - \E\left[\smooth(\pfcgt[\tend]) \right]
  \\
  \shoveleft{
  = 
  \E\Biggl[
    \int_0^\tend 
    \Big\langle \dual'(\pfen(\cdot;\mdpos{\no}),t)\; , \;
      \driftmd(\cdot;\mdpos{\no})-\driftcgx(\cdot)
    \Big\rangle_{L^2(D)} \;dt }
    \\
    + 
    \frac{1}{2}
    \int_0^\tend
    \bigg\langle \dual''(\pfen(\cdot;\mdpos{\no}),t)\; , \;
      \sumall{j}\sum_{k=1}^3(\diffumd_{j,k} \otimes \diffumd_{j,k})
      (\cdot,\cdot;\mdpos{\no}) - 
      \sumall[\nrw]{j}(\diffucgx_j \otimes \diffucgx_j)(\cdot,\cdot)
    \bigg\rangle_{L^2(D\times D)} \;dt
  \Biggr],
\end{multline*}
where  $\otimes$ denotes the tensor product 
$(\diffucgx_j\otimes\diffucgx_j)(x,x')=\diffucgx_j(x)\diffucgx_j(x')$,
and $\dual'$ and $\dual''$ denote the first and second variations of
$\dual(\mu,t)$ with respect to $\mu$. 
Assuming that $\dual'$ can be expanded in powers of
$\driftmd-\driftcgx$, the choice~\eqref{eq:def_driftcgx} cancels the
leading term in the error associated with $\dual'$. Similarly, 
\eqref{eq:def_diffucgx} corresponds to cancelling the dominating term
in the expansion of $\dual''$. 

In a practical computation the functions \driftmd\ and $\diffumd_j$ 
can only be evaluated in a discrete set of points
$D_K=\{x^1,\ldots,x^K\}\subset D$. The 
right hand sides in~\eqref{eq:def_driftcgx}
and~\eqref{eq:def_diffucgx} become a vector and a symmetric positive
semidefinite $K$-by-$K$ matrix, respectively. 
Hence $\driftcgx(x)$ becomes a vector of tabulated values for 
$x\in D_K$.
It is natural to have one Wiener process per point $x^k$
in the spatial discretisation, so that $K=\nrw$. 
The corresponding $K$ tabulated individual diffusion coefficient 
functions, $\diffucgx_j$, will be obtained by a square root
factorisation of the computed matrix, by means of an eigenvector
expansion; 
this choice of factorisation preserves the connection between
the evaluation point $x_k$ and the elements $k$ in $\diffucgx_j$ and 
produces spatially localised functions, consistent with the
association of individual Wiener processes and points in $D_K$.

In the second step, the initial configuration, $\mdpos[\no]{0}$, 
in~\eqref{eq:def_coeffcgx} is chosen so that the microscopic domain
$D$ includes a solid--liquid interface in equilibrium.
Since the interface is stationary no phase transformation occurs in
the simulation, and consequently the part of the reaction term in the
Allen-Cahn equation~\eqref{eq:AllenCahn} relating the speed of the
phase change to the deviation from the melting point,
$k_2k_3\monoton'(\pfgen) (\Tmelt-\temp)$,
can not be obtained; the simulation must be performed 
at the melting point, \Tmelt, under the given conditions.
The simulation of a travelling front, off the equilibrium temperature,
requires more advanced micro model simulations than the ones considered
here. 

The interface is assumed to be locally planar on the microscopic
scale and the spatially averaged properties are expected to vary
much more slowly in the directions parallel to the interface than in
the direction normal to the interface. 
Label the direction normal to
the interface as direction $x_1$ and let $x_2$, $x_3$ be orthogonal
directions in the plane of the interface.
Then the mollifier, \molli, in~\eqref{eq:def_pfen} can be chosen to
make the averages much 
more localised in the $x_1$ direction than in the $x_2$ and $x_3$
directions. In the microscopic domain, $D$, the averages in the $x_2$
and $x_3$ directions are chosen to be uniform averages over the entire 
domain, so that the phase-fields, \pfen\ and \pfcg, and the drift and
diffusion functions, \driftmd, $\diffumd_{j,k}$, \driftcgx, and
$\diffucgx_j$, become functions of one space variable, $x_1$. 
Hence the evaluation points in $D_K$ are only distinguished by their $x_1$
coordinates. 
As mentioned above, the drift coefficient, \driftmd, depends on the
derivatives up to second order of, \molli, and the potential energy
contributions $\pfen_i$. After averaging out the $x_2$ and $x_3$
dependence, it can be written as 
\begin{align*}
  \driftmd(x_1;\mdpos{\no}) & = 
  \kb\temp\frac{\partial^2}{\partial x_1^2}\pfen(x_1;\mdpos{\no})
  + \frac{\partial}{\partial x_1} A_1(x_1;\mdpos{\no})
  + A_0(x_1;\mdpos{\no}),
\end{align*}
for some functions $A_1$ and $A_0$. Keeping this form in the
averaging, the coarse-grained drift coefficient
in~\eqref{eq:def_driftcgx} can be written
\begin{align*}
  \driftcgx(x_1) & = 
  \kb\temp\frac{\partial^2}{\partial x_1^2}\pfen_\mathrm{av}(x_1)
  + \frac{\partial}{\partial x_1} \driftcgx_1(x_1)
  + \driftcgx_0(x_1),
\end{align*}
where the second order derivative of the averaged phase-field,
\begin{align}
  \label{eq:def_mdphaseav}
  \pfen_\mathrm{av}(x_1) &=
  \frac{1}{\tend}
  \E\biggl[\int_0^\tend\pfen(x_1;\mdpos{\no})\,dt\biggr],
\end{align}
corresponds to the diffusion term in~\eqref{eq:AllenCahn}.
Assuming that the averaged phase-field $\pfen_\mathrm{av}$
is a monotone function of $x_1$ in the interface, the explicit
dependence on the spatial variable can be eliminated by inverting
$\pfen_\mathrm{av}$ and defining
\begin{align}
  \label{eq:def_coeff_cg}
  \driftcg(\pfcg) & = \driftcgx(\pfen_\mathrm{av}^{-1}(\pfcg)),
  &
  \diffucg_j(\pfcg) & = \diffucgx_j(\pfen_\mathrm{av}^{-1}(\pfcg)),
\end{align}
which give drift and diffusion coefficients on the
form~\eqref{eq:sde_pfcg}. 

The present study is a practical test of the method described above. 
In particular the aims are to verify that Smoluchowski dynamics can be
used in practise, in the sense that the coarse grained drift and
diffusion coefficient functions can be determined together with the
phase-field model potential, \dwell, and that they seem reasonable.
For this purpose simulations are performed at just one temperature and
density (at the melting point) and with just two values of the angle
of the stationary interface with respect to the crystal structure in
the solid. 
An actual determination of the model functions in the phase field
model would require many more simulations with varying parameters.

\section{Computational Methods}
\label{sec:methods}

The numerical computations consist of molecular dynamics
computations, giving the microscopic description of the two-phase
system, and the extraction of model functions for a coarse
grained stochastic differential equation model. 

\subsection{Molecular Dynamics Models and Simulation}
\label{sec:models}

Two mathematical models of the material are used; both are
one component 
molecular dynamics models where the interaction between 
particles is determined by a pair potential of the exponential-6
(Exp-6) type.
The coarse graining is based on a stochastic model where the particle
trajectories on the diffusion time scale are given by the
Smoluchowski dynamics~\eqref{eq:Smoluchowski}.
The computations with this model are performed under constant
volume at the melting point where a liquid and a solid phase
coexist in the computational domain. 
The melting point is determined using constant pressure
simulations of the deterministic molecular dynamics model where
the particle trajectories are determined by Newton's second law
with forces given the by gradients of the model potential.
Both models and the corresponding simulations are described
below, after a description of the potential common to the
models.

\subsubsection{Pair Potential Defining the Total Potential
  Energy}

\begin{figure}[hbp]
  \centering
  \subfigure[Pair Potentials for Argon]{
    \label{fig:Argon_pots}
    \includegraphics[height=6.0cm]{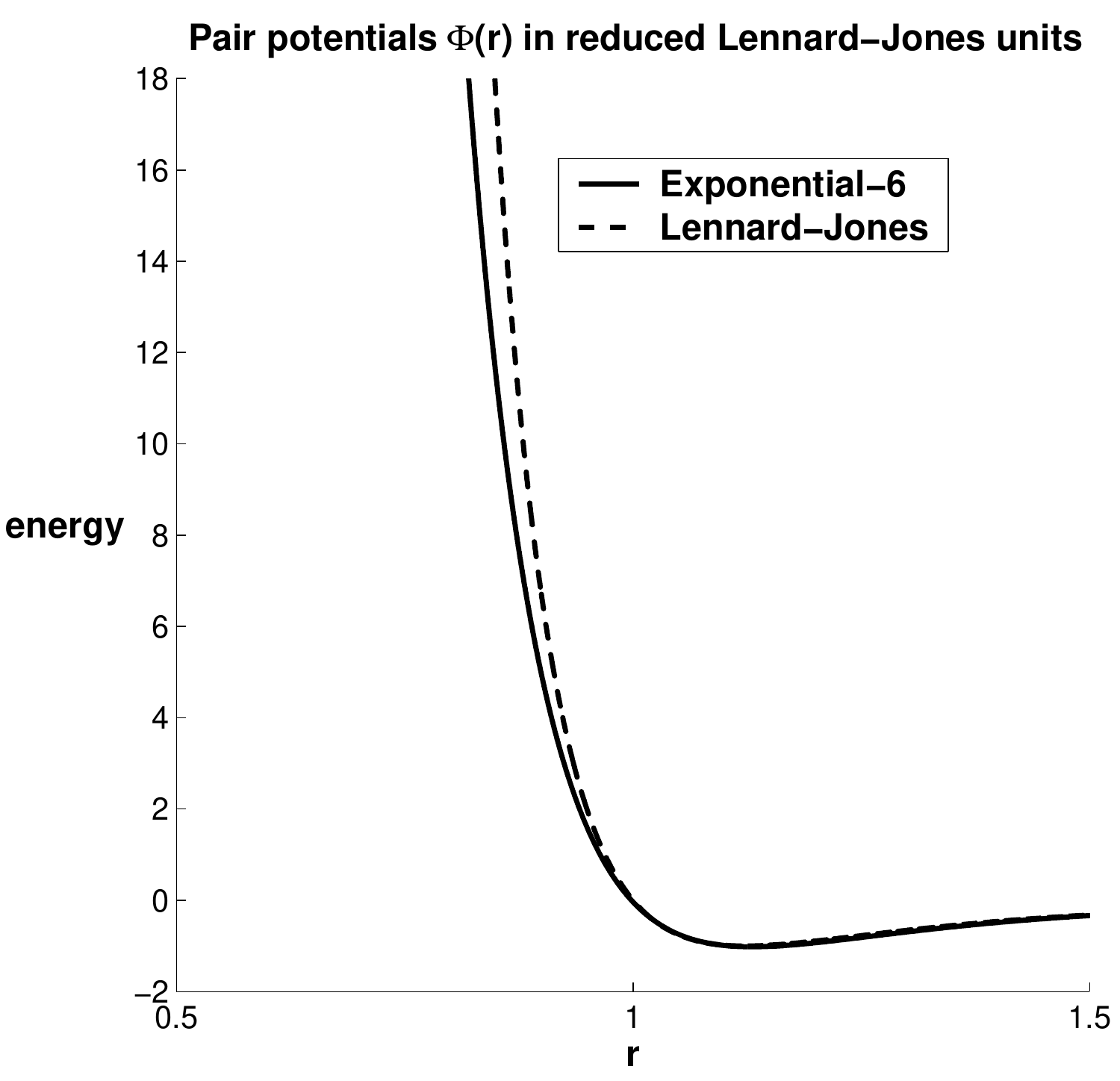}
  }    
  \hspace{1cm}
  \subfigure[General form of the Lennard-Jones potential]{
    \label{fig:LJ_general}
    \includegraphics[height=6.0cm]{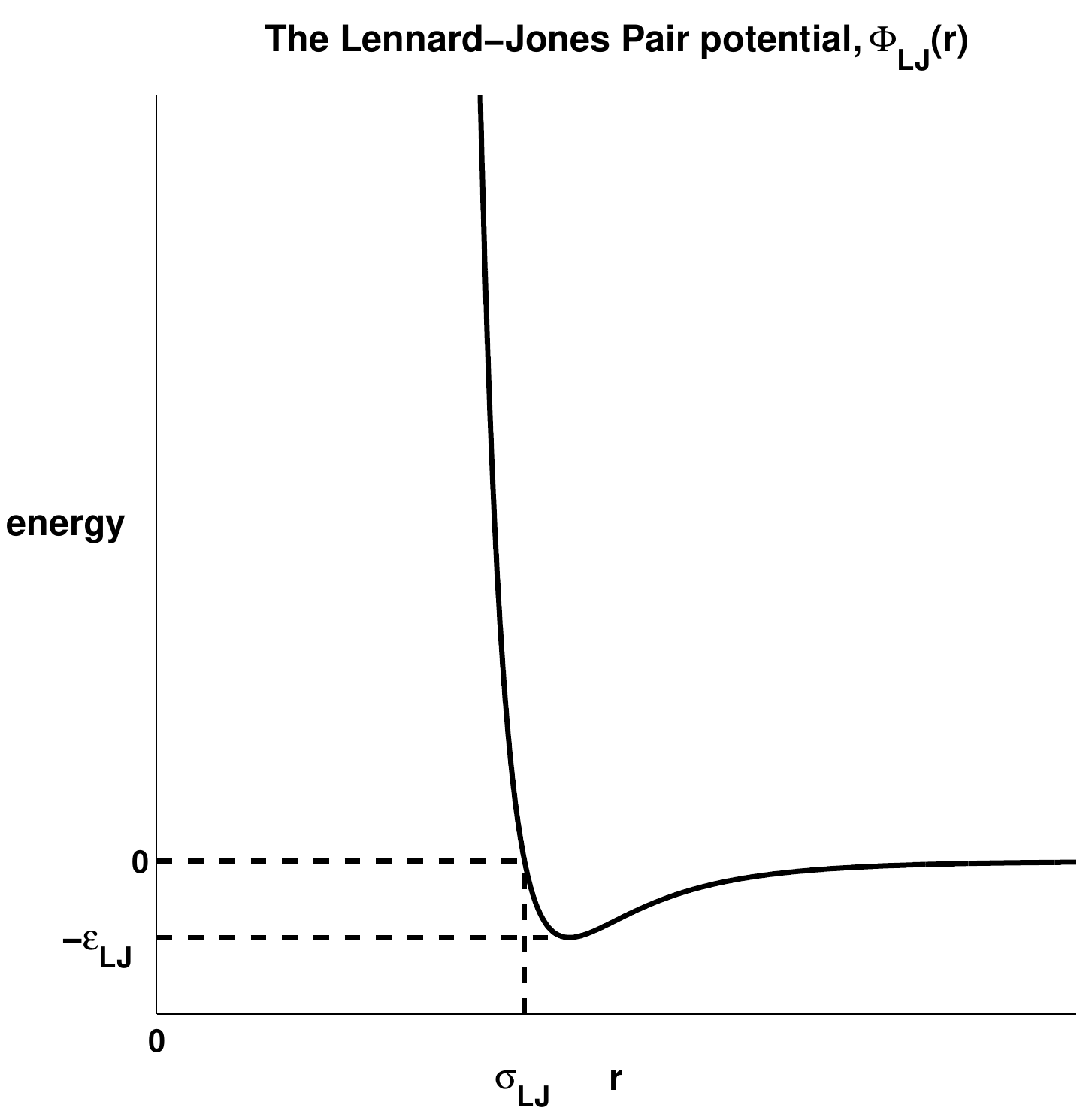}
  }
  \caption{(a): The Exp-6 pair potential is similar to
    the Lennard-Jones pair potential near the minimum, but the
    repulsion is slightly weaker in the Exp-6. The radius and the
    energy are measured in reduced Lennard-Jones units, where the
    Lennard-Jones parameters are 
    $\epsilon_{\mathrm{LJ}}=\kb 120\,\mathrm{K}$ and 
    $\sigma_{\mathrm{LJ}}=3.405\,$Å. 
    \newline
    (b): The parameter $\sigma_{\mathrm{LJ}}$ is the radius where
    the Lennard-Jones potential is 0, which is equal to the
    potential at infinite separation, and the parameter
    $\epsilon_{\mathrm{LJ}}$ is the depth of potential minimum.
  }
  \label{fig:pots}
\end{figure}

The microscopic system consists of \nrp\ identical particles at
positions $\mdposgen=(\mdpos[\no]{1},\ldots,\mdpos[\no]{\nrp})$
in three dimensions. 
The total potential energy, \totpot, of the system is determined
by the particle positions through
\begin{align}
  \label{eq:def_poten}
  \totpot(\mdposgen) 
  & = \frac{1}{2}\sumall{i} \sumneq{k}{i} 
  \pairpot(\mdpos[\no]{i}-\mdpos[\no]{k}),
\end{align}
using pairwise interactions only.
The pair potential is the spherically symmetric Exp-6 potential
\begin{align}
  \label{eq:def_exp6}
  \pairpot(r) & = 
  A\exp(-Br)-\frac{C}{r^6},
\end{align}
with $r$ denoting the distance between two particles, and $A$,
$B$, and $C$ being positive model parameters. 
The Exp-6 potential, like the similar Lennard-Jones pair
potential,  
$\pairpot_\mathrm{LJ}(r)=
4\epsilon_\mathrm{LJ}\left((\sigma_\mathrm{LJ}/r)^{12}
-(\sigma_\mathrm{LJ}/r)^6\right)$,
is a short range interaction that
can be used to model condensed noble gases. 
With the parameters used here, obtained from~\cite{exp6pars}, the
Exp-6 potential models Argon at high pressures. 
At pressures around 2~GPa, where the solid-liquid phase
transition will be simulated, the Exp-6 potential with
its slightly softer repulsive part describes the equation of
state of Argon better than the Lennard--Jones
potential does; see~\cite{exp6pars,exp6pars_II}. 
The shapes of the two pair potentials around the global minimum
of the Lennard--Jones potential can be compared in
Figure~\ref{fig:Argon_pots}; the typical inter atomic distances
between nearest neighbours in both the simulated solid and liquid
will be close to 1.
Note that, while the Lennard--Jones pair potential tends to
infinity as the interatomic distance tends to zero,
the Exp-6 pair potential, as stated in~\eqref{eq:def_exp6},
reaches a global maximum before turning down and approaching minus
infinity in the limit. This clearly illustrates that the model
based on the Exp-6 potential breaks down if two atoms come too
close, but neither one of the pair potentials is designed to
describe interactions of particles much closer than the typical
nearest neighbour separation.

For short range potentials, like the Exp-6 and the Lennard-Jones
potentials, the potential (and its derivative) decay sufficiently
fast for the combined effect on the total potential energy (and
the interatomic forces) of all atom pairs separated more than a
certain distance to be negligible compared to the 
effect of the pairs separated less than the same distance. 
To take advantage of this in computations a
cut-off radius is introduced and 
all interactions between particles separated by a distance larger
than the cut-off are neglected;
instead of summing over all $k\!\!\neq\!\!i$ in the inner sum
in~\eqref{eq:def_poten} the sum is only taken over particles in a
spherical neighbourhood of particle $i$.

All the physical quantities in this report are given in the
reduced Lennard-Jones units. Thus length is measured in units of
$\sigma_\mathrm{LJ}$, energy in units of $\epsilon_\mathrm{LJ}$,
and time in units of 
$\sqrt{m\sigma_\mathrm{LJ}^2/\epsilon_\mathrm{LJ}}$, where $m$ is
the mass of one atom. (The time unit is the inverse of the
characteristic frequency.)
A list of the dimensionless units in the Argon model as
well as the parameters in the Exp-6~potential can be found in  
Table~\ref{tab:LJunits}. 
At the temperatures and pressures considered here, the stable
phase of the Exp-6 potential is either the Face Centered Cubic (FCC)
lattice or a liquid phase.
\begin{table}[htp]
  \centering
  \label{tab:LJunits}
  \begin{tabular}{lr}
    \begin{tabular}{|l|l|}
      \hline
      Quantity & Unit \\
      \hline
      Energy & $1.6568\cdot10^{-21}$~J \\
      Time & $2.1557\cdot10^{-12}$~s \\
      Mass & $6.6412\cdot10^{-26}$~kg \\
      Length & $3.405\cdot10^{-10}$~m \\
      Temperature &  $120$~K \\
      Pressure & $4.1968\cdot10^{7}$~Pa \\
      \hline
    \end{tabular}
    & 
    \begin{tabular}{|l|l|}
      \hline
      Constant & Value \\
      \hline
      \kb & 1 \\
      \hline
      \multicolumn{2}{c}{~} \\
      \hline 
      Parameter & Value \\
      \hline
      $A$ & $3.84661\cdot 10^5$ \\
      $B$ & $11.4974$ \\
      $C$ & $3.9445$ \\
      \hline
    \end{tabular}
  \end{tabular}
  \caption{Atomic units and corresponding values of physical
    constants and parameters in the Exp-6 model~\eqref{eq:def_exp6}.
    Non dimensional molecular dynamics equations are
    obtained after normalising with the atom mass, $m$, and the
    Lennard-Jones parameters, $\sigma_\mathrm{LJ}$ and 
    $\epsilon_\mathrm{LJ}$; in this Argon model
    $m=6.6412\cdot10^{-26}$~kg (or 39.948 atomic mass units), 
    $\sigma_\mathrm{LJ}=3.405$~Å, and 
    $\epsilon_\mathrm{LJ}/\kb=120$~K, where \kb\ is the Boltzmann
    constant.} 
\end{table}

\subsubsection{Newtonian System Simulated at Constant Pressure}
\label{sec:NTP}

The purpose here is to approximately determine the melting point
at a high fixed pressure, to be able to set up and simulate
stationary (FCC-liquid) two-phase systems later. Determination of the
melting point follows the two-phase method described by
Belonoshko and co-authors in~\cite{Anatoly}.

\paragraph{The mathematical model}
is a classical system of
\nrp\ identical particles where the positions, 
$\mdpos{\no}=(\mdpos{1},\ldots,\mdpos{\nrp})$,
and the velocities, 
$\mdvel{\no}=(\mdvel{1},\ldots,\mdvel{\nrp})$,
evolve in time according to Newton's equations
\begin{subequations}
  \label{eq:Hamiltonian}
\begin{align}
  \label{eq:Hamiltonian_pos}
  \frac{d\mdpos{\no}}{dt}
  & = \mdvel{\no}, \\
  \label{eq:Hamiltonian_vel}
  \frac{d\mdvel{\no}}{dt}
  & = -\nabla_{\mdposgen} \totpot(\mdpos{\no}),
\end{align}
\end{subequations}
where the total potential energy of the system is given
by~\eqref{eq:def_poten}-\eqref{eq:def_exp6} using the parameter
values in Table~\ref{tab:LJunits}. 
Here $\nabla_{\mdposgen}$ denotes the gradient with respect to
the particle positions. The force acting on particle $i$ is 
$-\nabla_{\mdpos[\no]{i}}\totpot(\mdpos{\no})$ and, since all
particles have unit mass in the non-dimensional units, the
acceleration is equal to the force.
Particle positions are restricted to a finite computational
box with periodic boundary conditions, corresponding to an
infinite system where the same configuration of particles is
repeated periodically in all three directions; 
a particle leaving the computational cell on one side enters the
cell again from the opposite side and particles interact with
periodic images of particles in the cell.

For a fixed volume of the computational cell the
equations~\eqref{eq:Hamiltonian} 
will preserve the total energy, \toten, (the sum of
potential and kinetic energy) of the system as well as the number
of particles. It will approximately sample the $(N,V,E)$ ensemble.
In the determination of the melting point the simulations are
instead performed in an approximation of the $(N,T,P)$ ensemble,
using a constant 
number of particles, \nrp, a constant temperature, \temp, and a 
constant pressure, \pressure. This must allow for the volume of the
computational cell to change during the simulation. There must
also be mechanisms for keeping the temperature and the pressure
constant, thus 
modifying~\eqref{eq:Hamiltonian} 
so that the total energy varies.

\paragraph{Numerical computations}
of the (N,T,P) molecular dynamic simulations were
performed using Keith Refson's publicly available software package
Moldy,~\cite{Moldy}. 
Constant temperature was enforced using the Nos\'{e}-Hoover
thermostat, where the equations of motions~\eqref{eq:Hamiltonian}
are modified, and extended, to include an additional degree of
freedom modelling a thermal reservoir. The fictitious inertia
associated with the thermal reservoir was
$100~\mathrm{kJ}\,\mathrm{mol}^{-1}\,\mathrm{ps}^2$, 
corresponding to $21.57$ in the dimensionless equation.
The pressure was kept constant using the Parinello-Rahman
equation, controlling the dynamics of the vectors (three edges) 
that define the computational cell. The fictitious mass parameter
in the Parinello-Rahman equation was $300\,\mathrm{amu}$ 
corresponding to $1.20\cdot10^4$ in the reduced Lennard--Jones
units. 
A short description of the Nos\'{e}-Hoover thermostat and the
Parinello-Rahman equation, with references to papers with
theoretical foundations of the methods, can be found in the
manual~\cite{Moldymanual}. 

The time stepping method in Moldy is a modification of Beeman's
algorithm using predictor-corrector iterations in the computation
of the velocities; see~\cite{Moldymanual} for details. The
simulations described here used the constant time step 
$4.639\cdot10^{-5}$ and the potential cut-off $2.937$.

\paragraph{In the two-phase method} for determination of the 
melting point the molecular dynamics simulation starts from an 
initial configuration that is part solid and part liquid. As the
$(\nrp,\temp,\pressure)$ simulation proceeds the whole liquid
part will solidify, if $\temp<\Tmelt$ for the given pressure, or
the solid will melt, if $\temp>\Tmelt$, resulting in a single
phase system. Starting from a coarse estimate of the temperature
interval containing the melting temperature, that interval can be 
narrowed down by running simulations at temperatures in the
interval and noting whether they equilibrate to an all solid or
an all liquid system. The validity of this two-phase approach has
been verified in~\cite{Anatoly} for determining, among other
things, the melting point of a molecular dynamics model of Xenon,
similar to the Argon model used here.

The initial configuration in a two-phase simulation was
composed of pre-simulated solid and liquid configurations. 
The solid part was prepared by taking a perfect FCC configuration and
performing a short molecular dynamics run at the 
temperature and pressure of the intended two-phase simulation 
to adapt the size of the computational cell. 
\begin{figure}[htbp]
  \centering
  \subfigure[The FCC unit cube]{
    \includegraphics[width=1.80in]{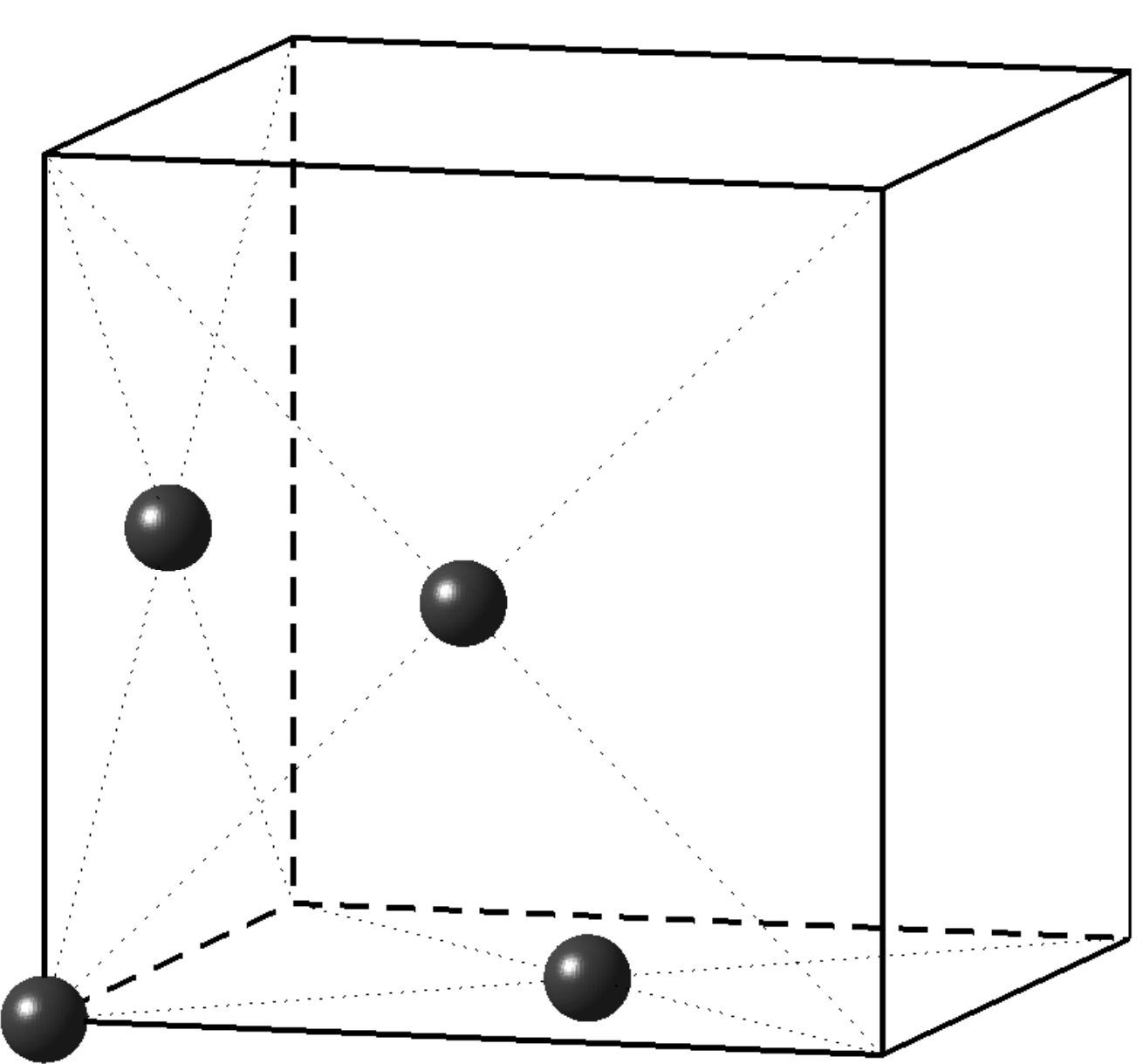}
    \label{fig:fcc_unit}
  }
  \subfigure[Eight FCC unit cubes]{
    \includegraphics[width=1.80in]{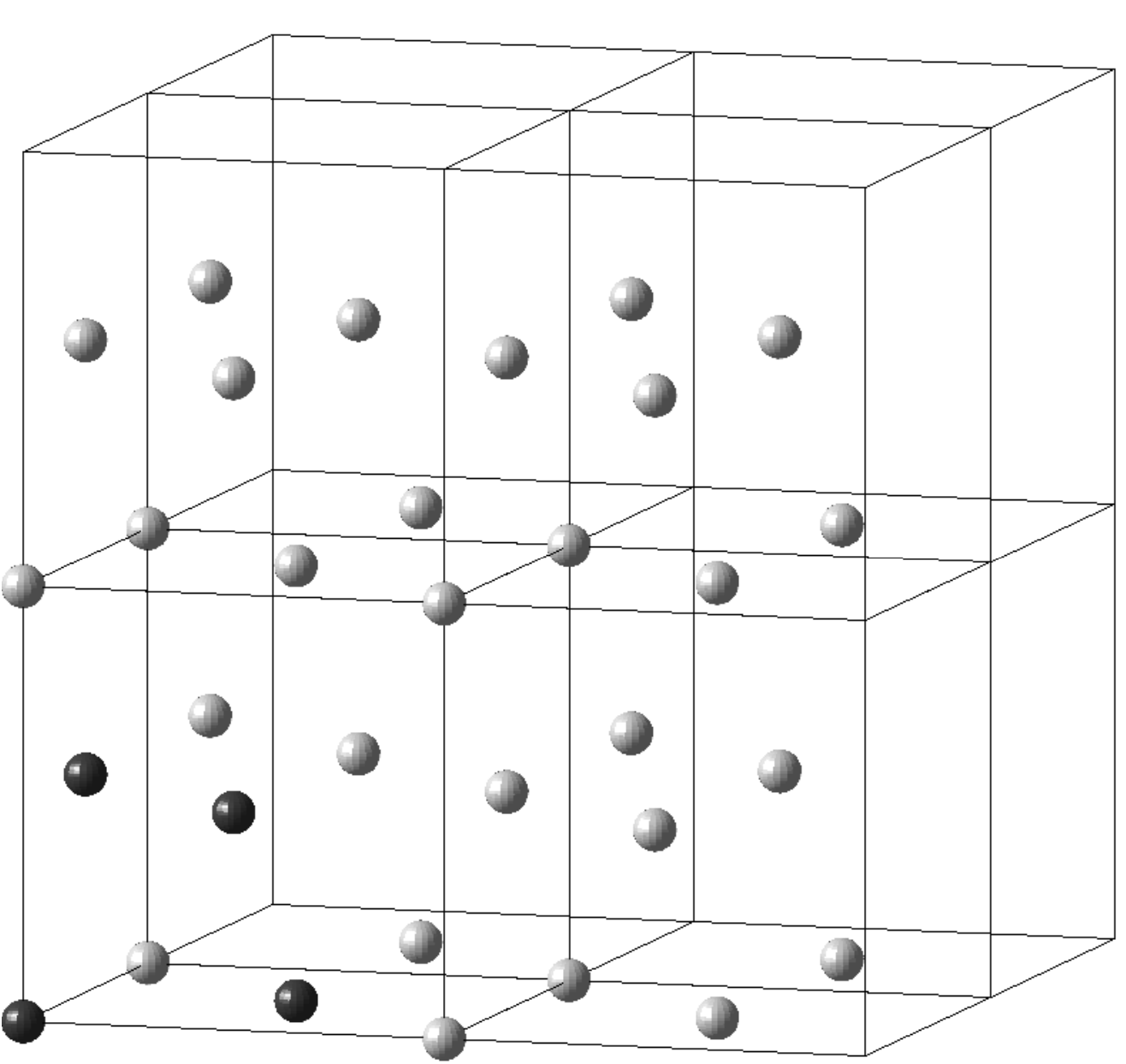}
    \label{fig:fcc_cubes}
  }
  \caption{A perfect FCC lattice consists of FCC unit cubes,~(a),
    stacked next to each other in three dimensions,~(b). With one atom
    in the $(0,0,0)$ corner of the unit cube the three other atoms
    are placed at the centres of the cubic faces
    intersecting in $(0,0,0)$.}
  \label{fig:lattice}
\end{figure}
Initially the sides of the computational cell were aligned with
the sides of the unit cube in the perfect FCC lattice; see
Figure~\ref{fig:lattice}. While in 
general the dynamics of the cell edges in the Parinello-Rahman
equations allow the cell to take the shape of any
parallelepiped, here the dynamics were restricted to only allow
rescaling, without rotation, of the three edges and thus keeping
the rectangular box shape of the cell.
The preparation of the liquid part started from the configuration
of the already prepared FCC-solid and a run was performed at a
temperature well over the estimated melting point, where the
sample would melt quickly; after equilibrating at the higher
temperature the sample was quenched to the temperature of the
two-phase simulation.
Only one side of the computational cell was allowed to change
while preparing the liquid part and thus the orthogonal cross
section of the simulation cell was preserved from the FCC
simulation. 
The solid and liquid parts were joined in the two-phase initial
configuration by placing them next to each other, letting the
cell faces of identical shape face each other. 
The general appearance is similar to the configurations shown in 
Figure~\ref{fig:init_confs} on page~\pageref{fig:init_confs},
even though those configurations belong to the constant
volume Smoluchowski simulations where the set up procedure is
slightly modified. 
Periodic boundary conditions were still applied in all
directions, so that each part (solid or liquid) corresponded to a
semi-infinite slab surrounded on two sides by the other phase
with the effect of simulating a periodic, sandwiched, material.
Voids of thickness of approximately one nearest neighbour
separation were introduced in both solid--liquid interfaces to
make sure that no pair of particles ended up to close in the
initial configuration. 
Since the two-phase simulations were performed at constant
pressure, the voids would fill in the beginning of the run as
the length of the computational cell decreased.

In the two-phase simulations the lengths of all three vectors
defining the cell edges were allowed to change. Starting from an
initial two-phase configuration the molecular dynamics simulation
was run until the system was considered equilibrated. After
equilibration the computational cell was filled with either the
solid or the liquid phase. 
The density of the FCC solid is higher than that of the liquid
phase. If the phase change was solidification of the liquid, then
the volume of the computational cell would decrease during the
equilibration stage before assuming an approximately constant
value; if the solid was melting, the total volume would grow
during equilibration. The density of the stable phase at the
given pressure and temperature was obtained by time averages of
the simulation after equilibration.

\begin{figure}[htp]
  \centering
  \subfigure[Simulation data]
  {
    \label{graph:melt_point}
    \includegraphics[width=7.9cm]{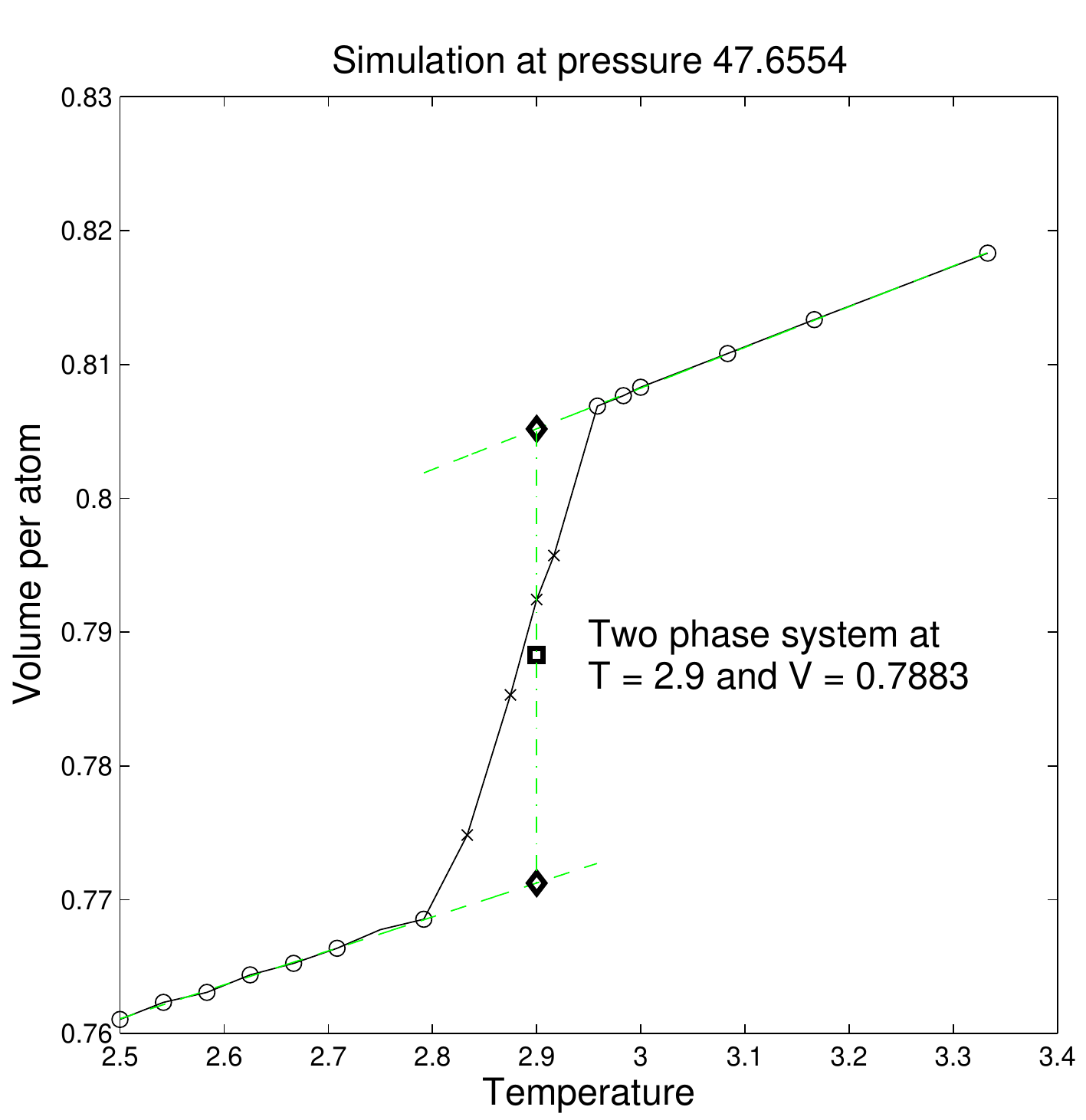}
  }
  \hspace{0.5cm}
  \subfigure[Extrapolated data at $\temp=2.9$]
  {
    \label{tab:melt_point}
    \begin{picture}(136,220)(0,0)
      \put(66,42){\makebox(0,0){
          \begin{tabular}{|l|c|c|}
            \hline
            & $V_a$ & \density \\
            \hline
            Liquid & 0.8060 & 1.241 \\ 
            \hline
            FCC & 0.7714 & 1.296 \\ 
            \hline
            Combined & 0.7883 & 1.269 \\ 
            \hline
          \end{tabular}}}
    \end{picture}
  }
  \caption{(a) Volume per atom as a function of temperature at a
    the pressure 47.6554 (or $2.0$~GPa) in $(N,P,T)$ simulations. 
    Data points from simulations that are considered equilibrated
    are marked with $\circ$ and those from simulations that are
    not equilibrated are marked with $\times$. 
    The two regions of equilibrated values where the volume per
    atom varies approximately linearly correspond to solid (FCC),
    at lower temperatures, and liquid, at higher temperature,
    respectively. The melting point at the given pressure is 
    somewhere in between; the approximate value $\temp=2.9$ is used
    below and in the constant volume simulations.
    \newline
    (b) The volume per atom of solid and liquid have been
    extrapolated to $\temp=2.9$ by least square fits of 
    straight lines to the simulation data and the corresponding
    number densities, \density, have been computed. 
    If $\temp=2.9$ is sufficiently close to the melting point at
    this pressure, then the two phases will coexist in
    \emph{constant volume}, $(N,V,T)$, simulations provided that
    the total density is between the estimated densities of pure
    solid and pure liquid. The ratio of the volumes of the solid
    and the liquid part is determined by the total density of the
    combined system. The tabulated value of the density for a
    combined system gives approximately equal volumes of both
    parts at a pressure close to the one in the constant pressure
    simulations.} 
  \label{fig:melt_point}
\end{figure}

When the volume per particle is shown as a function of the
temperature, at constant pressure, it will display a sharp change
at the melting point; see Figure~\ref{graph:melt_point} on
page~\pageref{fig:melt_point}. 
The procedure will obtain an interval around the melting point
and the accuracy can be improved by performing simulations at
more temperatures to shorten the interval of uncertainty. 
However, the equilibration requires longer time when close to the
melting point and the cost for refining the approximation grows, not
only because the number of simulations grows, but more importantly
because every single simulation takes longer to perform.

The main purpose here is to investigate the possibility of obtaining
the model functions in a coarse grained phase-field model from
$(\nrp,\vol,\temp)$ Smoluchowski dynamics simulations, as
described next. Therefor the accuracy in the determination of the
melting point at the given pressure is critical only to the
extent that it must be possible to perform the constant volume
simulations at this temperature; that is, it must be possible to
perform simulations on a two-phase system with stable interfaces
between the solid and liquid parts. If the purpose were to
perform computations at the melting point at this very pressure,
then more computational effort would have to be spent on the
accuracy of the melting point and the corresponding densities. 

The numerical simulations were performed with $\nrp=8000$
particles; the initial solid configuration consisted of 4000
particles, corresponding to $10\!\times\!10\!\times\!10$ FCC unit
cells with four atoms each, and the liquid had the same number of
particles.
From simulations at the pressure 47.7 in the reduced
Lennard-Jones units (corresponding to $2.0$~GPa) 
an approximate value of 2.9 for the melting point was obtained
together with number densities for the liquid and solid 
extrapolated to this temperature; see Figure~\ref{fig:melt_point}
on page~\pageref{fig:melt_point}. 
Fixing the temperature and the number density $\nrp/\vol$,
only one degree of freedom remains in the triple
$(\nrp,\vol,\temp)$, allowing the system size to vary.

\subsubsection{Smoluchowski System Simulated at Constant Volume}
\label{sec:NTV}

The constant volume and temperature Smoluchowski dynamics
two-phase simulations described here were used to compute 
the functions~\eqref{eq:def_coeff_cg} defining the coarse-grained
phase-field dynamics~\eqref{eq:sde_pfcg}, 
as described in the introduction. 
This meant computing time averaged quantities 
like the time averaged potential energy
phase-field~\eqref{eq:def_mdphaseav} and the corresponding
coarse-grained drift and diffusion coefficient
functions~\eqref{eq:def_coeffcgx}.

\paragraph{The mathematical model} is that of \nrp\ particles
whose positions $\mdpos{\no}$ follow the Smoluchowski dynamics
\begin{align}
  \label{eq:Smol_impl}
  d\mdpos{\no} & = -\nabla_{\mdposgen} \totpot(\mdpos{\no})\;dt 
  + \sqrt{2\kb \temp}\;dW^t,
\end{align}
introduced on page~\pageref{eq:Smoluchowski}. 
There are no velocities in the Smoluchowski dynamics. Instead the
positions of all particles in the system give a complete
description of the system at a particular time. Such a
description, $\mdpos{\no}$, will be refered to as a configuration
of the system.
The particles are contained in a computational cell, shaped like
a rectangular box, of fixed dimensions and the boundary
conditions are periodic in all directions. Hence the volume,
\vol, and the number of particles, \nrp, are fixed.
Without velocities there is no kinetic energy, but the
temperature, \temp, enters directly in the dynamics. The
temperature parameter is held fixed, which can be viewed as a
kind of thermostat built into the dynamics. 

Since the volume of the computational cell is constant, 
unlike in the $(\nrp,\temp,\pressure)$ simulations above, 
the overall density of the system remains constant over time,
which allows for stationary two-phase configurations where part
of the domain is solid and part is liquid.

\paragraph{The numerical simulations} The discrete time
approximations  
\mdposd{\no} of \mdpos[t_n]{\no}, were computed using the 
explicit Euler-Maruyama scheme
\begin{align}
  \label{eq:algo_Smol}
  \mdposd{\no} & = \mdposd[n-1]{\no}
  -\nabla_{\mdposgen} \totpot(\mdposd[n-1]{\no})\;\Delta t^n
  + \sqrt{2\kb \temp}\;\Delta W^n,
\end{align}
where $\Delta t^n = t^n-t^{n-1}$ is a time increment
and $\Delta W^n = W(t^n) - W(t^{n-1})$ is an increment in the
$3\nrp$-dimensional Wiener process.
Each run was performed using constant time step size, 
$\Delta t^n\equiv\Delta t$, but the time step could change
between different runs depending on the purpose; in the
equilibration phase the typical step size was $\Delta t=10^{-4}$,
but in the production phase the step size had to be taken
smaller, as discussed later.

The computation of 
$\nabla_{\mdposgen}\totpot(\mdposd[n-1]{\no})$ in every time
step is potentially an $\mathcal{O}(\nrp^2)$ operation since the
potential is defined by pairwise interactions.
The computations described here used the potential cut-off radius
$3.0$, which meant that each particle only interacted directly
with a relatively small number of neighbours (independent of
\nrp\ since the density was approximately constant). 
To avoid the $\mathcal{O}(\nrp^2)$ task of computing all pairwise
distances in each time step, the computational cell is divided
into smaller sub cells, where the size is defined in terms of the
cut-off radius so that two particles only can interact if they
are in the same sub cell or in two neighbouring sub cells; 
information about particles migrating between sub cells is
exchanged in each time step.
The computations use a two dimensional grid of sub cells, where
the particle positions within each sub cell are sorted with
respect to the third coordinate dimension in every time step.
When the particles are sorted the sweep over all particles in a
sub cell can be efficiently implemented and the sorting procedure
is not too expensive since the particles do not move far in one
time step.
A more thorough description of this algorithm can be found
in~\cite{MD_algo}. 
The actual code used here is a
modification of a parallelised code for Newtonian molecular
dynamics obtained from Måns Elenius in Dzugutov's
group\cite{Dzugutov}; the main modifications
when adapting to Smoluchowski dynamics is the removal of
velocities from the system and the introduction of a pseudo
random number generator for the Brownian increments, $\Delta W^n$.

With the cut-off radius $3.0$ used in the computation and the
model parameters in Table~\ref{tab:LJunits} on
page~\pageref{tab:LJunits}, the Exp-6 pair potential and its
derivatives are small at the cut-off radius. 
Still the potential will be discontinuous at the cut-off, unless
it is slightly modified. A small linear term is added to make the
potential continuously differentiable at the cut-off radius. 
In the practical computations, both the pair potential and the 
derivatives were obtained by linear interpolation from tabulated
values. 

The random number generator for normally distributed random
variables was the Ziggurat method, described in~\cite{ziggurat},
in a Fortran~90 implementation by Alan Miller, accessible from
Netlib~\cite{Netlib}. The underlying
32-bit integer pseudo random number generator is the 3-shift
register SHR3. 
Since the purpose of the simulations only is to
investigate if the coarse-graining procedure gives reasonable
results just one pseudo random number generator was used, while 
several different random number generators ought to be used in a
practical application. 
The generator was initialised with different seeds on different
processors in the parallel computations, but it does not have 
distinct cycles simulating independent random variables. 
The hope is that the nature of the molecular dynamics simulations
is enough to avoid the danger of correlated random numbers on the
different processors, but this could be tested by comparing with
other pseudo random generators that actually simulate independent
random variables on different processors.

\paragraph{The two-phase systems}
for the Smoluchowski dynamics simulations were set up to obtain a
two-phase system at temperature $\temp=2.90$ with approximately
equal volumes of solid and liquid and with stationary interfaces.
To achieve this two equal volumes of FCC-solid and liquid were
pre-simulated with the densities tabulated in
Figure~\ref{fig:melt_point}, on page~\pageref{fig:melt_point}. 
The preparation of the initial configurations for the
Smoluchowski dynamics two-phase simulations was similar to the
procedure described above, but some adjustments must be made
because of the constant volume restriction.
The shape of the computational cell used when generating the
solid part was chosen to match the periodic structure of the FCC
lattice at the tabulated density for the FCC part. 
A short equilibration run, at $\temp=2.90$, starting from a
perfect FCC lattice at this density gave the initial solid
configuration. 
The computational cell for the initial liquid part was chosen to
be the same as the one in FCC simulation 
and the initial configuration when pre-simulating the liquid part
was obtained from the FCC configuration by distributing 
vacancies to get the correct density in the liquid.
In a simulation of~\eqref{eq:algo_Smol} using a temperature, \temp,
above the melting point, \Tmelt, the sample was melted and
equilibrated. 
Afterwards the liquid was cooled to desired temperature using a 
subsequent simulation with $\temp=\Tmelt$.  

\begin{figure}[hbp]
  \centering
  \begin{tabular}{ccc}
    \begin{picture}(20,10)(0,0)
      \put(0,0){\makebox(0,0)}
    \end{picture}
    &
    \subfigure[Initial configuration, orientation 1]{
      \label{fig:conf_O1_init}
      \includegraphics[width=9.85cm]
      {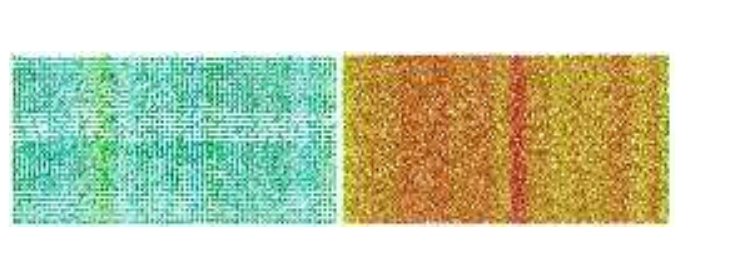}
    }
    & \multirow{4}{*}{
      \includegraphics[width = 1.2cm]
      {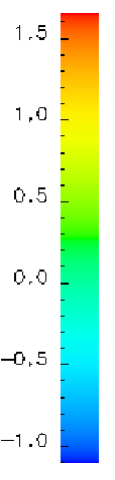}}
    \\
    &
    \subfigure[Configuration at a later time, orientation 1]{
      \label{fig:conf_O1_later}
      \includegraphics[width=10.0cm]
      {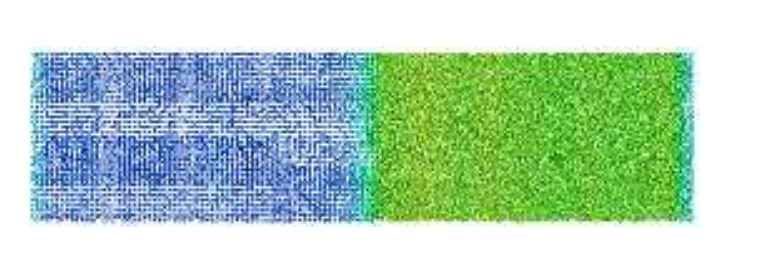}
    }
    \\  
    &
    \subfigure[Initial configuration, orientation 2]{
      \label{fig:conf_O2_init}
      \includegraphics[width=10.7cm]
      {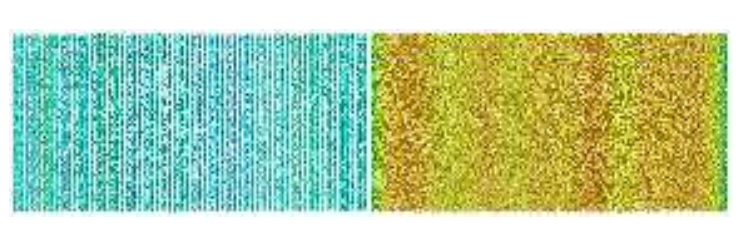}
    }
    \\  
    &
    \subfigure[Configuration at a later time, orientation 2]{
      \label{fig:conf_O2_later}
      \includegraphics[width=10.8cm]
      {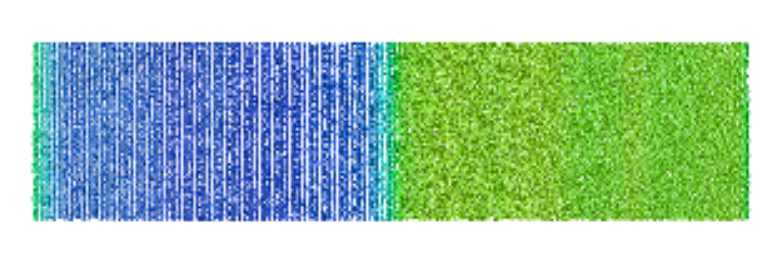}
    }
  \end{tabular}
  \caption{Snapshots of the process of setting up initial
    configurations for the two-phase simulations O1 and O2.
    The left part is solid (FCC) and the right part liquid.
    In the initial configurations, (a) and~(c), the 
    individual parts have been equilibrated at $T_{\mathrm{melt}}$
    (for the combined system), and slightly 
    compressed in one direction (to allow for two gaps).
    Subfigures~(b) and~(d) show configurations at later times when the 
    parts have expanded to fill the voids and form two interfaces.
    The atoms are coloured according to a computed phase variable;
    in~(a) and~(b) the phase variable is just the instantaneous field
    $\pfen(x_1;\mdpos[0]{\no})$, whereas~(b) and~(d) use discrete time
    averages approximating
    $\frac{1}{t_2-t_1}\int_{t_1}^{t_2}\pfen(x_1;\mdpos{\no})\,dt$. 
    \newline ~
    \newline
    Simulation O1 used 64131 particles in a computational cell of 
    dimensions $93.17\times23.29\times23.29$, while simulation O2
    used 78911 particles in a cell of dimensions 
    $100.86\times24.71\times24.96$.
  }
  \label{fig:init_confs}
\end{figure}
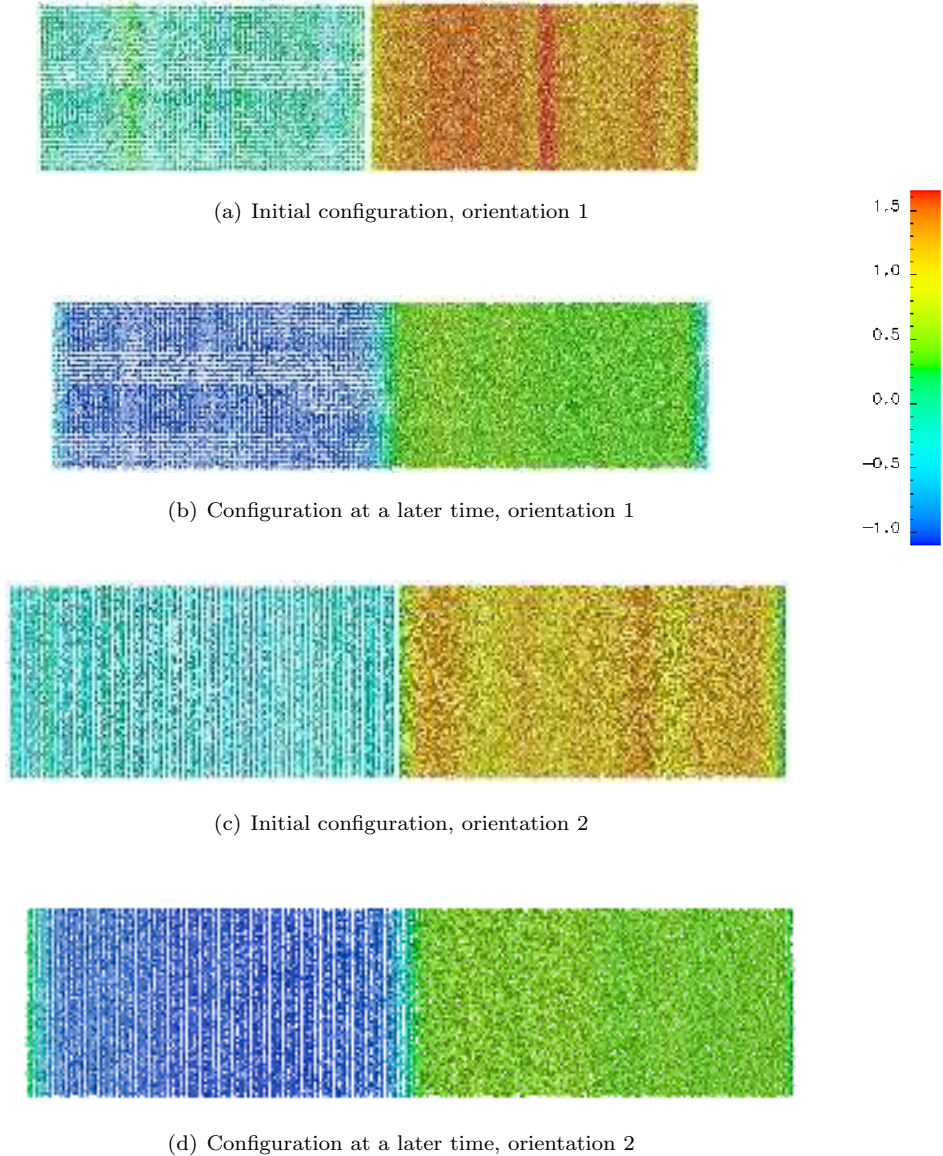
\begin{figure}[htp]
  \centering
  \subfigure[Unit cell]{
    \label{fig:fcc_unit_2}
    \includegraphics[height=1in]{./figures/unitcube.pdf}
  }
  \hspace{1cm}
  \subfigure[Orientation 1]{
    \label{fig:orient1}
    \includegraphics[height=1in]{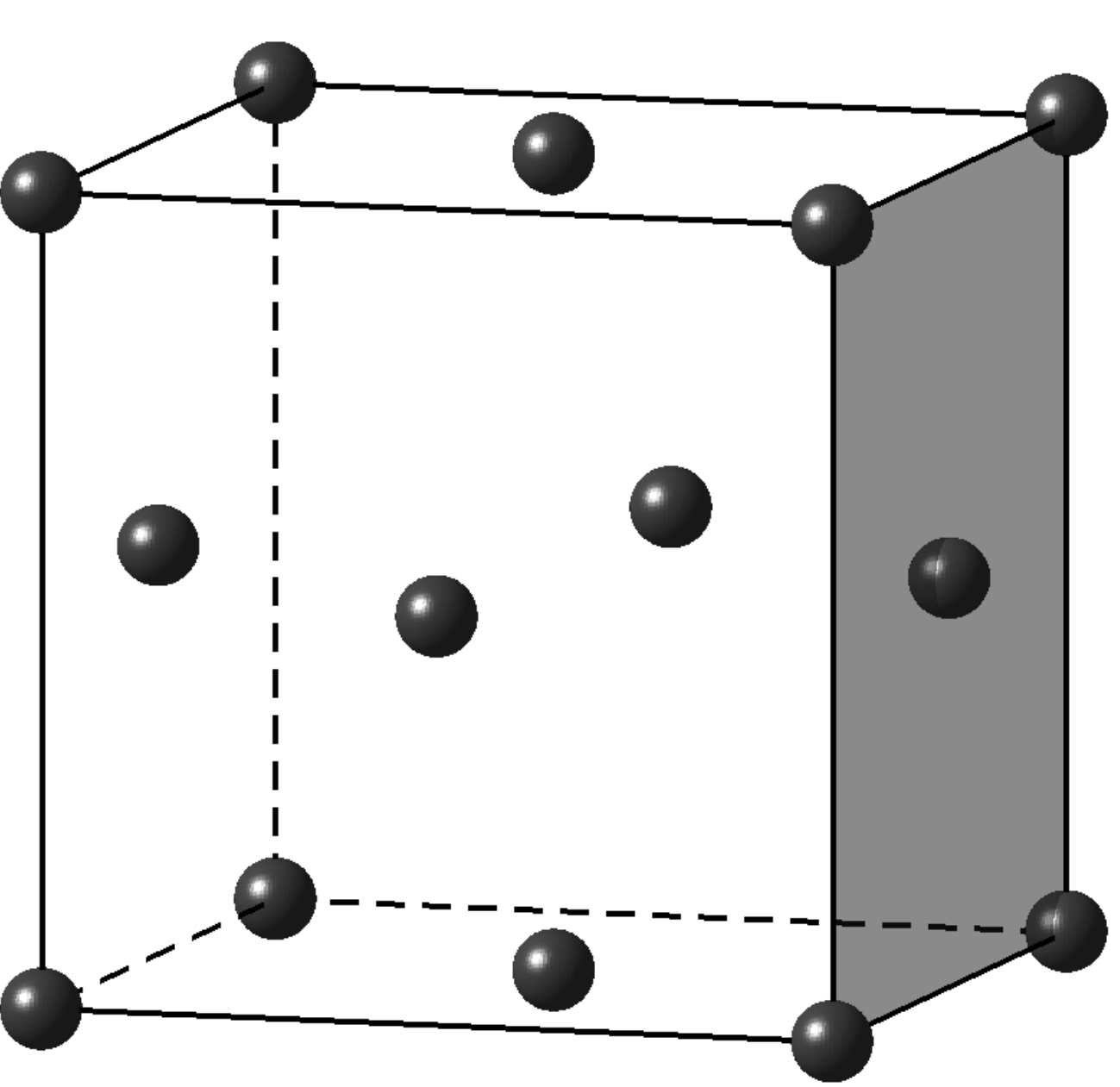}
  }
  \hspace{1cm}
  \subfigure[Orientation 2]{
    \label{fig:orient2}
    \includegraphics[height=1in]{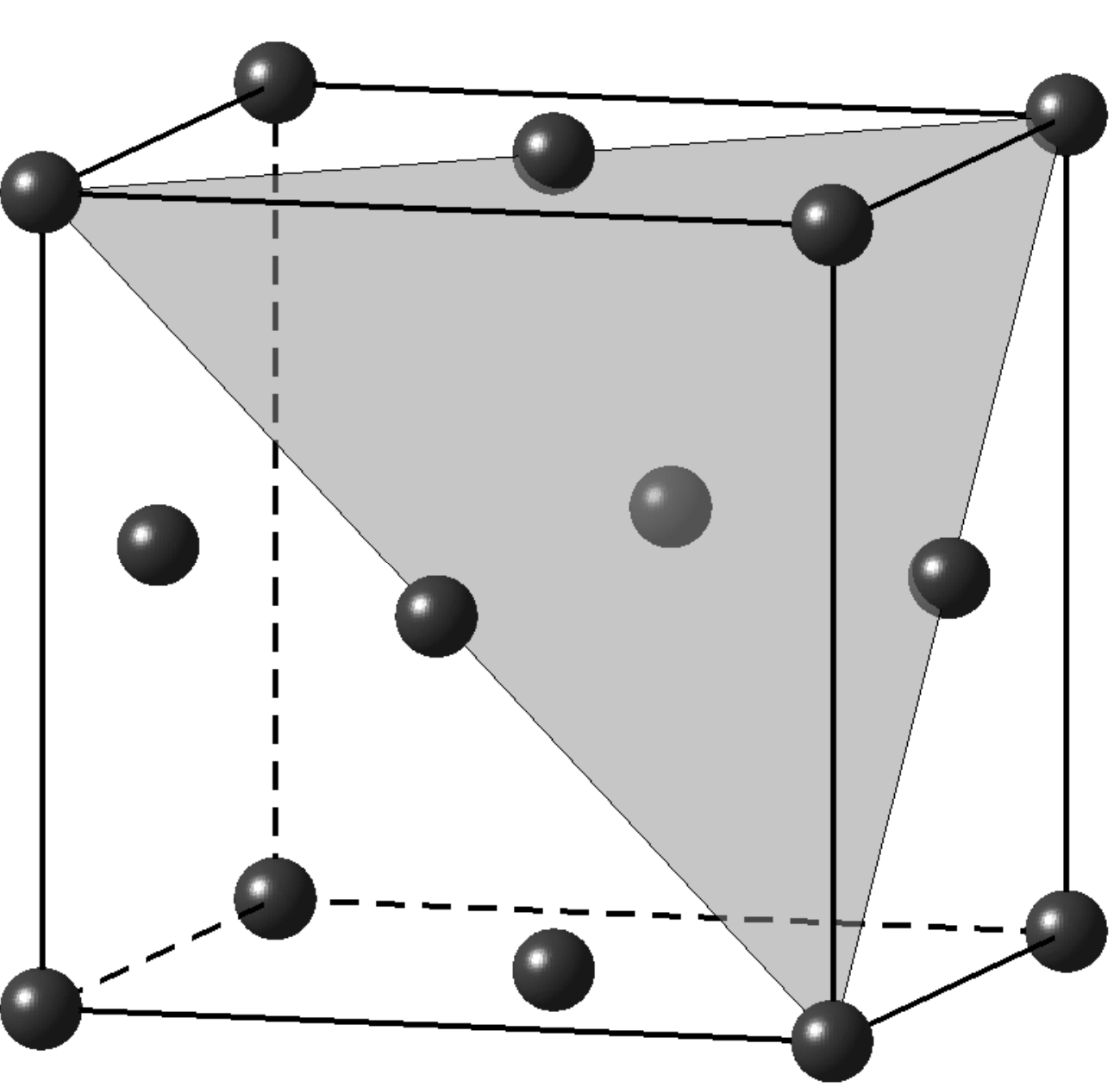}
  }
  \caption{The shaded planes in~(b) and~(c) show the two orientations of
    the solid-liquid interface with respect to the FCC lattice
    treated in the numerical simulations.} 
  \label{fig:orientations}
\end{figure}

Since no pair of atoms can be too close in the initial
configuration, 
gaps had to be introduced between the solid and liquid parts, 
but the voids could not be introduced as additional volumes in
the computational cell;  
the individual parts were equilibrated at $(N,V,T)$
corresponding to the expected densities for solid and liquid in
the combined system, so increasing the total volume would reduce 
the overall density, resulting in partial or total melting of the
solid part. To make room for the voids both the solid and the
liquid parts were compressed slightly in the direction normal to
the solid--liquid interfaces, before inserting them in their
respective volumes in the computational cell for the two-phase
simulation. 
Initial configurations obtained by this procedure are shown 
as configurations (a) and (c) in
Figure~\ref{fig:init_confs}, on page~\pageref{fig:init_confs}.
The orientation of the solid--liquid interfaces with respect to
the FCC lattice differ between the two initial configurations
shown, and these orientations with the corresponding numerical
simulations will be labelled Orientation~1 (O1) and Orientation~2
(O2) in the following.
The shaded plane in Figure~\ref{fig:orientations}(b) shows the
orientation of the interface in O1 and the shaded plane in
Figure~\ref{fig:orientations}(c) shows the orientation in O2.

Even though the compression in one direction was small,
it introduced an artificial internal stress in the system. 
The higher value of the phase-field in the subfigures~(a) and~(c)
in Figure~\ref{fig:init_confs} compared to the corresponding
regions in the subfigures~(b) and~(d) is an effect of the
compression. 
In the initial phase of the equilibration of the two-phase
system, the compressed parts expand to fill the voids.
The phase-fields in the interiors of the solid and
liquid parts in subfigures~(b) and~(d) have reached the levels
seen in the corresponding single phase systems, which shows at
least that the local potential energy contributions had returned
to normal before the production runs started.

\begin{figure}[hbp]
  \centering
  \subfigure[FCC]{
    \label{fig:rdf_fcc}
    \includegraphics[width=6cm]{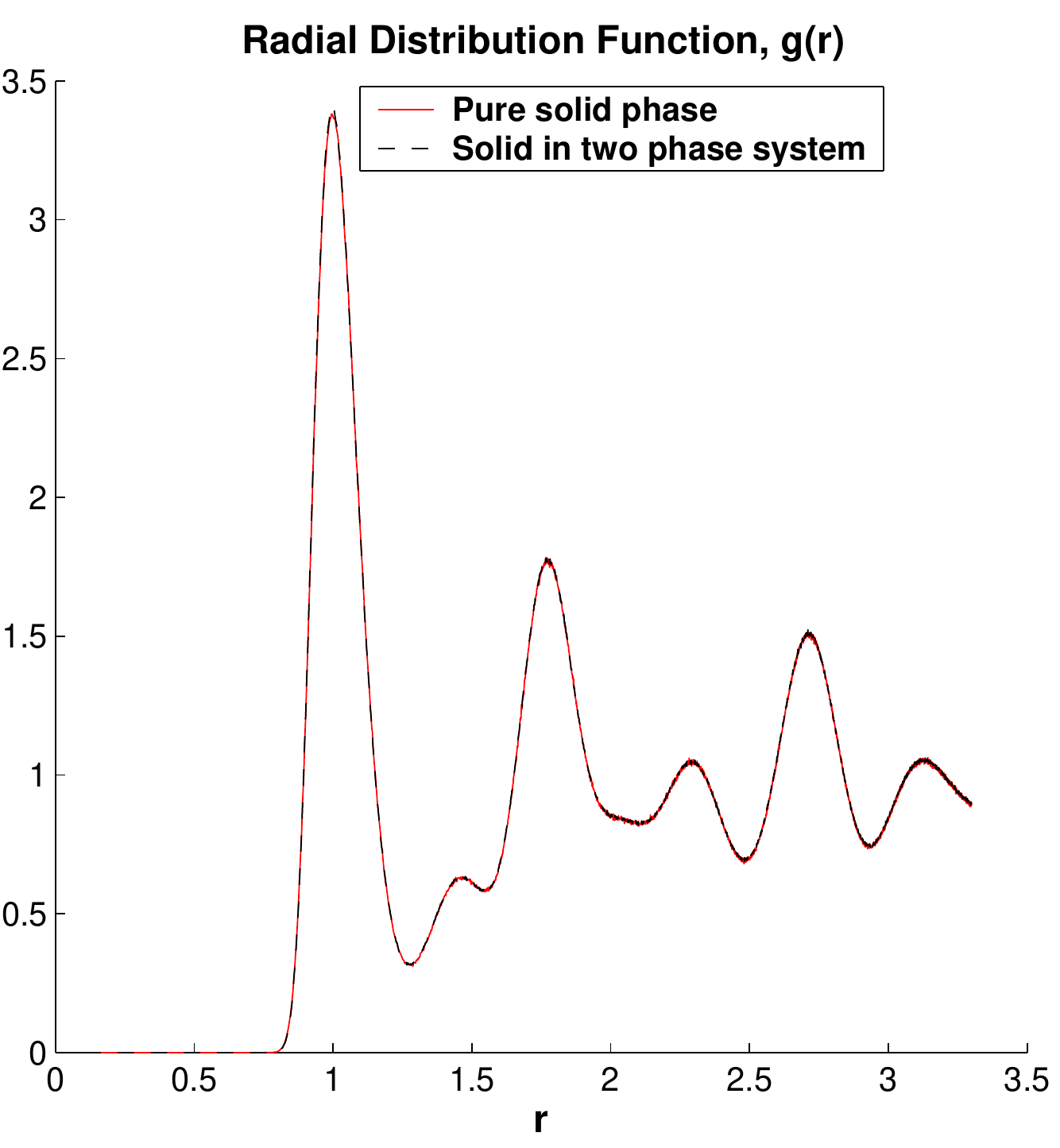}
  }
  \subfigure[Liquid]{
    \label{fig:rdf_liq}
    \includegraphics[width=6cm]{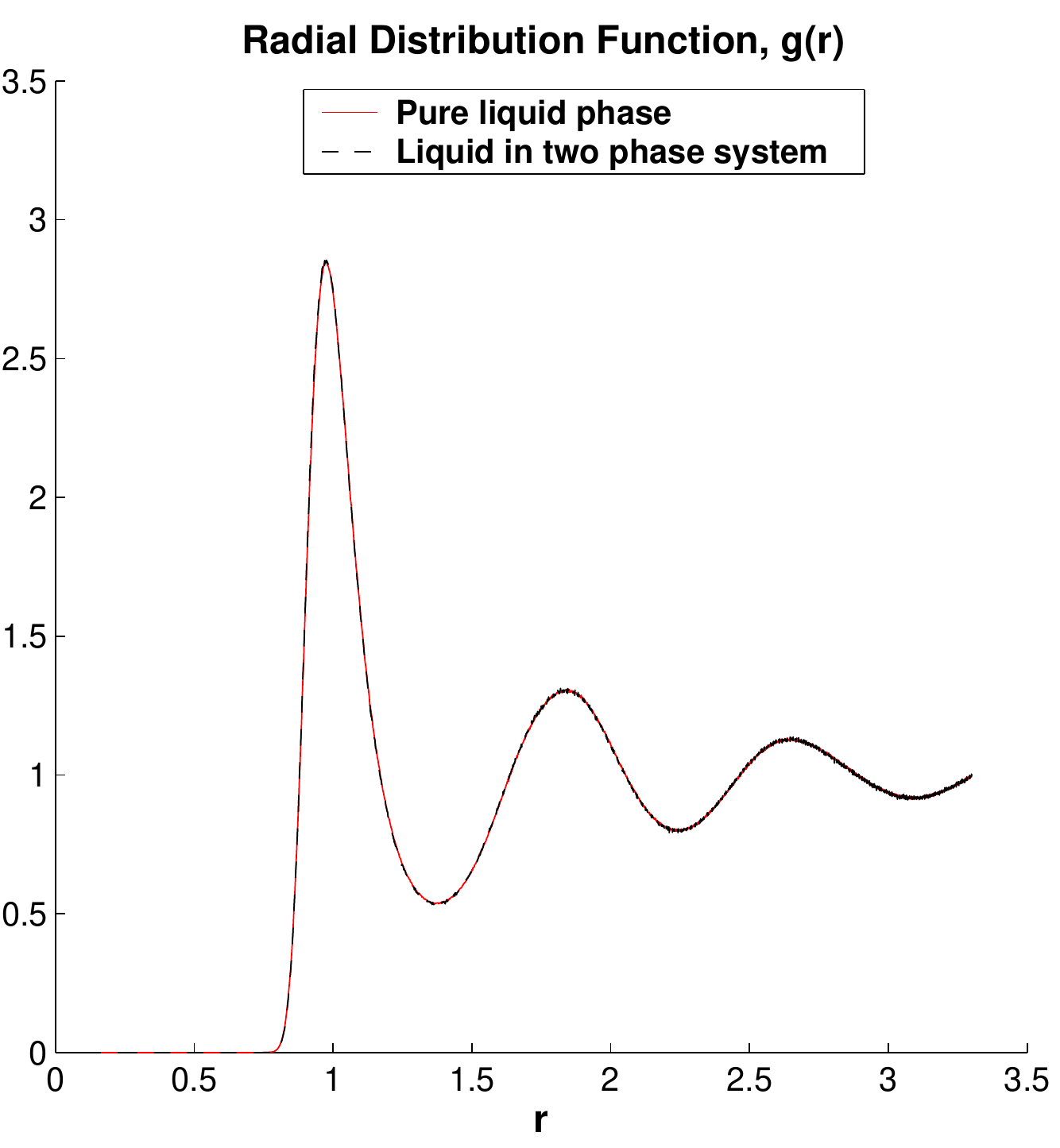}
  }
  \caption{The radial distribution function, \rdf, computed from
    several configurations, separated in time, in the process of
    setting up the two-phase system in simulation O2. 
    The solid curve shows \rdf\ computed as an average over all
    particles in the computational cell used while pre-simulating the
    solid and the liquid part, in subfigure~(a) and~(b) respectively. 
    The dashed curves show \rdf\ computed as an average over 
    particles in two slices of the computational cell of the two-phase
    system; subfigure~(a) shows \rdf\ obtained from the slice
    $5.0\leq x_1\leq 45.43$, inside the solid phase, and subfigure~(b)
    shows \rdf\ from the slice $55.43\leq x_1\leq 95.86$, inside the
    liquid phase. 
    The configurations are taken from an equilibration run,
    after the closing of the initial gaps between the pre-simulated
    phases, but before the ``production'' run.
    The radial distribution functions show good agreement between the
    single phase systems and the corresponding solid and liquid
    subdomains away from the interface.
  }
  \label{fig:rdfs}
\end{figure}

As a test of the two-phase configuration serving as initial data
in the production run, the radial distribution functions in the
interior of the two phases were computed.
The radial distribution function, \rdf, is useful for identifying
the phase of a single-phase system. For a single component system
\rdf, where $r\in\R^+$, is implicitly defined by the condition
that the average number of atoms in a spherical shell between the
radii $r_1$ and $r_2$ from the centre of any atom is  
\begin{align*}
  \density\int_{r_1}^{r_2}\rdf[r]4\pi r^2\;dr,
\end{align*}
where \density\ is the global particle density. 
In other words, the radial distribution function is the average
particle density, as a function of the separation $r$, normalised
by overall density. 
Figure~\ref{fig:rdfs}, on page~\pageref{fig:rdfs}, shows good
agreement for simulation~O2 between \rdf\ corresponding to single
phase solid and liquid configurations and \rdf\ computed in the
interior of the two phases, excluding two intervals of length
10.0 in the interface regions. 

An effect of the finite size of the computational cell is that
periodic boundary conditions may interact with the solid and
affect the results; here the computational cell was chosen to
match the FCC structure in a specific orientation with respect to
the box and thus stabilises the structure and orientation. 
It is important to know that the density in the FCC part
(and hence the box cross section) is consistent with constant
pressure simulations close to the melting point.
A related question is whether the length of the computational box
is large enough for properties around the interfaces in the
infinitely layered structure to be good approximations of those
near an interface between a solid and liquid on the macroscopic
scale.

\subsection{Computation Of the Coarse-Grained Model Functions}
\label{sec:coarse}

The coefficient functions~\eqref{eq:def_coeff_cg} in the
stochastic differential equation~\eqref{eq:sde_pfcg} for the
coarse-grained phase-field are defined in terms of the time
averaged expected values~\eqref{eq:def_coeffcgx}
and~\eqref{eq:def_mdphaseav} on the form
\begin{align*}
  \frac{1}{\tend}
  \E\left[\int_0^\tend\funcxX{\mdpos{\no}}\;\bigg|\;
    \mdpos[0]{\no}=\mdpos[\no]{0}
  \right],
\end{align*}
where $\mdpos[\no]{0}$ is a configuration of a stationary
two-phase system. 
By setting up an initial configuration, $\mdpos[\no]{0}$, as
described in the previous section, and simulating discrete sample
trajectories using the Euler-Maruyama
method~\eqref{eq:algo_Smol}, a sequence of configurations 
$\{\mdposd[k]{\no}\}_{k=1}^K$ 
approximating the sequence $\{\mdpos[t_k]{\no}\}_{k=1}^K$
for some times $0<t_1<\cdots<t_K=\tend$, is obtained.
In a post processing step a set of configurations
$\sampleset\subseteq\{\mdposd[k]{\no}\}_{k=1}^K$
is selected and
averages 
\begin{align*}
  \average{\funcx}{\sampleset} & =
  \sum_{\mdposgen\in\sampleset}\funcxX{\mdposgen}w_\mdposgen,
\end{align*}
consistently weighted with weights $w_\mdposgen$, 
are computed as approximations of the corresponding expected
values in the continuous time model.
It is usually more efficient not to include every configuration
in the averages. 
This will be discussed in Section~\ref{sec:results}.

As described in the introduction, the averages are functions of
the coordinate direction $x_1$, normal to the planar interface,
since the mollifier in the definition~\eqref{eq:def_pfen} of the 
microscale phase-field, \pfen, is chosen to take uniform averages
in the planes parallel to the interface. The mollifier used in
the computations is 
\begin{align}
  \label{eq:mollifier}
  \molli(x) & = \molli(x_1) = c
  \exp{
    \left(-\frac{1}{2}\left(\frac{x_1}{\epsilon}\right)^2\right)
  }\I{|x_1|<R_c},
\end{align}
where $c$ is a normalising constant, \molliscale\ is a smoothing
parameter, and $R_c$ is a cut-off. The smoothing parameter is on
the order of typical nearest neighbour distances,
$\molliscale\approx1$, and 
$R_c=6\molliscale$, for all choices of \molliscale, which gives
$\molli(R_c)\approx1.5\cdot10^{-8}\molli(0)$; the shape of
\molli\ can be seen in Figure~\ref{fig:mollifier}, on
page~\pageref{fig:mollifier}. 

An explicit derivation of expressions for the drift and
the diffusion is given in Appendix~\ref{sec:calculations}. 
Separating the drift in terms containing two, one, and zero,
derivatives of the mollifier, the right hand side
of~\eqref{eq:def_driftcgx} is approximated by 
\begin{align*}
  \kb\temp \ddxett \average{\pfen}{\sampleset}
  + \dxett \average{\driftone}{\sampleset} 
  + \average{\driftzero}{\sampleset},
\end{align*}
where 
\begin{align}
  \label{eq:driftone}
  \driftone(x;\mdposgen) & = 
  \sumall{j}
  (\kb\temp - \potenp{j}(\mdposgen))[\force_j(\mdposgen)]_1
  \molli(x-\mdpos[\no]{j})
  \intertext{and}
  \nonumber
  \driftzero(x;\mdposgen) & = 
  - \sumall{j}
  \left(\kb\temp \dX\cdot\force_j(\mdposgen) 
    + \frac{1}{2} ||\force_j(\mdposgen)||^2 \right) 
  \molli(x-\mdpos[\no]{j}) 
  \\ & \quad 
  \label{eq:driftzero}
  - \frac{1}{2} 
  \sumall{j} \sumneq{i}{j}
  \pairforce_{ij}(\mdposgen) \cdot 
  \force_j(\mdposgen) \molli(x-\mdpos[\no]{i}).
\end{align}
Here $\force_j$ is the total force acting on particle $j$, 
$[\force_j(\mdposgen)]_1$ is the $x_1$-component of the force, 
and $\pairforce_{ij}$ are the contributions from individual pairs,
\begin{align*}
  \force_j(\mdposgen) = 
  - \dX \totpot(\mdposgen)
  & = \sumneq{i}{j}  
  \pairpot'(|| \mdpos[\no]{i}-\mdpos[\no]{j} ||) 
  \frac{\mdpos[\no]{i}-\mdpos[\no]{j}}
  {|| \mdpos[\no]{i}-\mdpos[\no]{j} ||}
  = \sumneq{i}{j} \pairforce_{ij}(\mdposgen).
\end{align*}
The right hand side in equation~\eqref{eq:def_diffucgx}, for the
coarse grained diffusion, is approximated by 
\begin{align}
  \label{eq:diffumat}
  \diffumat(\cdot,\cdot) & = 
  \average{2\kb\temp\sumall{j}
    \bigl(p_j(\cdot,\cdot;X)+q_j(\cdot,\cdot;X)\bigr)}{\sampleset},
\end{align}
where
\begin{align*}
  p_j(x,y;\mdposgen) & =
  \left(\frac{\potenp{j}(\mdposgen)}{\epsilon^2}\right)^2
  \left[x-\mdpos[\no]{j}\right]_1 \left[y-\mdpos[\no]{j}\right]_1
  \molli(x-\mdpos[\no]{j})\molli(y-\mdpos[\no]{j})\\
  & \quad - 
  \frac{\potenp{j}(\mdposgen)}{2\epsilon^2} 
  [x-\mdpos[\no]{j}]_1 \molli(x-\mdpos[\no]{j})
  \biggl(
  [\force_j(\mdposgen)]_1 \molli(y-\mdpos[\no]{j})
  + \!\!\sumneq{i}{j}
  [\pairforce_{ij}(\mdposgen)]_1 \molli(y-\mdpos[\no]{i})
  \biggr)
  \\
  & \quad - 
  \frac{\potenp{j}(\mdposgen)}{2\epsilon^2} 
  [y-\mdpos[\no]{j}]_1 \molli(y-\mdpos[\no]{j})
  \biggl(
  [\force_j(\mdposgen)]_1 \molli(x-\mdpos[\no]{j})
  + \!\!\sumneq{i}{j}
  [\pairforce_{ij}(\mdposgen)]_1 \molli(x-\mdpos[\no]{i})
  \biggr)
  \intertext{and}
  q_j(x,y;\mdposgen) & = 
  \frac{1}{4}
  \biggl( 
  \force_j(\mdposgen) \molli(x-\mdpos[\no]{j})
  + \sumneq{i}{j} \pairforce_{ij}(\mdposgen) 
  \molli(x-\mdpos[\no]{i}) 
  \biggr) \\
  & \quad\quad
  \cdot
  \biggl( 
  \force_j(\mdposgen) \molli(y-\mdpos[\no]{j})
  + \sumneq{i}{j} \pairforce_{ij}(\mdposgen) 
  \molli(y-\mdpos[\no]{i}) 
  \biggr).
\end{align*}

The functions
$\average{\funcx}{\sampleset}$ are computed in a
discrete set of points $D_K=\{x_1^i\}_{i=1}^K$ along the $x_1$
axis of the molecular dynamics domain. 
This makes the computed components,
$\average{\pfen}{\sampleset}$, $\average{\driftone}{\sampleset}$, 
and $\average{\driftzero}{\sampleset}$, of the drift coefficient
function $K$-vectors and the 
computed \diffumat\ a $K$-by-$K$ matrix. 
The individual diffusion coefficient functions $\diffucgx_j$ are
obtained by taking the square root of the computed diffusion
matrix, $\diffumat=\diffumat^{1/2}\transpose{(\diffumat^{1/2})}$,
and letting the $j$:th column of $\diffumat^{1/2}$ define
$\diffucgx_j$. 
While an exact computation would produce a symmetric positive
semi definite matrix \diffumat, finite precision effects make
some computed eigenvalues negative, but small in absolute value. 
In an eigenvector factorisation of \diffumat, 
let $\eigenvalmat$ denote a diagonal matrix with all eigenvalues
of \diffumat\ and $\poseigenval$ a smaller diagonal matrix
containing the dominant, possibly all, of the positive
eigenvalues but no negative ones. Let $\eigenbasemat$ and
$\poseigenmat$ be the matrices of the corresponding eigenvectors. 
Then the square root of the matrix $\poseigenval$ is a real
diagonal matrix which can be used in the approximation 
\begin{align}
  \label{eq:def_diffucgxmat}
  \diffumat & = 
  \eigenbasemat \eigenvalmat \transpose{\eigenbasemat}
  \approx
  \poseigenmat \poseigenval \transpose{\poseigenmat}
  = \left(
    \poseigenmat \poseigenval^{1/2} \transpose{\poseigenmat}
  \right)
  \transpose{\left(
      \poseigenmat \poseigenval^{1/2} \transpose{\poseigenmat}
    \right)}
  \backdefeq \diffucgxmat \transpose{\diffucgxmat}.
\end{align}
With one Wiener process $\indepw_j$ in the coarse-grained
stochastic differential equation~\eqref{eq:sde_pfcg} per
evaluation point, $K=\nrw$, the component vectors, $\diffucgx_j$,
of the diffusion in coarse-grained equation can be defined as the
column vectors of the matrix \diffucgxmat, to obtain
\begin{align*}
  \sumall[\nrw]{j}\diffucgx_j\transpose{\diffucgx_j}
  & \approx
  \diffumat.
\end{align*}
If two grid points, $x_1$ and $y_1$, are further apart than twice
the sum of the cut-off in the potential and the cut-off in the
mollifier, then $p_j(x,y;\cdot)$ and $q_j(x,y;\cdot)$ is zero;
hence a natural ordering $x_1^1<x_1^2<\cdots<x_1^K$ of the grid
points makes \diffumat\ a band matrix. 
The definition of \diffucgxmat\ in~\eqref{eq:def_diffucgxmat}
preserves the connection between grid points and diffusion
functions and the dominating terms in a tabulated vector
$\diffucgx_j$ are those of nearby grid points.

\section{Results}
\label{sec:results}

This section describes results from numerical simulations
performed to compute the coarse-grained model functions.
The value of the smoothing parameter \molliscale\ in the
mollifier is $1.0$, unless another value is specified.

\subsection{The averaged phase-field 
  $\pfen_\mathrm{av}\approx\average{\pfen}{\sampleset}$}
\label{sec:computed_phase-field}

\begin{figure}[hbp]
  \centering
  \subfigure[Orientation 1]{
    \label{subfig:pf_O1}
    \includegraphics[width=6.5cm]
    {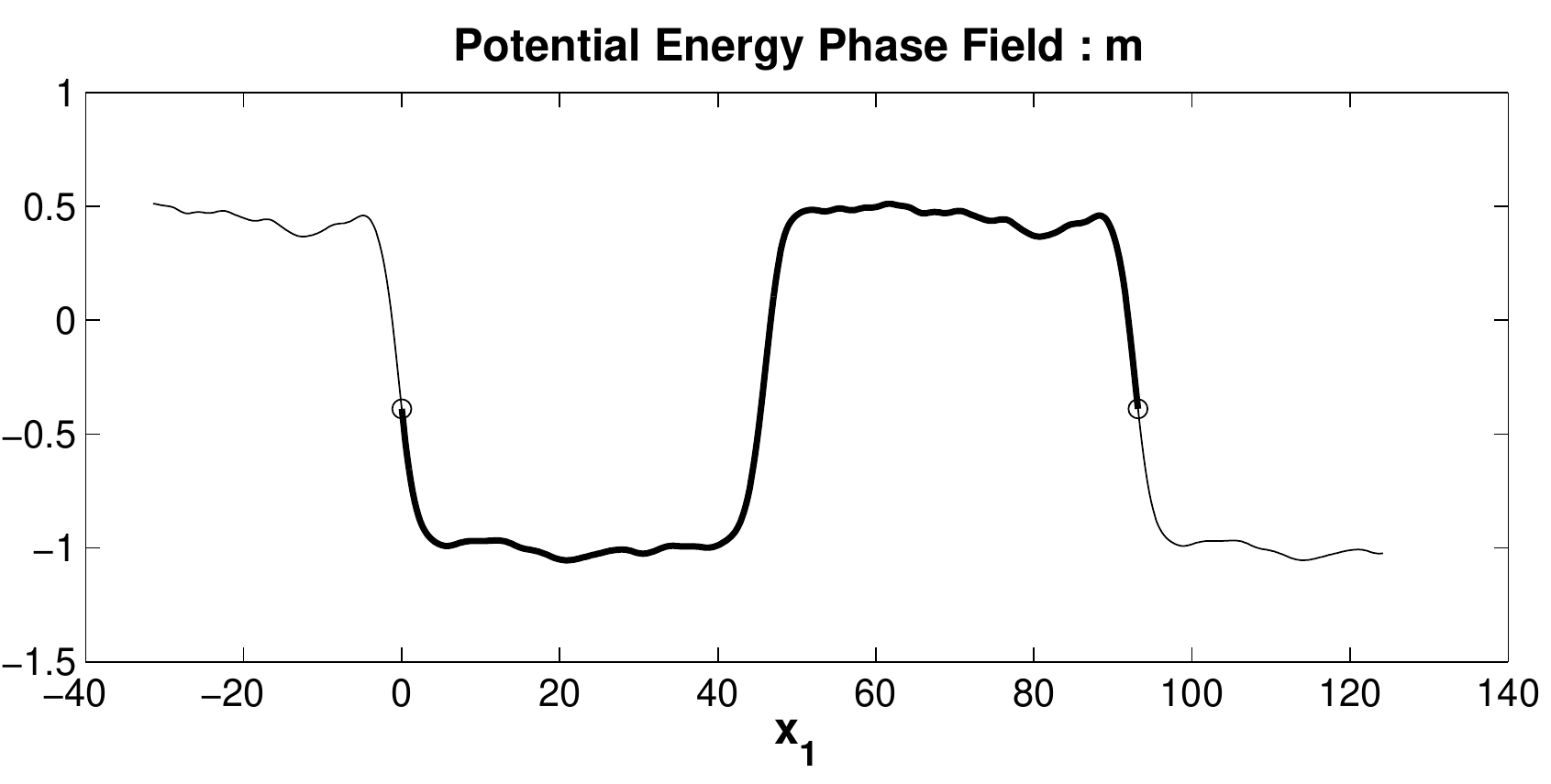}
  }
  \subfigure[Orientation 1]{
    \label{subfig:rho_pf_O1}
    \includegraphics[width=6.5cm]
    {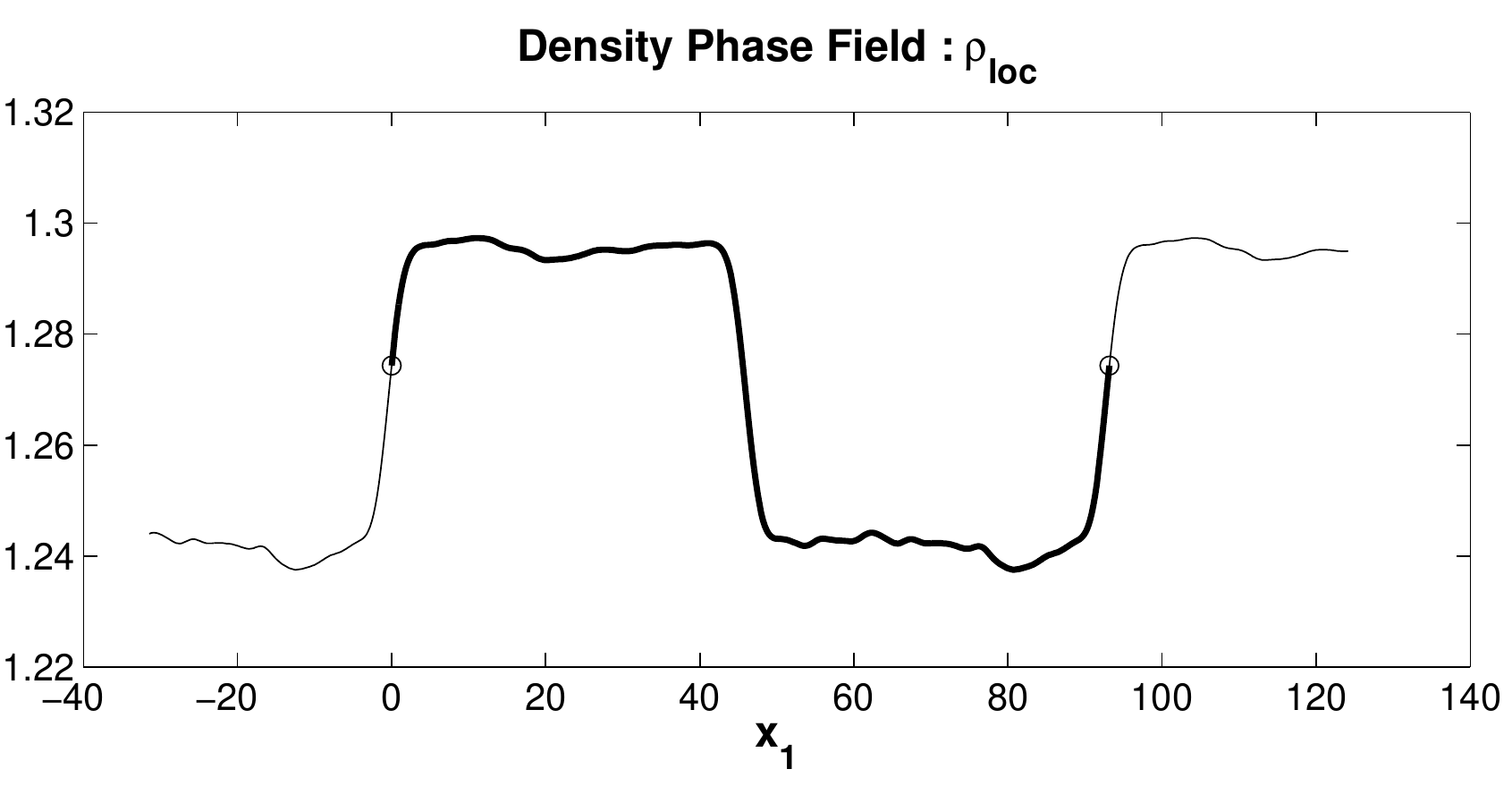}
  }
  \subfigure[Orientation 2]{
    \label{subfig:pf_O2}
    \includegraphics[width=6.5cm]
    {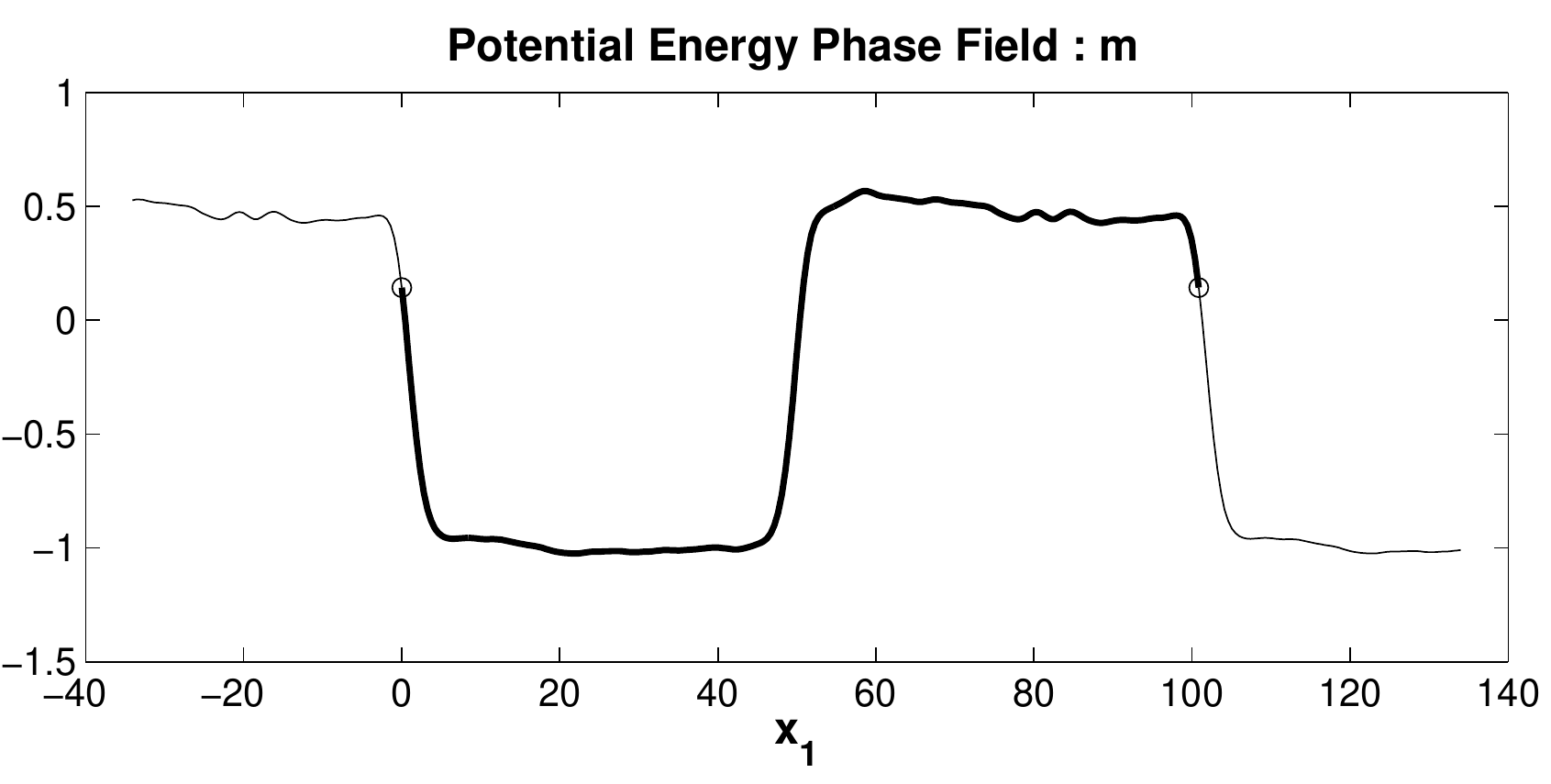}
  }
  \subfigure[Orientation 2]{
    \label{subfig:rho_pf_O2}
    \includegraphics[width=6.5cm]
    {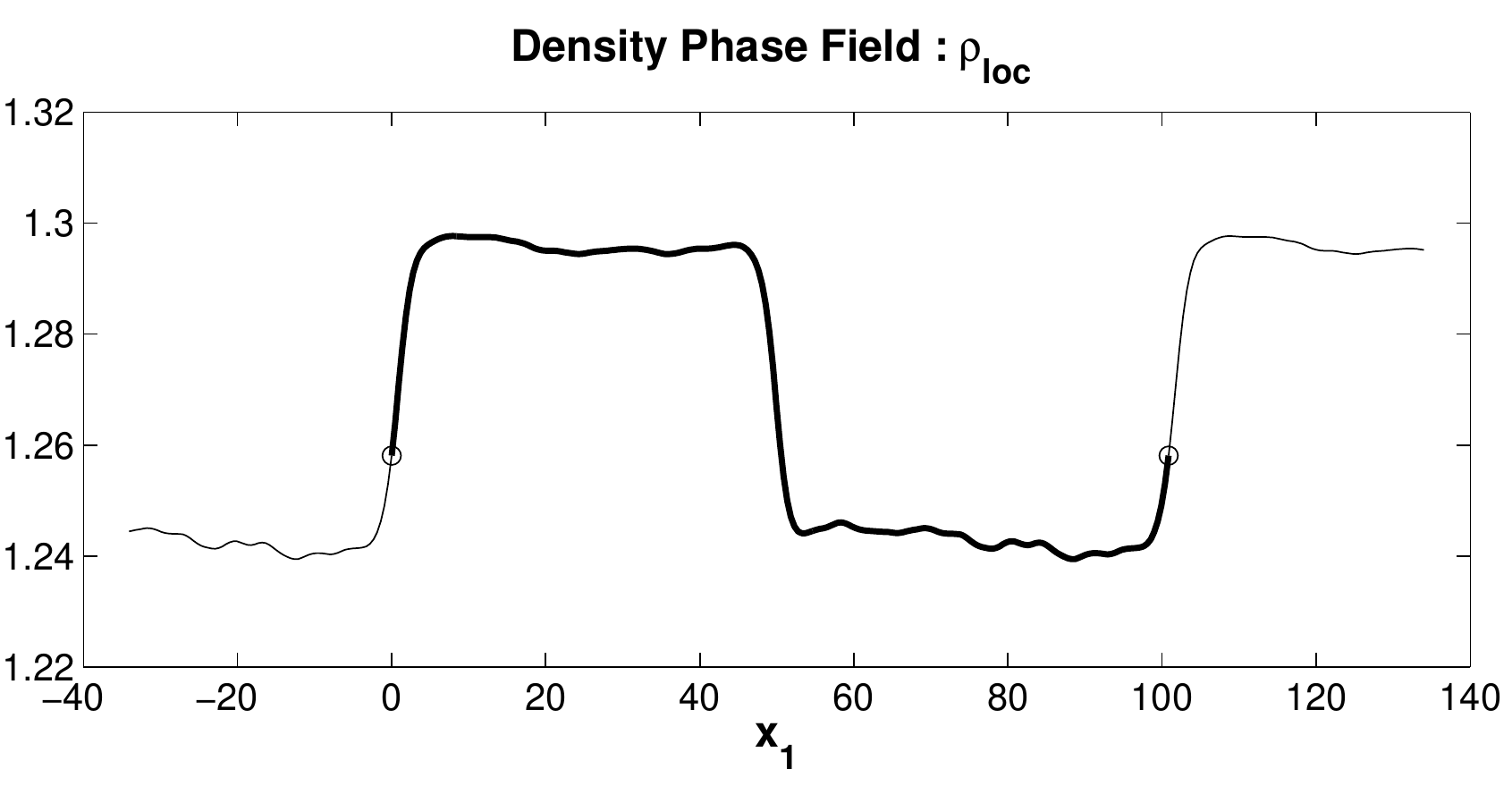}
  }
  \subfigure[Orientation 2]{
    \label{subfig:var_pf_O2}
    \includegraphics[width=6.5cm]
    {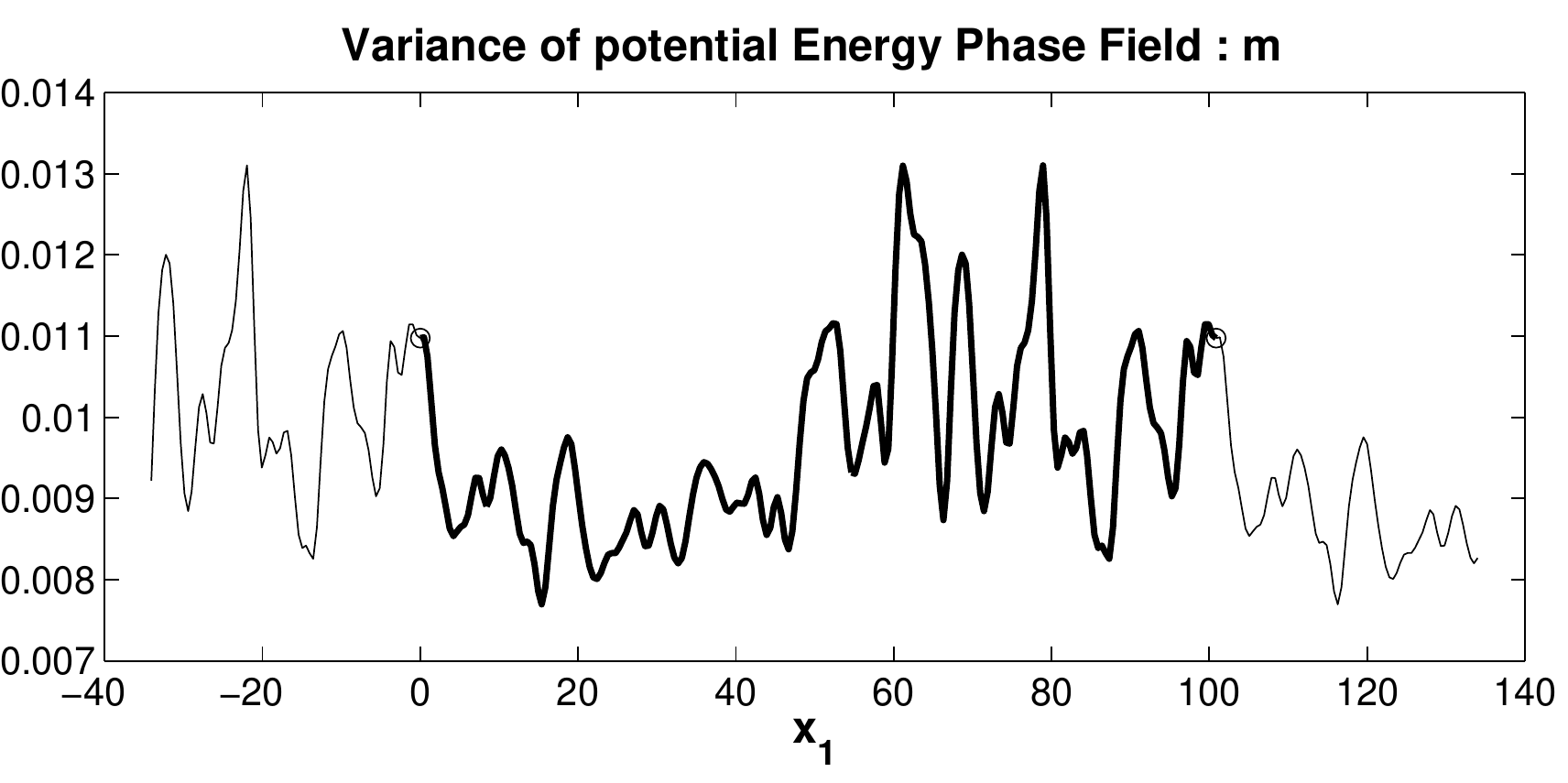}
  }
  \subfigure[Orientation 2]{
    \label{subfig:var_rho_pf_O2}
    \includegraphics[width=6.5cm]
    {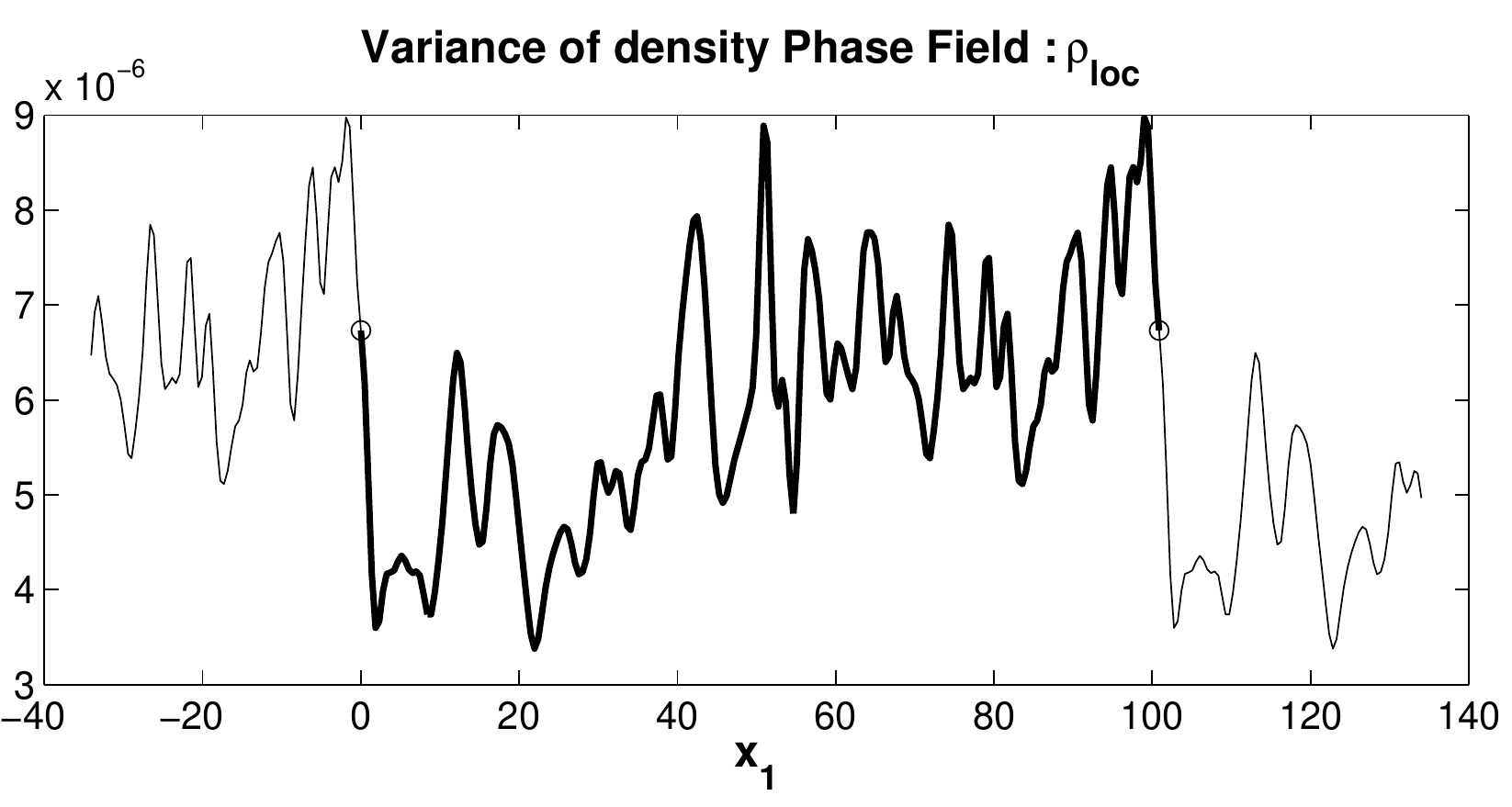}
  }
  \caption{Subfigures (a) and (c) show the potential energy
    phase-field, $\average{\pfen}{\sampleset}$, computed from
    simulations O1 and O2, respectively.
    Subfigures (b) and (d) show the corresponding spatially
    averaged particle densities.
    Subfigures (e) and (f) show the pointwise sample variance
    associated with the averages in (c) and (d).
    The thick parts of the curves show the computed functions
    in molecular dynamics cell.
    The thinner parts show the periodic continuations across the
    boundaries of the cell, marked by circles.    
    The averages in simulation O1, and O2, were formed over 1721,
    and 1775, configurations separated in time by $5\cdot10^{-4}$, so that
    the total time from first to last configuration was 0.860, 
    and 0.8875, respectively.
    The high frequency fluctuations are small after averaging on
    this time scale, but larger fluctuations remain in both
    phases. This suggests that the two phase system is not yet
    equilibrated. 
    Still the computed phase-fields appear qualitatively correct.}  
  \label{fig:phasefield_not_eq}
\end{figure}
\begin{figure}[htp]
  \centering
  \subfigure[Density Phase-Field, O2]{
    \label{subfig:density_levels_O2}
    \includegraphics[width=6.5cm]
    {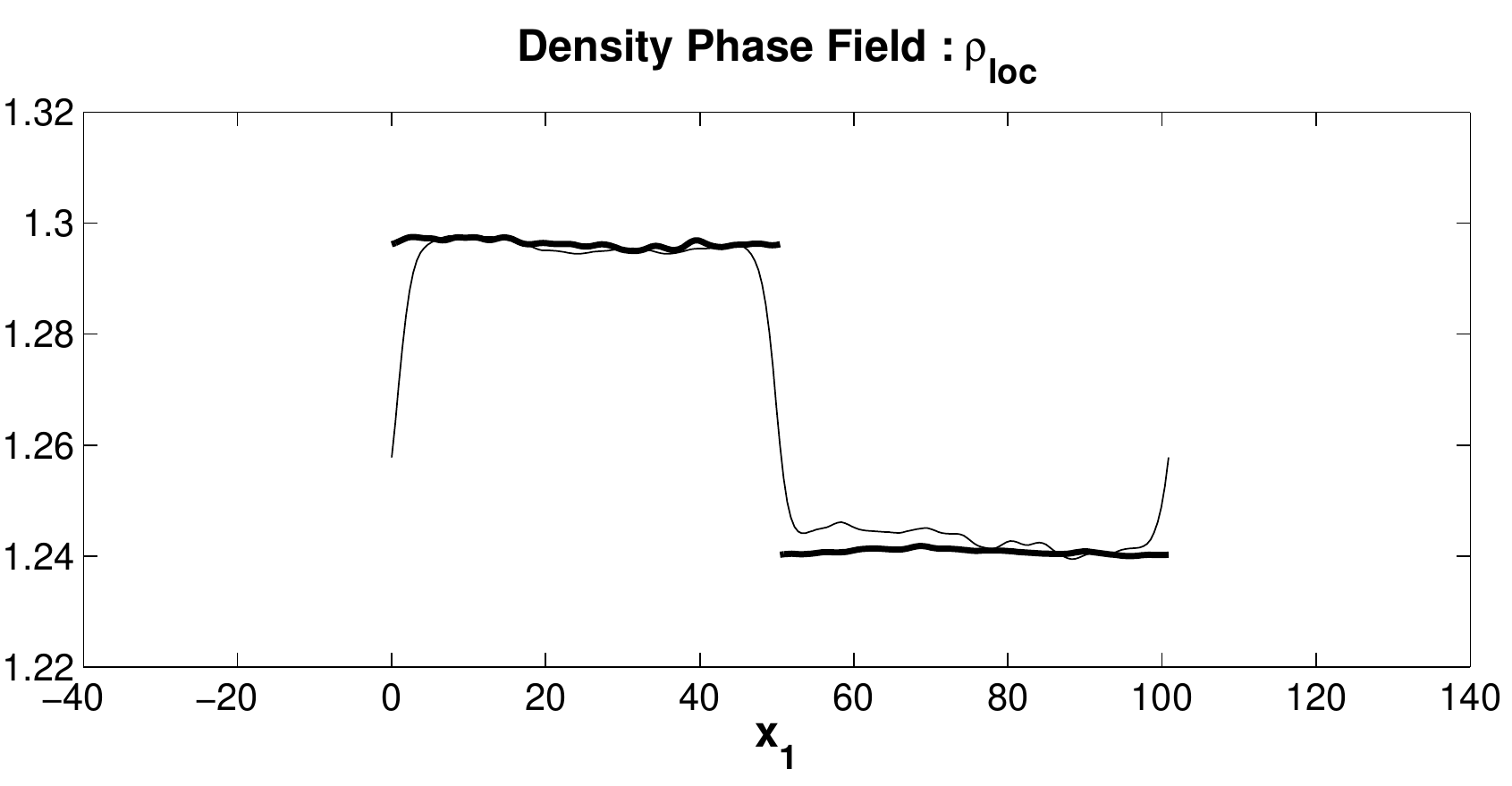}
  }
  \subfigure[Potential Energy Phase-Field, O2]{
    \label{subfig:poten_levels_O2}
    \includegraphics[width=6.5cm]
    {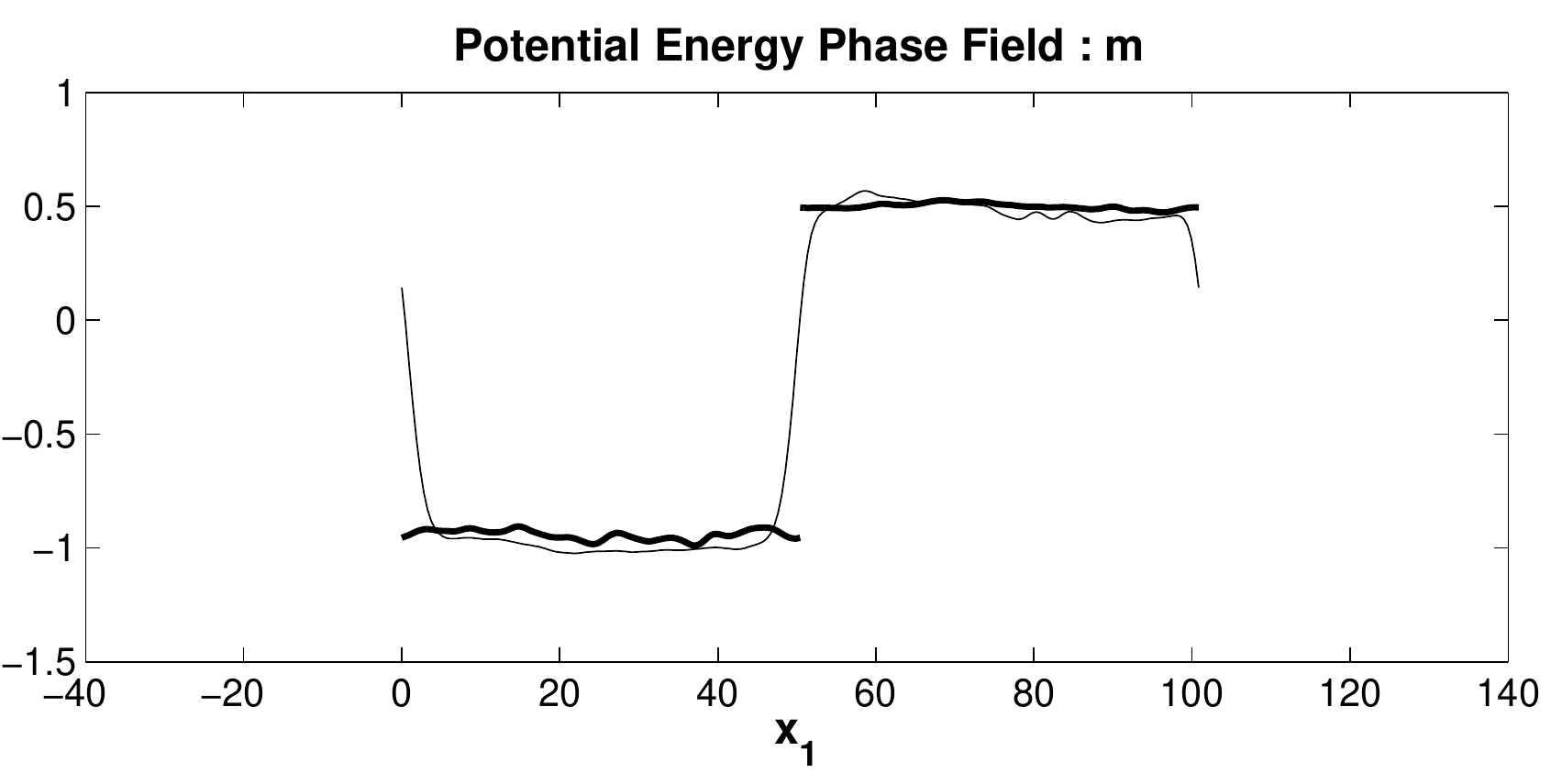}
  }
  \caption{The computational cell in the molecular dynamics
    simulations must be sufficiently large for the infinitely
    layered structure to resemble a system with a single
    solid--liquid interface on the macroscopic scale.
    In simulation O2 the total length of the computational cell
    was 100.86; subfigure~(b) shows that this was sufficient for
    the averaged phase-field, $\average{\pfen}{\sampleset}$, to
    obtain values in the interior of each phase that are similar
    to the functions, marked by thick curves, obtained in the
    single phase configurations simulated during the setup of
    simulation O2. 
  }  
  \label{fig:levels_O2}
\end{figure}

The first observation is that during the time intervals of the
molecular dynamics simulations, the interfaces between the solid
and the liquid subdomains were sufficiently stable for the
averaged potential energy phase-fields,
$\average{\pfen}{\sampleset}$, to appear qualitatively right. The 
phase-field appears to have two distinct equilibrium values,
corresponding to the solid and liquid subdomains, and the
transitions between the two regions are smooth and occur over
distances of a few nearest neighbour distances; see
Figure~\ref{fig:phasefield_not_eq}. 
Figure~\ref{subfig:poten_levels_O2} shows that the computational
cells in the molecular dynamics simulations are large enough for
the phase-field in the interior of the two phases to attain
values similar to the values in the corresponding single phase
simulations. In simulations with a cubic,
$23.29\times23.29\times23.29$, computational cell the gap between the
phase-field levels in the solid and the liquid was significantly
smaller, which indicates that the length of the computational
cell can not be taken much smaller than in simulations O1 and O2.
It is still possible that further increasing the size of the
computational cell may affect the results.

\subsection{The averaged drift 
  $\driftcgx\approx\average{\driftmd}{\sampleset}$}
\label{sec:computed_drift}

\begin{figure}[hbp]
  \centering
  \includegraphics[width=6.5cm]
  {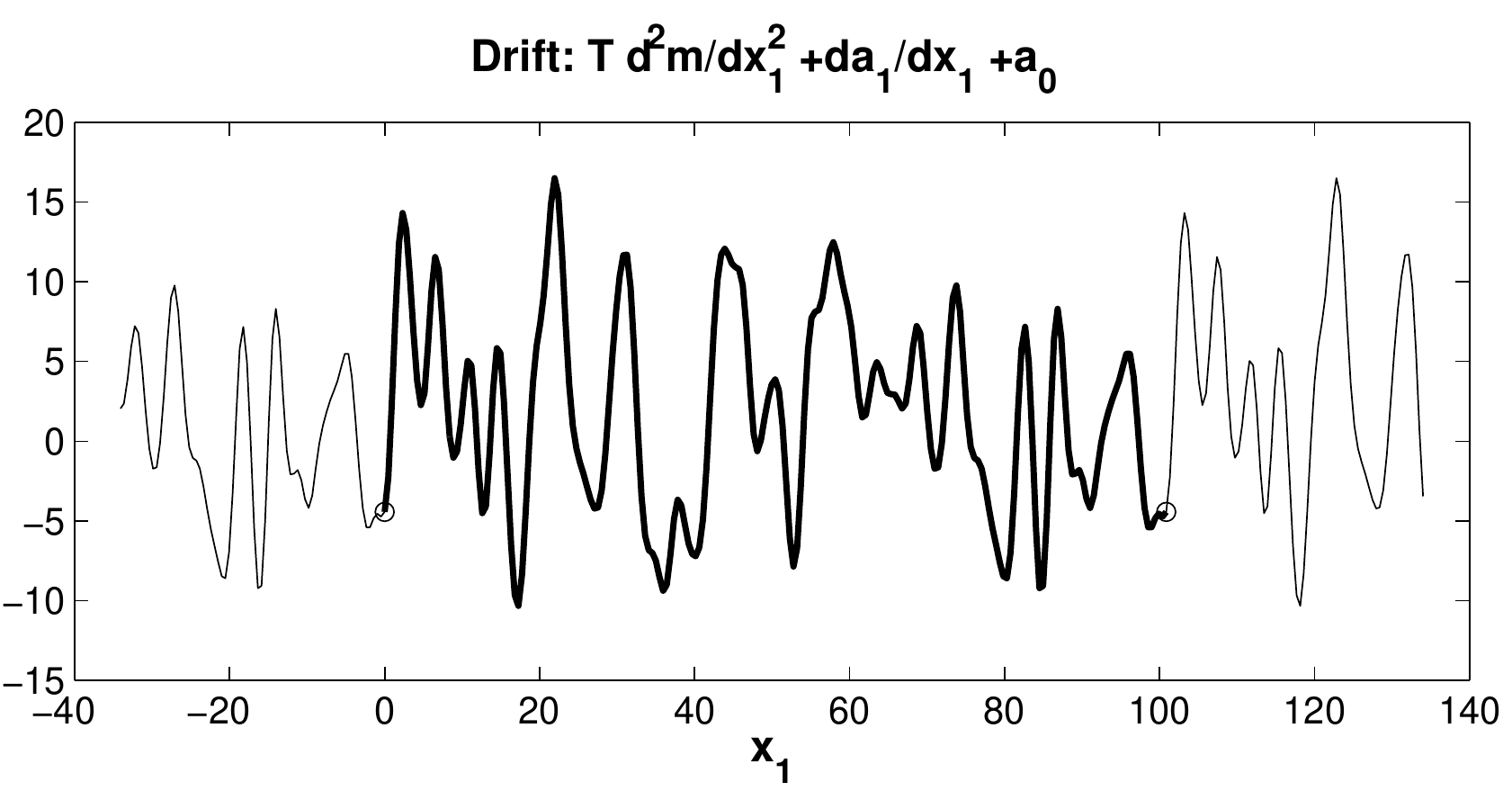}
  \caption{The average total drift,
  $\average{\driftmd}{\sampleset}$, based on the same 1775
  configurations from simulation O2 as
  $\average{\pfen}{\sampleset}$ in Figure~\ref{subfig:pf_O2}, is
  still dominated by large oscillations.} 
  \label{fig:drift_not_zero}
\end{figure}

The average $\average{\pfen}{\sampleset}$ approximates the
expected time average~\eqref{eq:def_mdphaseav}. 
The next expected value to study is the one defining the coarse
grained drift in~\eqref{eq:def_driftcgx}.
In a stationary situation, where the interfaces do not move
during the simulation and the averaged phase-field converges to a
stationary profile, the average total drift in the stochastic
differential equation describing the phase-field variable must
converge to zero.
Still the time averaged total drift corresponding to the
simulation O2, whose averaged phase-field was discussed above, is
far from zero; see Figure~\ref{fig:drift_not_zero}.
The computed time averaged drift
\begin{align*}
  \average{\driftmd(x;\mdposd{\no})}{\sampleset}
  & \approx
  \frac{1}{\tend}
  \E\biggl[\int_0^\tend \driftmd(x;\mdpos{\no})\,dt
  \;\Bigl\lvert\;
  \mdpos[0]{\no}=\mdpos[\no]{0}
  \biggr]
\end{align*}
depends both on the length of the time interval where the average
is computed, the number of configurations used in the average,
and on the discrete approximation \mdposd{\no}\ of \mdpos[t_n]{\no};
these potential error sources must be analysed to explain the
result. 

\subsubsection{The effect of discrete time dynamics}

\begin{figure}[hbp]
  \centering
  \subfigure[\rdf, FCC, different $\Delta t$]{
    \label{subfig:rdf}
    \includegraphics[width = 5cm]{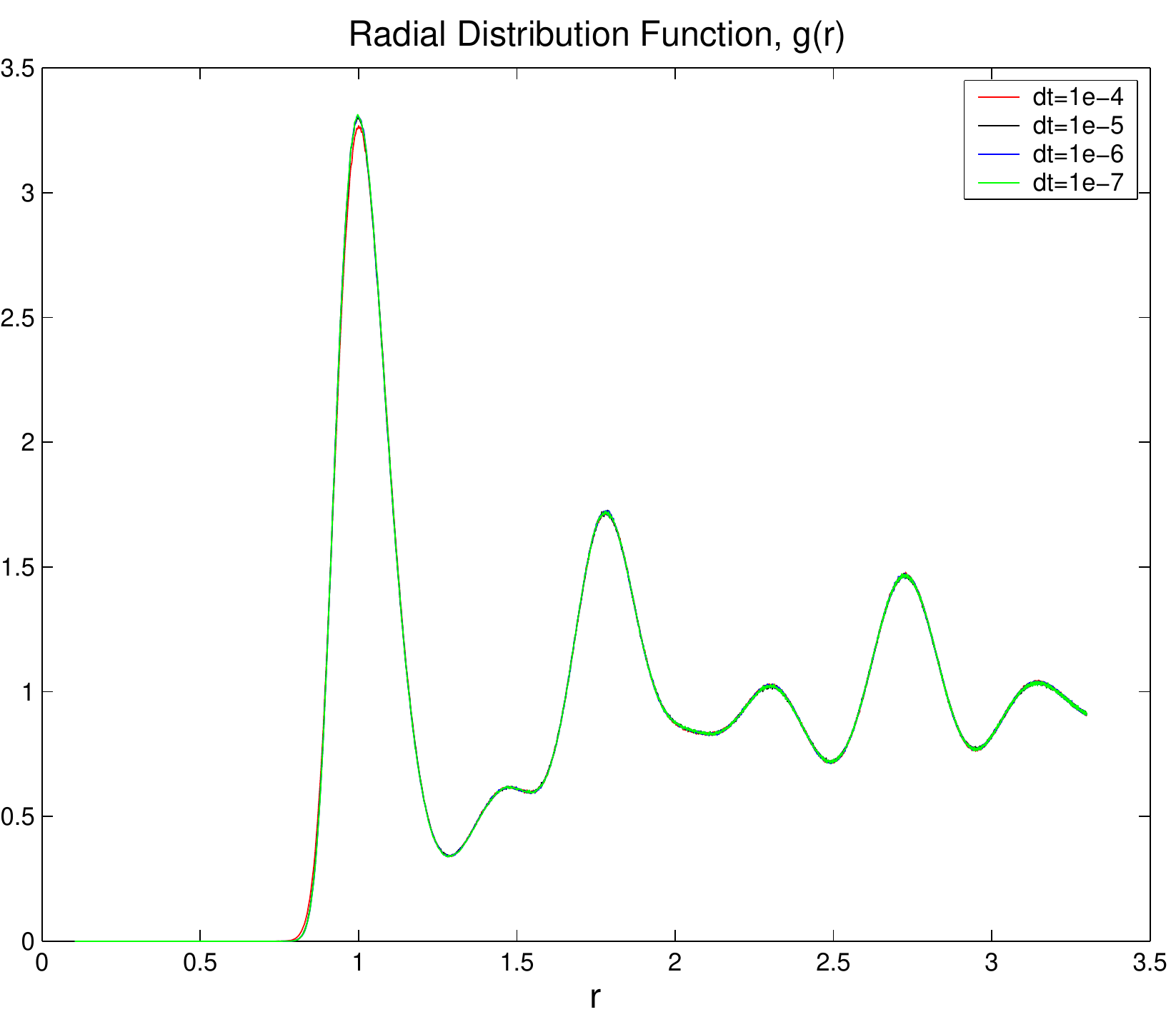}
  }
  \subfigure[Detail: first slope]{
    \label{subfig:rdf_detail}
    \includegraphics[width = 5cm] 
    {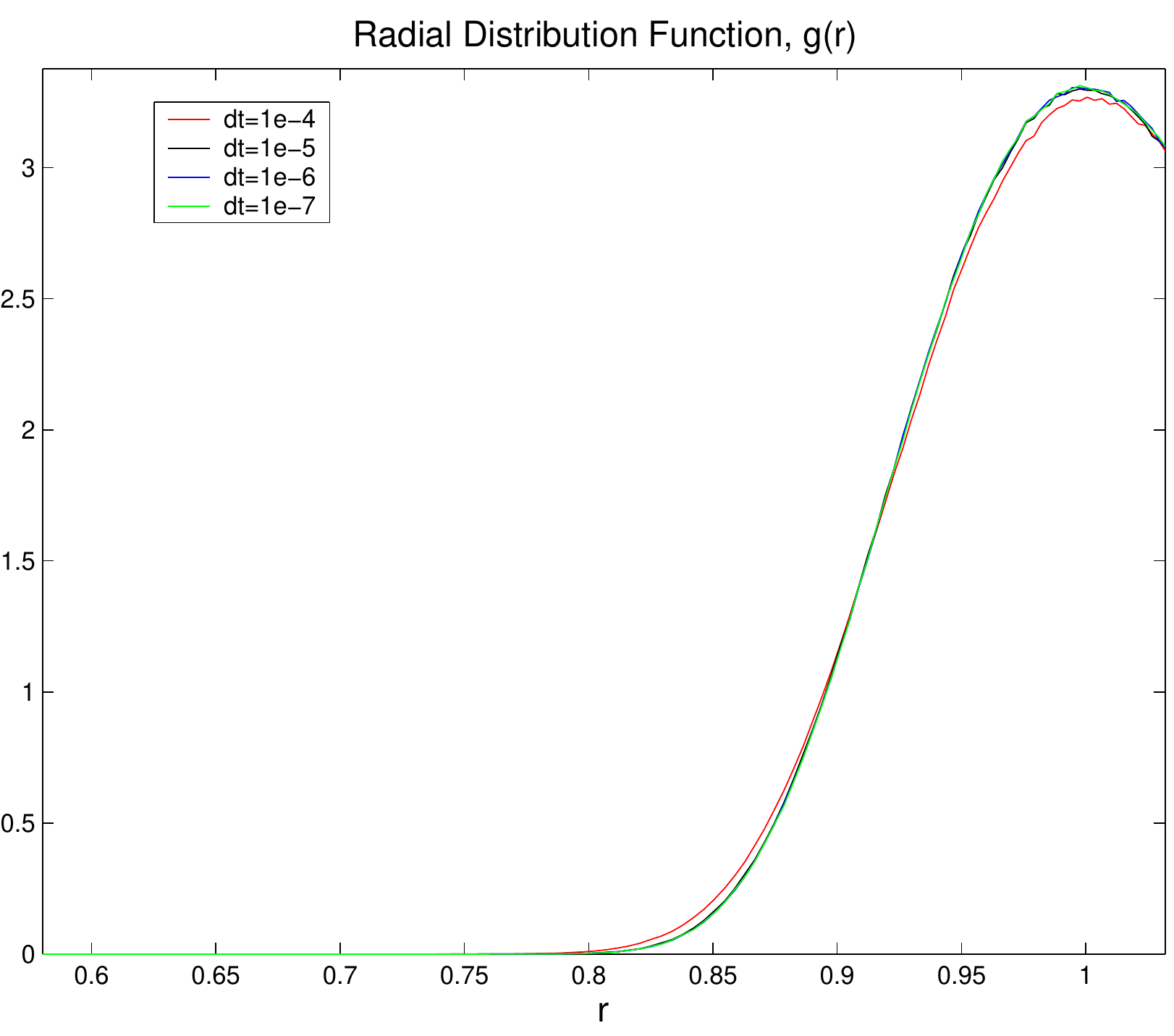}
  }
  \caption{The radial distribution function, \rdf, computed using four 
    different step sizes in a single phase FCC simulation. 
    The difference between the curves is small~(a), 
    even if the one obtained for $\Delta t=10^{-4}$ differs visibly 
    from the others in the first peak~(b). 
    In spite of the good approximation in the radial distribution
    function, the larger step sizes give very poor results in the
    computed dynamics of \pfen.}
  \label{fig:ref_dt}
\end{figure}

\begin{figure}[hbp]
  \centering
  \subfigure[$\average{\pfen}{\sampleset}$]{
    \label{subfig:pfen_dt1em5}
    \includegraphics[width=6.5cm]
    {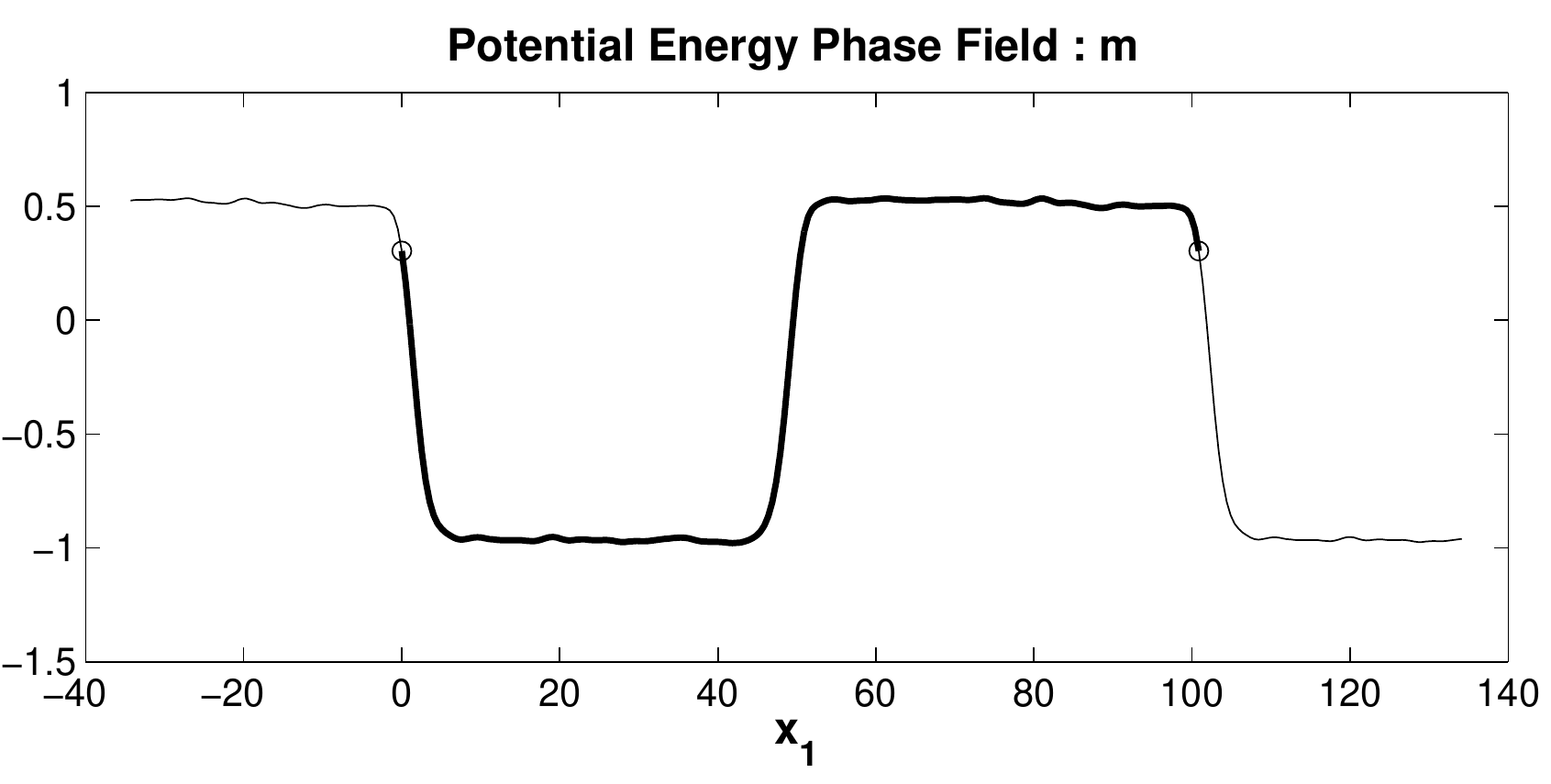}}
  \subfigure[$\average{\driftmd}{\sampleset}$]{
    \label{subfig:drift_dt1em5}
    \includegraphics[width=6.5cm]
    {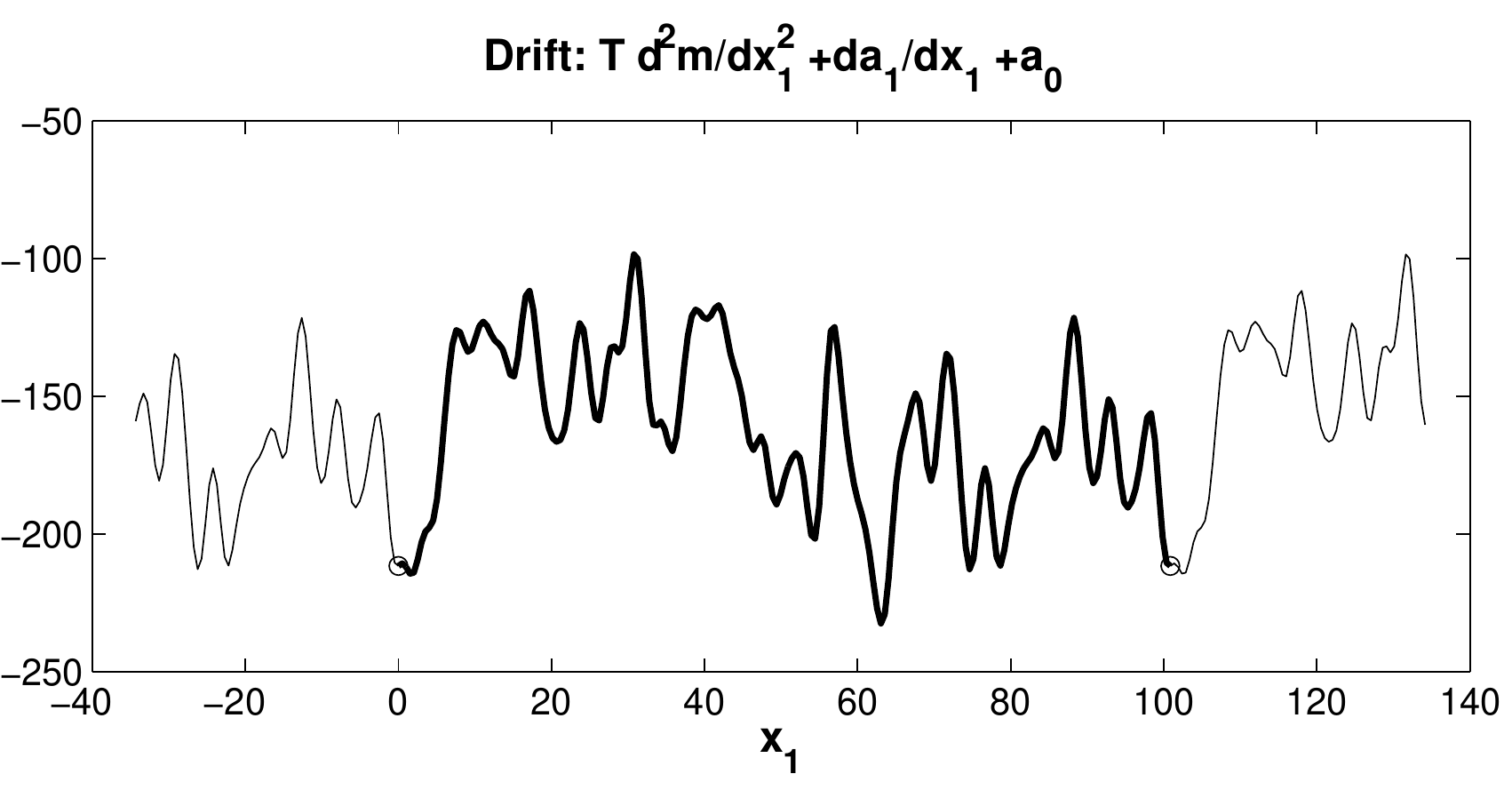}}
  \subfigure[$\average{\driftone}{\sampleset}$]{
    \label{subfig:a1_dt1em5}
    \includegraphics[width=6.5cm]
    {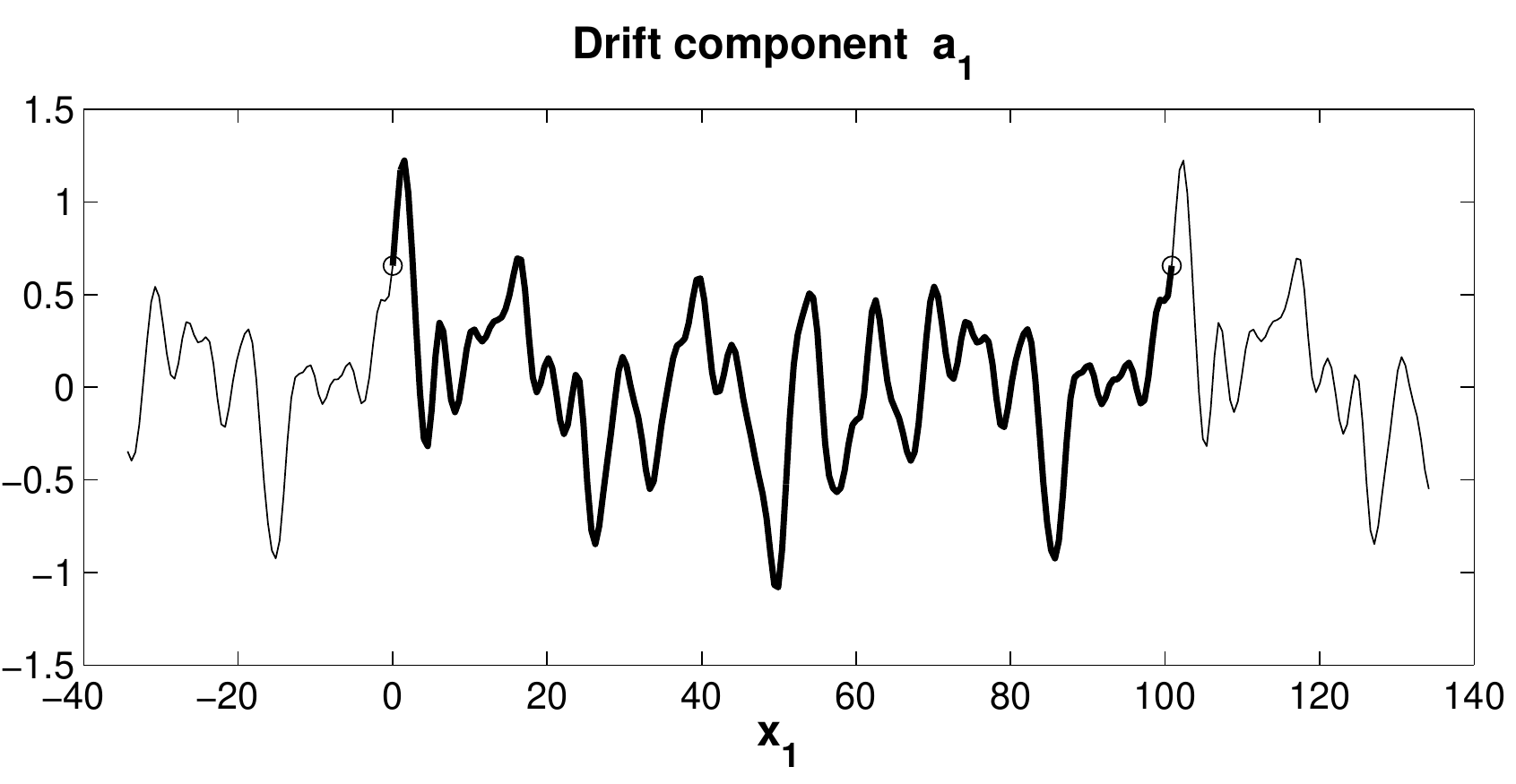}}
  \subfigure[$\average{\driftzero}{\sampleset}$]{
    \label{subfig:a0_dt1em5}
    \includegraphics[width=6.5cm]
    {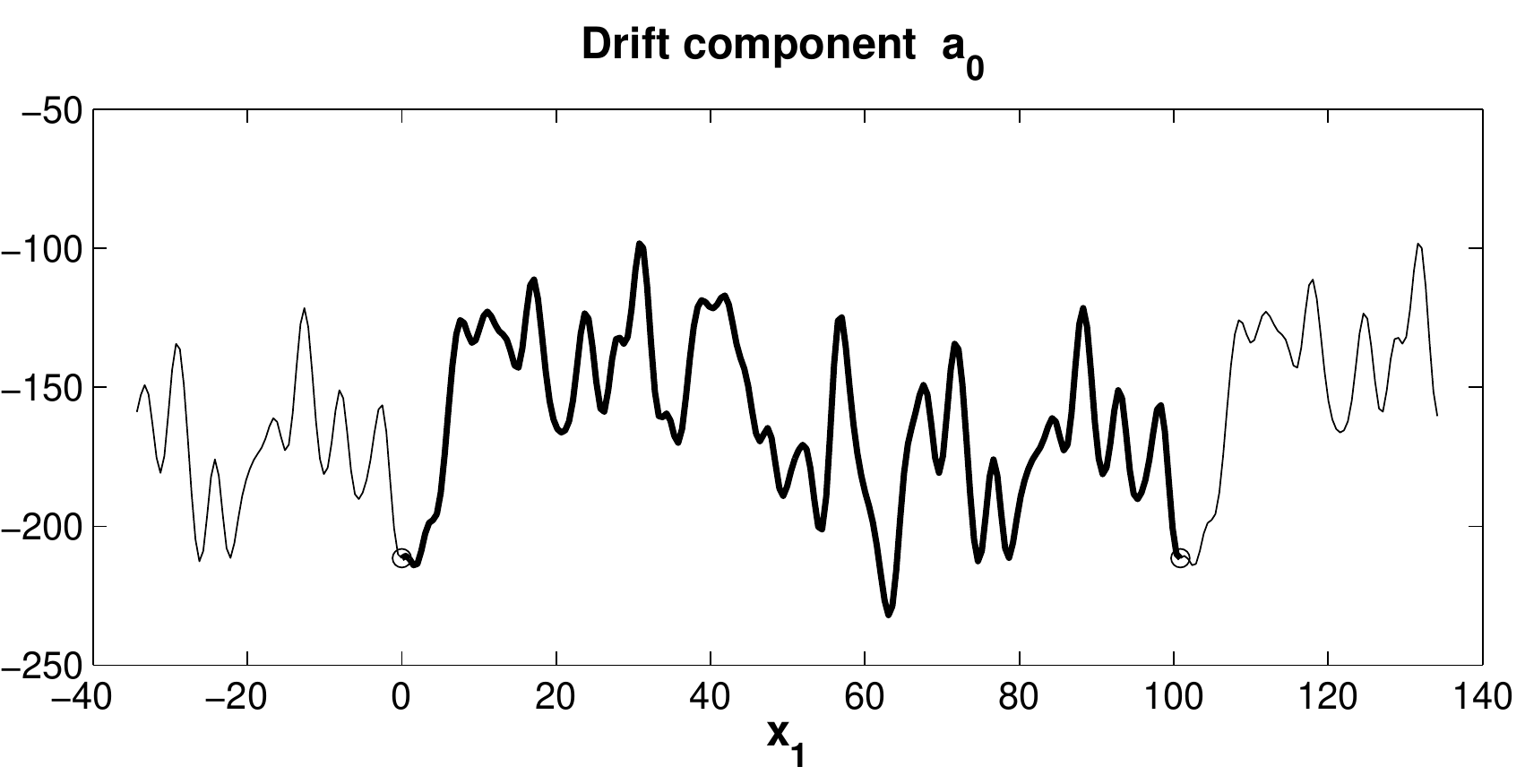}}
  \subfigure[$\kb\temp\ddxett\average{\pfen}{\sampleset}$]{
    \label{subfig:d2m_dt1em5}
    \includegraphics[width=6.5cm]
    {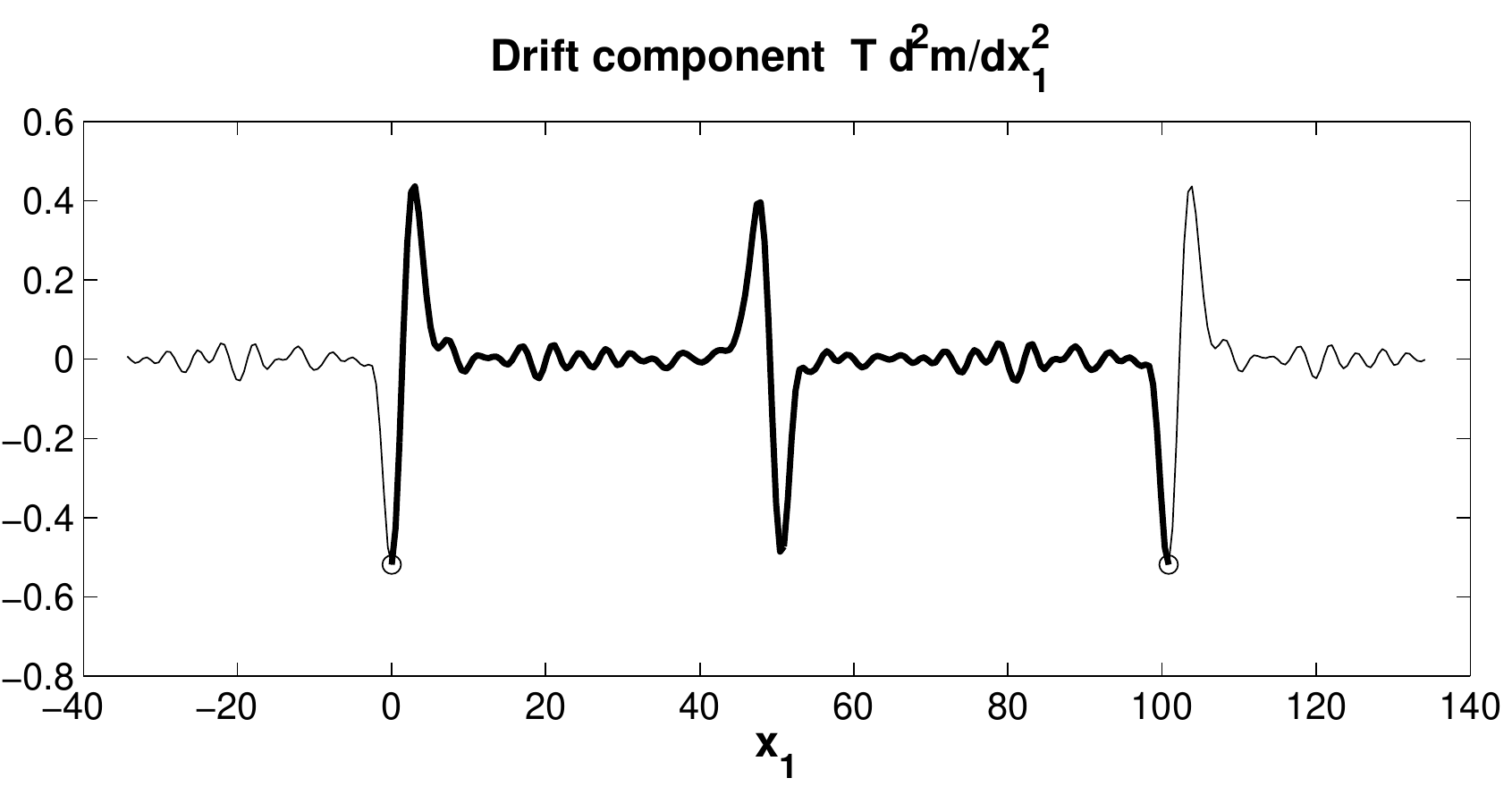}}
  \caption{Using the step size $\Delta t=1\cdot10^{-5}$ in the
    Euler-Maruyama scheme, 
    the computed average phase-field
    $\average{\pfen}{\sampleset}$ is approximately stationary
    during the time interval of the averaging. In subfigure~(a)
    the average is based on 123 configurations, sampled at every  
    ten thousandth time step, corresponding to a total time
    interval of 12.3. 
    Still, the computed average drift
    $\average{\driftmd}{\sampleset}$ is far from zero during this
    time interval.
    The large deviation from zero is entirely due to the
    term~$\average{\driftzero}{\sampleset}$.} 
  \label{fig:dt_1em5}
\end{figure}

First consider the error associated with the discrete dynamics. 
The explicit form of the drift 
is derived for the continuous time mathematical model with the
Smoluchowski dynamics~\eqref{eq:Smol_impl}, 
and not the discrete time Euler-Maruyama
dynamics~\eqref{eq:algo_Smol} that is used in the numerical
simulations. 
For a fixed size of the time step this means that, even if the
state of the numerical simulation is stationary on the time scale
of the simulation so that time averaged phase-field 
converges to an equilibrium profile, the time averaged total
drift 
will not go zero because of the time discretisation error.
Figure~\ref{fig:ref_dt} shows that the computed radial
distribution functions, here from single phase solid
configurations, are close when the time steps used vary from
$10^{-7}$ to $10^{-4}$; still the larger time steps give 
average computed drifts
$\average{\driftmd(x;\mdposd{\no})}{\sampleset}$ that are
inconsistent with the observed time evolution of the average
phase-field $\average{\pfen(x;\mdposd{\no})}{\sampleset}$. 
As shown in Figure~\ref{fig:dt_1em5}, the time step 
$\Delta t=1\cdot10^{-5}$ gives an average drift 
that oscillates between -100 and -250, even when the computed
phase-field $\average{\pfen(x;\mdposd{\no})}{\sampleset}$  is
approximately constant over times of the order 10.
For this reason, the time step used in simulations O1 and O2,
generating configurations for the computation of 
$\average{\pfen(x;\mdposd{\no})}{\sampleset}$ and
$\average{\driftmd(x;\mdposd{\no})}{\sampleset}$, was 
$\Delta t=5\cdot10^{-7}$ , while the time step used in the setup
of the initial configurations often was a thousand times larger. 
With this small time step the fluctuations in the computed
average drift outweighs the deviation from the expected zero
mean; see Figure~\ref{fig:drift_not_zero}.

The choice of the time step size $\Delta t=5\cdot10^{-7}$ was
guided by a rough error estimate, taking into account the maximal 
absolute value of second order derivatives of the Smoluchowski
drift $-\dX\totpot(\mdpos{\no})$ when the nearest neighbours
don't come closer than approximately 0.8, as indicated by
Figure~\ref{fig:ref_dt}. 
Then the time step was adjusted so that the slow convergence of
the time averaged drift in terms of \tend\ and the number of
configurations, \mdposd{\no}, was the dominating error source in
the results.
This over-killing of the time discretisation error in the
molecular dynamics wastes computer power and could possibly be
avoided by more accurate error estimates, allowing a matching of
the different error contributions.
Using a reasonable number of grid points, $K$, in the computation
of the drift coefficient $K$-vectors and the diffusion
 $K$-by-$K$ matrix \diffumat, in~\eqref{eq:diffumat}, the
computational cost for obtaining \diffumat\ in particular, far
exceeds the cost of actually making a time step in the molecular
dynamics simulation. 
Hence the additional cost of over-killing the time step error is
not very significant, provided that not every configuration
in the time stepping is included in the averages
\average{\pfen}{\sampleset}, \average{\driftmd}{\sampleset},
and \diffumat. 
In the averages shown in Figure~\ref{fig:phasefield_not_eq} and
Figure~\ref{fig:drift_not_zero}, for example, the configurations
were sampled at time intervals $5\cdot10^{-4}$,
corresponding to 1000 time steps in the molecular dynamics
simulation. 

\begin{figure}[hbp]
  \centering
  \includegraphics[width=7cm]{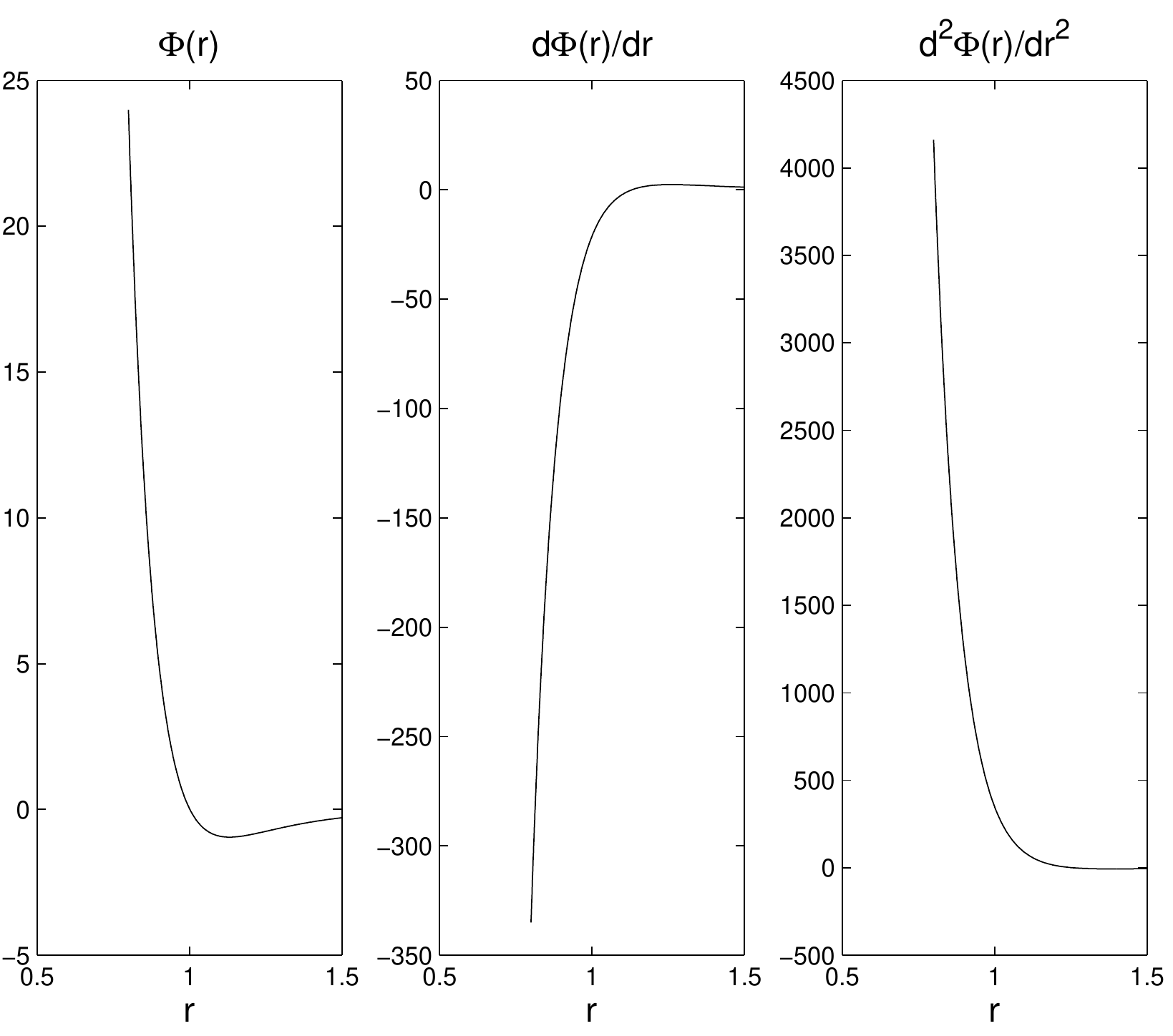}
  \caption{The absolute value of the Exp-6 potential and its
    derivatives grow very quickly with decreasing $r$, in
    the range with positive \rdf\ in Figure~\ref{subfig:rdf_detail}.
    The potential and its two first derivatives using
    the model parameters in Table~\ref{tab:LJunits}, on
    page~\pageref{tab:LJunits}, are shown here.}
  \label{fig:exp6_with_derivatives}
\end{figure}

A further improvement may be to incorporate finite step-size
effects in the expressions for the components of the drift.
The higher order derivatives of the pair potential attain large
values when two particles come closer than 1; see
Figure~\ref{fig:exp6_with_derivatives}. 
Hence the time step must be taken very small for \Ito's formula
to be a good approximation of the dynamics of the discrete
system. 
Instead of a direct application of \Ito's formula in the
derivation of the drift and diffusion terms
in~\eqref{eq_app:driftmd_j} and~\eqref{eq_app:diffumd_j} on
page~\pageref{eq_app:driftmd_j} one could include higher
order terms in the expansion to improve the accuracy of the
computed drift. 

\begin{figure}[hbp]
  \centering
  \subfigure[Mean based on 111 configurations, $\tend=0.0555$]
  {\label{subfig:drift_111}
    \includegraphics[width=6.5cm]
    {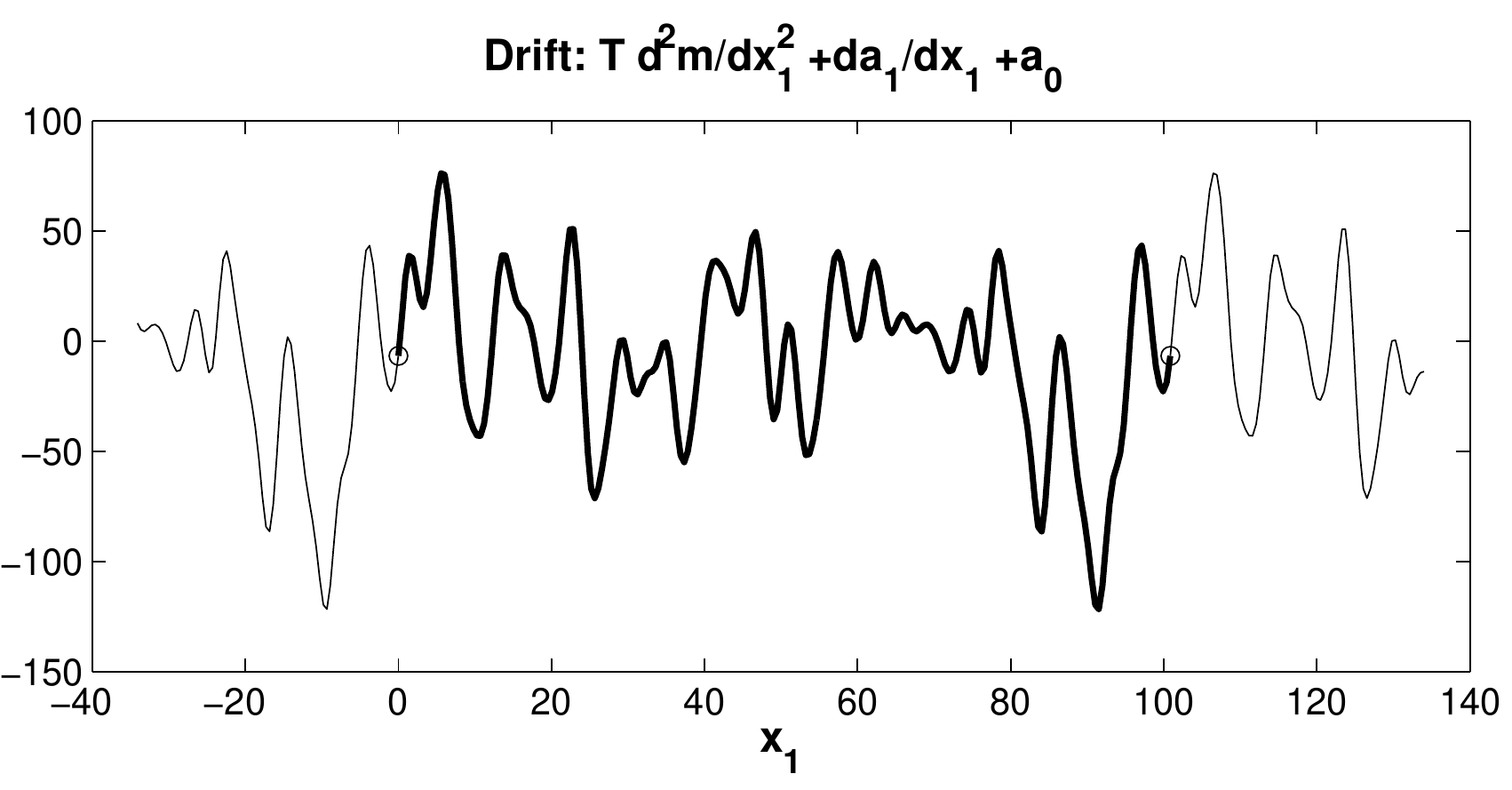}}
  \subfigure[Mean based on 444 configurations, $\tend=0.2220$]
  {\label{subfig:drfit_444_dense}
    \includegraphics[width=6.5cm]
    {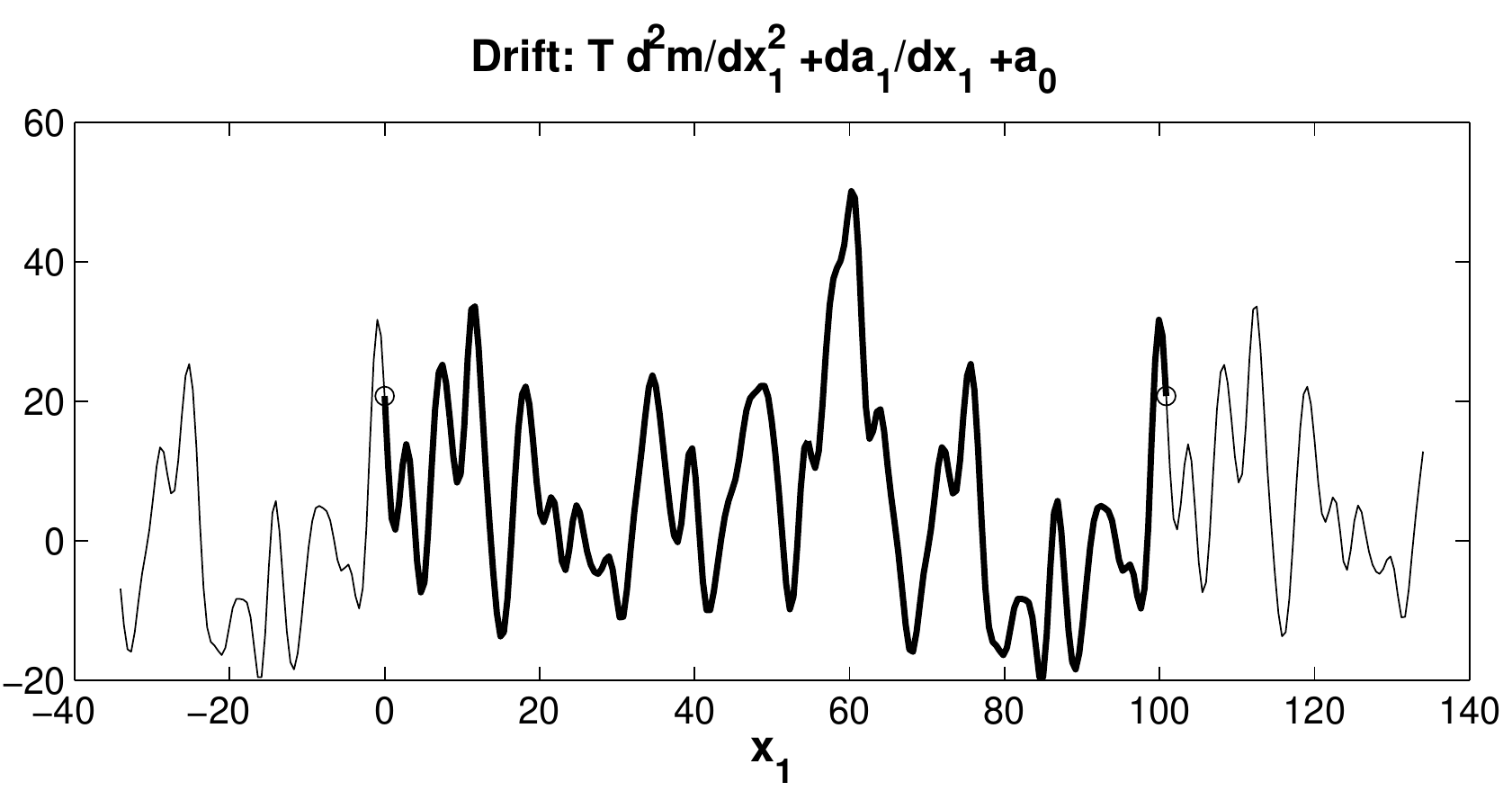}}
  \subfigure[Mean based on 1775 configurations, $\tend=0.8875$]
  {\label{subfig:drift_1775}
    \includegraphics[width=6.5cm]
    {figures/drift_O2_eps_1_0_1775_cfgs.pdf}}
  \subfigure[The spatial ($x_1$) mean and, within parentheses,
  variance of $\average{\driftmd}{\sampleset}$]{
    \begin{tabular}{|l|c|c|c|}
      \hline
      & $\tend=0.0555$ & $\tend=0.2220$ & $\tend=0.8750$  \\
      \hline
      111 cfgs. & -5.7 ($1.3\cdot10^3$) 
      & -10.3 ($1.2\cdot10^3$) & 4.3 ($1.1\cdot10^3$) \\  
      \hline
      444 cfgs. &  & 6.1 ($2.1\cdot10^2$) &  0.67 ($2.7\cdot10^2$) \\
      \hline
      1775 cfgs. &  &  & 1.9 ($3.8\cdot10^1$) \\
      \hline
    \end{tabular}
  }
  \caption{The total drift $\average{\driftmd}{\sampleset}$,
    decays slightly faster with \tend\ than the predicted 
    $1/\sqrt{\tend}$ in the examples (a), (b), and (c) above.
    Here the number of configurations in the averages grows with
    \tend\ and the means and variances of
    $\average{\driftmd}{\sampleset}$ tabulated in~(d) suggest
    that the number of configurations still restricts the rate of
    convergence. 
    The average in subfigure~(c) is based on the same 1775
    configurations from simulation O2 as
    $\average{\pfen}{\sampleset}$ in Figure~\ref{subfig:pf_O2}.
    The averages in subfigures~(a) and~(b) are based on the first
    111 and 444 configurations, respectively.
  }
  \label{fig:conv_total_drift}
\end{figure}

\subsubsection{Dependence on the length of the time averaging
    interval}

Next consider the dependence of the computed coarse-grained drift
coefficient function on the length of the time interval \tend.
Introducing the time averaged drift over a sample path as
\begin{align*}
  \driftcgxs_\tend & = 
  \frac{1}{\tend} \int_0^\tend \driftmd(\cdot;\mdpos{\no})\;dt, 
\end{align*}
the coarse-grained drift~\eqref{eq:def_driftcgx} is 
$\driftcgx=\E[\driftcgxs_\tend]$. 
The rate of convergence of $\driftcgx$, as $\tend\to\infty$,
in the continuous time mathematical model can be estimated by
integration of the stochastic differential
equation~\eqref{eq:sde_pfen} for the phase-field \pfen.
Integrating from 0 to \tend\ 
gives 
\begin{align}
  \label{eq:integrated_sde}
  \pfen(\cdot;\mdpos[\tend]{\no}) - \pfen(\cdot;\mdpos[0]{\no}) 
  & = 
  \int_0^\tend \driftmd(\cdot;\mdpos{\no})\;dt
  + 
  \int_0^\tend 
  \sumall{j} \sum_{k=1}^3\diffumd_{j,k}(\cdot;\mdpos{\no})\;dW_{j,k}^t,
\end{align}
so that, by taking the expectation and using that, since
\mdpos{\no} is $W^t$-adapted, the expectations of the
\Ito-integrals vanish 
\begin{align}
  \label{eq:meandrift}
  \E\left[
    \int_0^\tend \driftmd(\cdot;\mdpos{\no})\;dt 
  \right]
  & = 
  \E\left[
    \pfen(\cdot;\mdpos[\tend]{\no}) - \pfen(\cdot;\mdpos[0]{\no}) 
  \right].
\end{align}
Hence, if the phase-field is stationary, then the
expected mean drift over time is zero. 
Normalising~\eqref{eq:integrated_sde} and~\eqref{eq:meandrift} by
\tend, 
\begin{multline*}
  \driftcgxs_\tend - 
  \E\left[\driftcgxs_\tend\right]
  = \\
  \frac{1}{\tend}
  \left( 
    \pfen(\cdot;\mdpos[\tend]{\no}) - 
    \pfen(\cdot;\mdpos[0]{\no})
    - 
    \E\left[
      \pfen(\cdot;\mdpos[\tend]{\no}) - 
      \pfen(\cdot;\mdpos[0]{\no})
    \right]
    -
    \sumall{j} \sum_{k=1}^3
    \int_0^\tend 
    \diffumd_{j,k}(\cdot;\mdpos{\no})\;dW_{j,k}^t
  \right)
\end{multline*}
and the variance of
$\driftcgxs_\tend$ is obtained as
\begin{align*}
  \Var[\driftcgxs_\tend]
  & =
  \E\left[
    \left( \driftcgxs_\tend - 
      \E\left[\driftcgxs_\tend\right]
    \right)^2
  \right] 
  \\
  & =
  \frac{1}{\tend^2}\Var\biggl[ 
    \pfen(\cdot;\mdpos[\tend]{\no}) - 
    \pfen(\cdot;\mdpos[0]{\no})
  \biggr]
  \\
  & \quad +
  \frac{1}{\tend^2}
  \sumall{j} \sum_{k=1}^3
  \E\left[
    \left( 
      \int_0^\tend 
      \diffumd_{j,k}(\cdot;\mdpos{\no})\;dW_{j,k}^t,
    \right)^2
  \right] 
  \\
  & \quad -
  \frac{2}{\tend^2}
  \E\left[
    \bigl( 
      \pfen(\cdot;\mdpos[\tend]{\no}) - 
      \pfen(\cdot;\mdpos[0]{\no})
    \bigr)  
    \left(
      \sumall{j} \sum_{k=1}^3
      \int_0^\tend 
      \diffumd_{j,k}(\cdot;\mdpos{\no})\;dW_{j,k}^t
    \right)  
  \right],
\end{align*}
where last expression was simplified using the independence of the
different components of $W^t$, and the zero expected value of
\Ito\ integrals.
Assuming that both the phase-field and all the diffusion
coefficients are bounded, the dominating term in the expression
for the variance is 
\begin{align*}
  \frac{1}{\tend^2}
  \sumall{j} \sum_{k=1}^3
  \E\left[
    \left( 
      \int_0^\tend 
      \diffumd_{j,k}(\cdot;\mdpos{\no})\;dW_{j,k}^t,
    \right)^2
  \right] 
  & = \mathcal{O}\left(\frac{1}{\tend}\right).
\end{align*}
In the two phase simulations considered here, the values of the
computed phase-field varies between a lower level in the solid
a higher in the liquid.
Because of the small positive probability for two particles, with 
trajectories computed using the Euler-Maruyama
dynamics~\eqref{eq:algo_Smol}, to get within an arbitrarily small
distance of each other, there is no guarantee that computed
phase-field always will stay in this range.
However, if the minimum interatomic distance becomes to small,
that is a breakdown of the whole microscopic model and not just a
problem when computing the drift; this situation has not been
observed to happen in the simulations here and the observed
values of the phase-field are all in the range $(-1.5,1.0)$. 
Hence the assumption that \pfen\ is bounded seems reasonable
here; a bound on the absolute value of the diffusion coefficients
$\diffumd_{j,k}$ is less certain, and it will have to be larger
than the bound on \pfen.

For the average drift to be small compared to the
stationary values of the phase-field itself, it must be at least
a factor 100 smaller than the computed average shown in
Figure~\ref{fig:drift_not_zero}. 
Based on the rough analysis above, the expected time average of
the total drift can be expected to decay as $1/\sqrt{\tend}$ with
a large constant factor. 
When the computed drift \average{\driftmd}{\sampleset} in
Figure~\ref{fig:drift_not_zero} is compared to averages computed
using two smaller subsequences of configurations, the convergence
to zero appears to be slightly faster than $1/\sqrt{\tend}$; see
Figure~\ref{fig:conv_total_drift}. 
Even when extrapolating with the measured convergence rate,
decreasing the average drift by a factor 100 would require
increasing the averaging time interval by more than a factor
1000, which is beyond reach within the present project.
With increasing accuracy in the time average, eventually 
the time step in the molecular dynamics simulations must be
decreased, further increasing the computational cost.

Since the total drift coefficient function,
$\driftcgx(x_1)\approx\average{\driftmd(x_1;\cdot)}{\sampleset}$,
where 
\begin{align}
  \label{eq:meandriftmd}
  \average{\driftmd(x_1;\cdot)}{\sampleset} & =
  \kb\temp \ddxett \average{\pfen(x_1;\cdot)}{\sampleset}
  + \dxett \average{\driftone(x_1;\cdot)}{\sampleset} 
  + \average{\driftzero(x_1;\cdot)}{\sampleset},
\end{align}
in the coarse grained
model is expected to be zero in a stationary situation, a more
accurate computation would serve primarily as a consistency test. 
On the other hand, the individual terms in the right hand side
are not all expected to vanish independently. Indeed, it is clear
from the results on \average{\pfen(x_1;\cdot)}{\sampleset} in
Section~\ref{sec:computed_phase-field} 
that the term with two differentiations with respect to $x_1$
will not be identically zero. 
This also shows that while the total drift is far from 
$\average{\driftmd(x_1;\cdot)}{\sampleset}$ converged, at least
one term is reasonably accurate.

\begin{figure}[hbp]
  \centering
  \subfigure[Mean based on 111 configurations, $\tend=0.8875$]
  {\label{subfig:no_diff_111}
    \includegraphics[width=6.5cm]
    {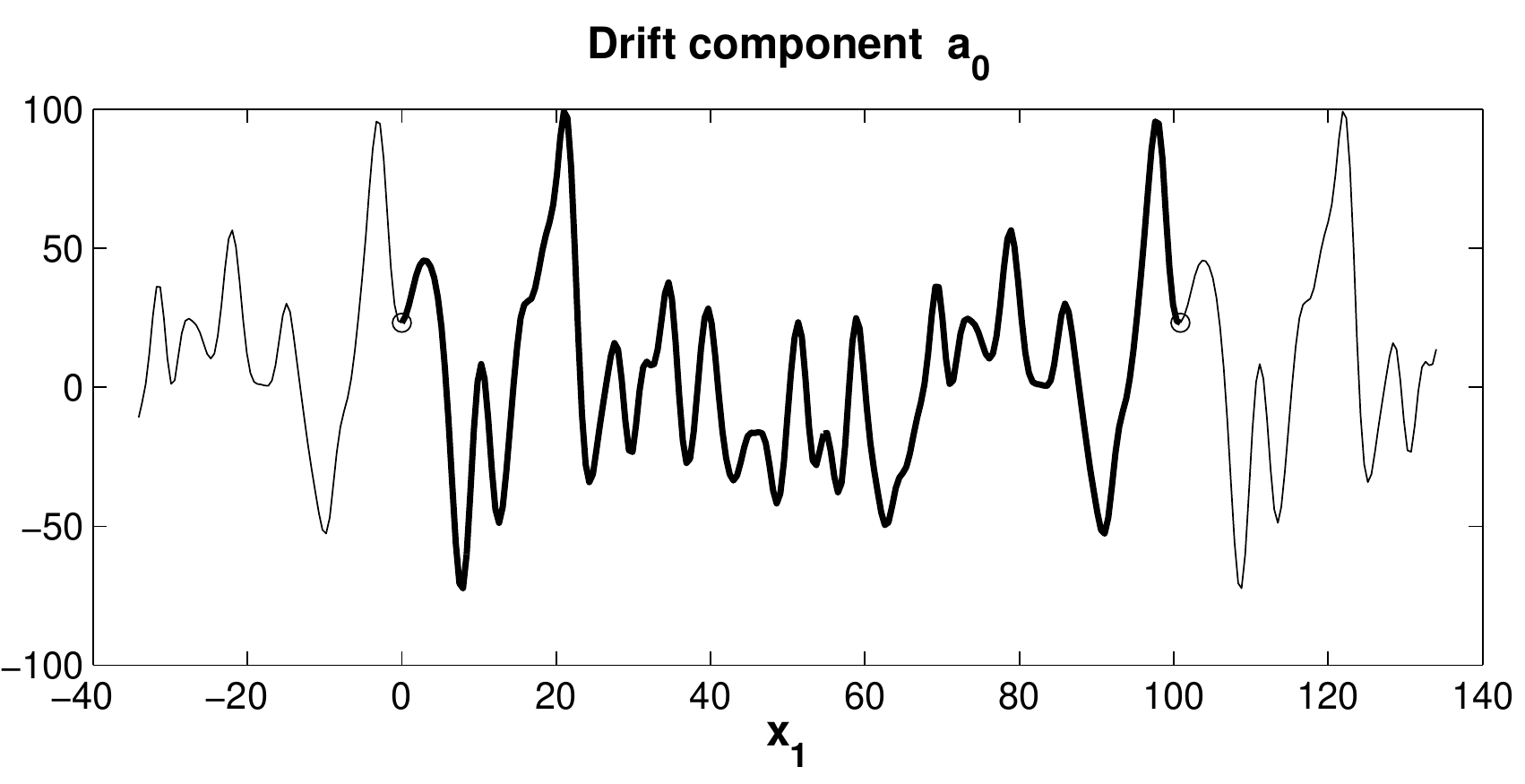}}
  \subfigure[Mean based on 444 configurations, $\tend=0.8875$]
  {\label{subfig:no_diff_444}
    \includegraphics[width=6.5cm]
    {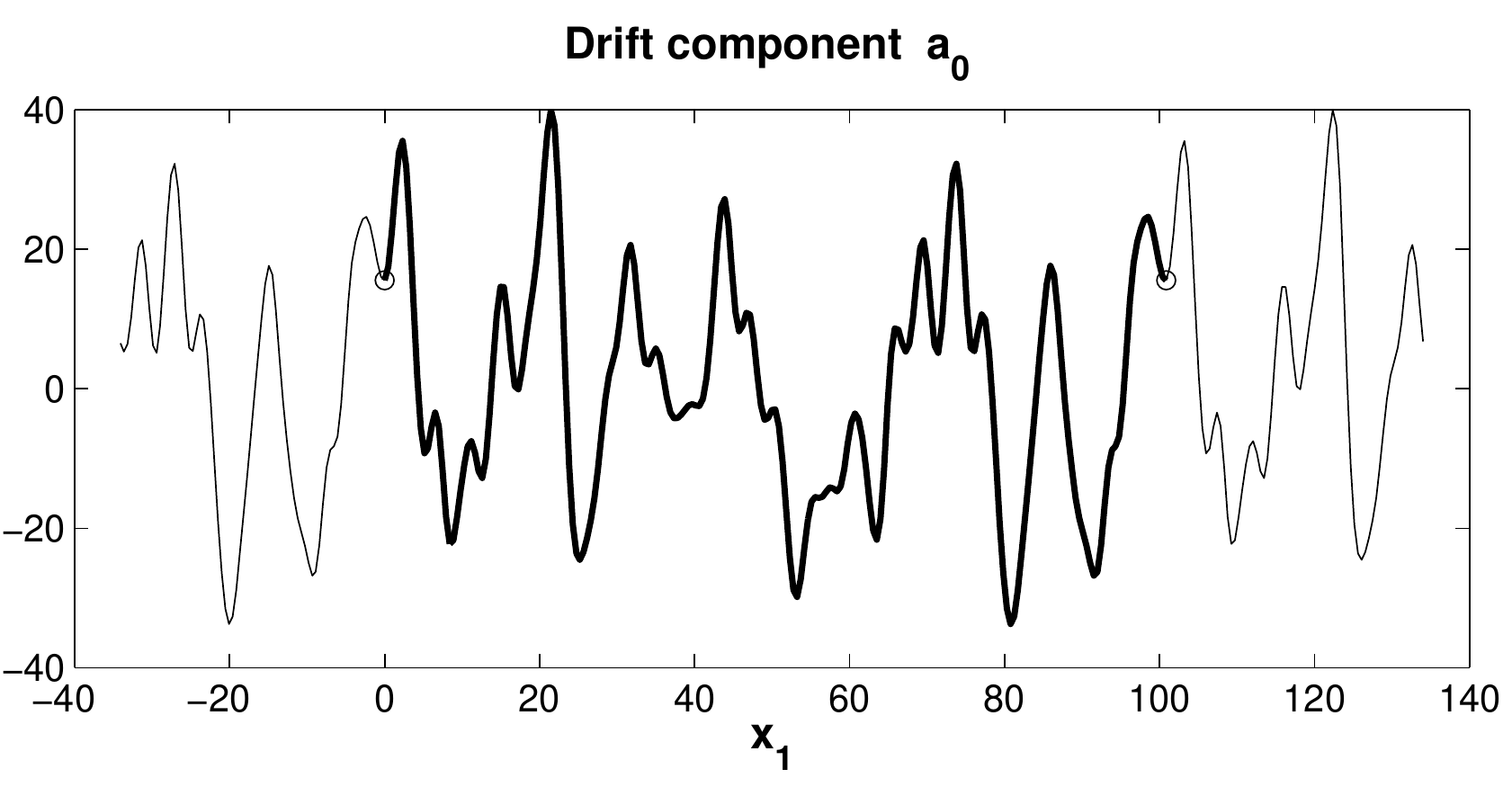}}
  \subfigure[Mean based on 444 configurations, $\tend=0.2220$]
  {\label{subfig:no_diff_444_dense}
    \includegraphics[width=6.5cm]
    {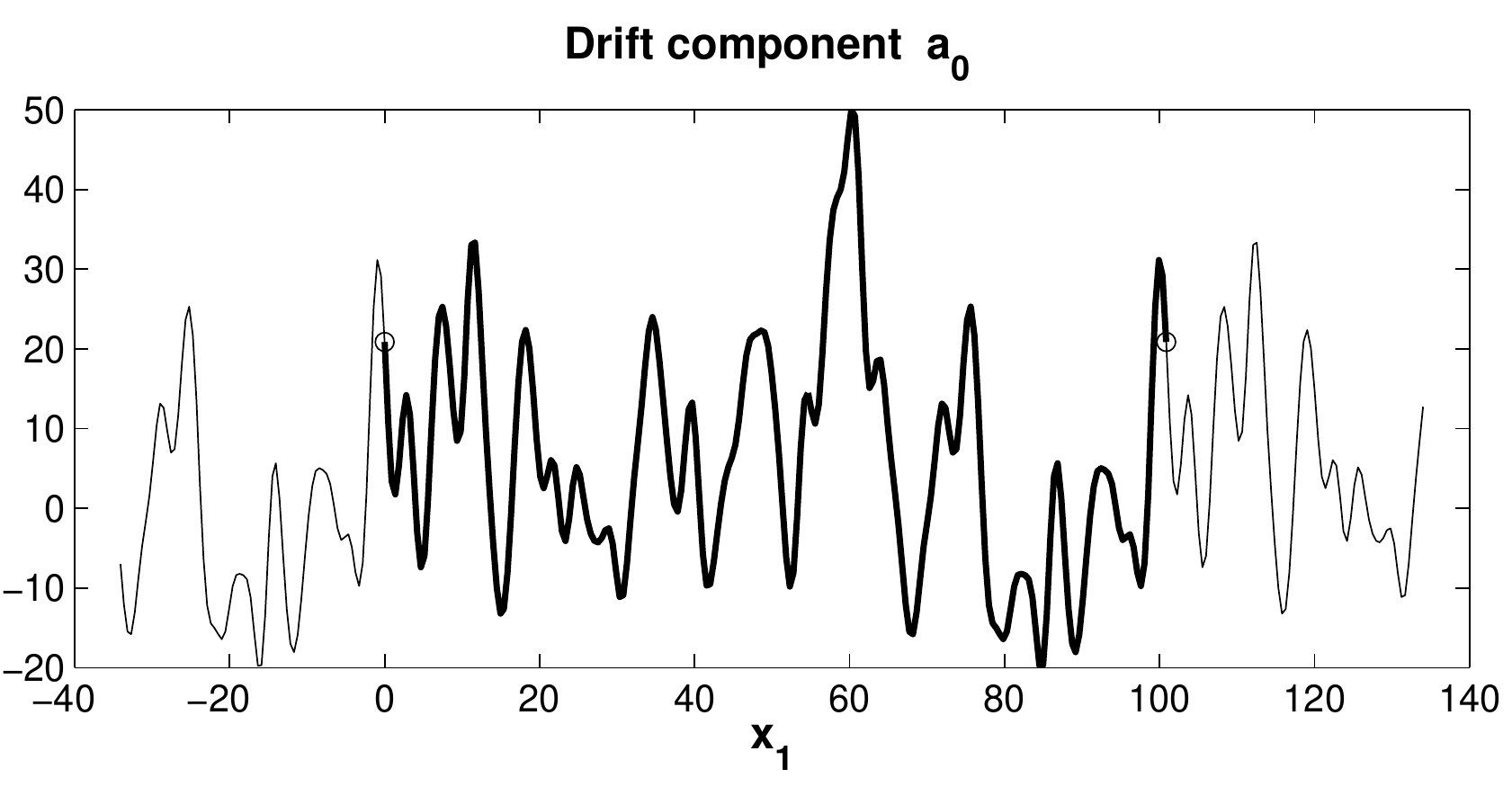}}
  \subfigure[Mean based on 1775 configurations, $\tend=0.8875$]
  {\label{subfig:no_diff_1775}
    \includegraphics[width=6.5cm]
    {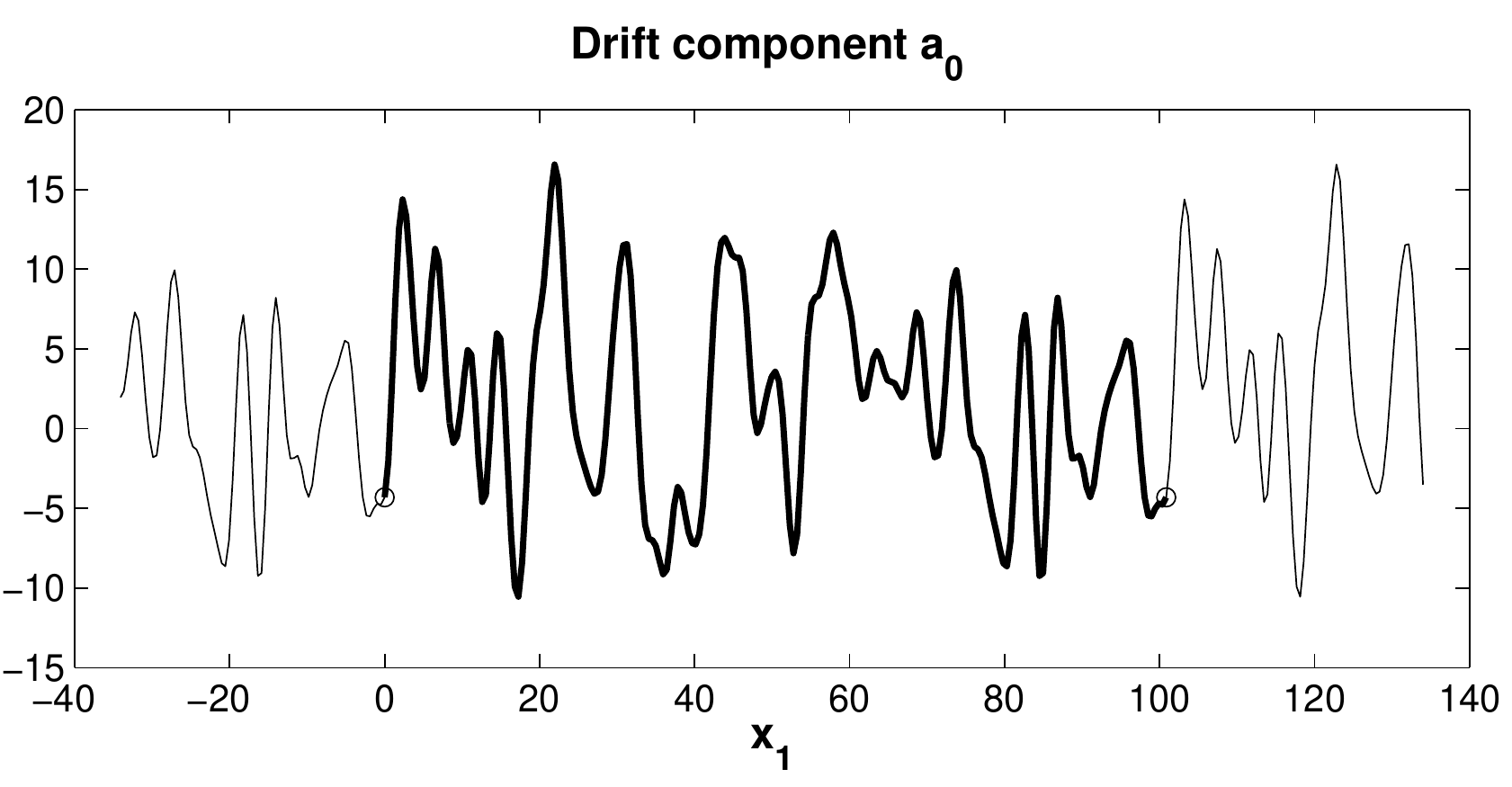}}
  \subfigure[One single configuration, \mdpos{\no}]{ 
    \label{subfig:instant_no_diff}
    \includegraphics[width=6.5cm]
    {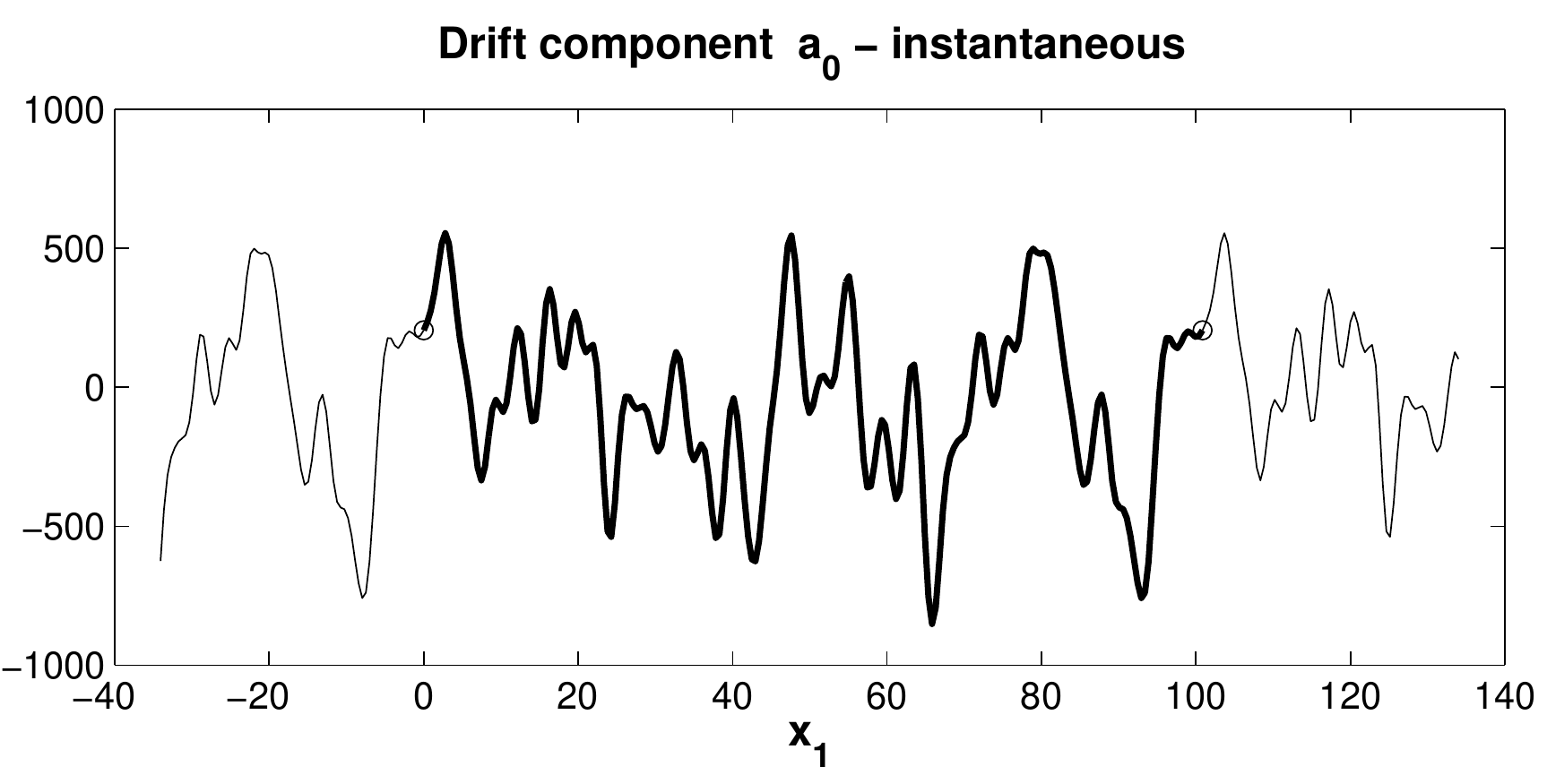}}
  \subfigure[Variance based on the configurations in (d)]{ 
    \label{subfig:var_no_diff}
    \includegraphics[width=6.5cm]
    {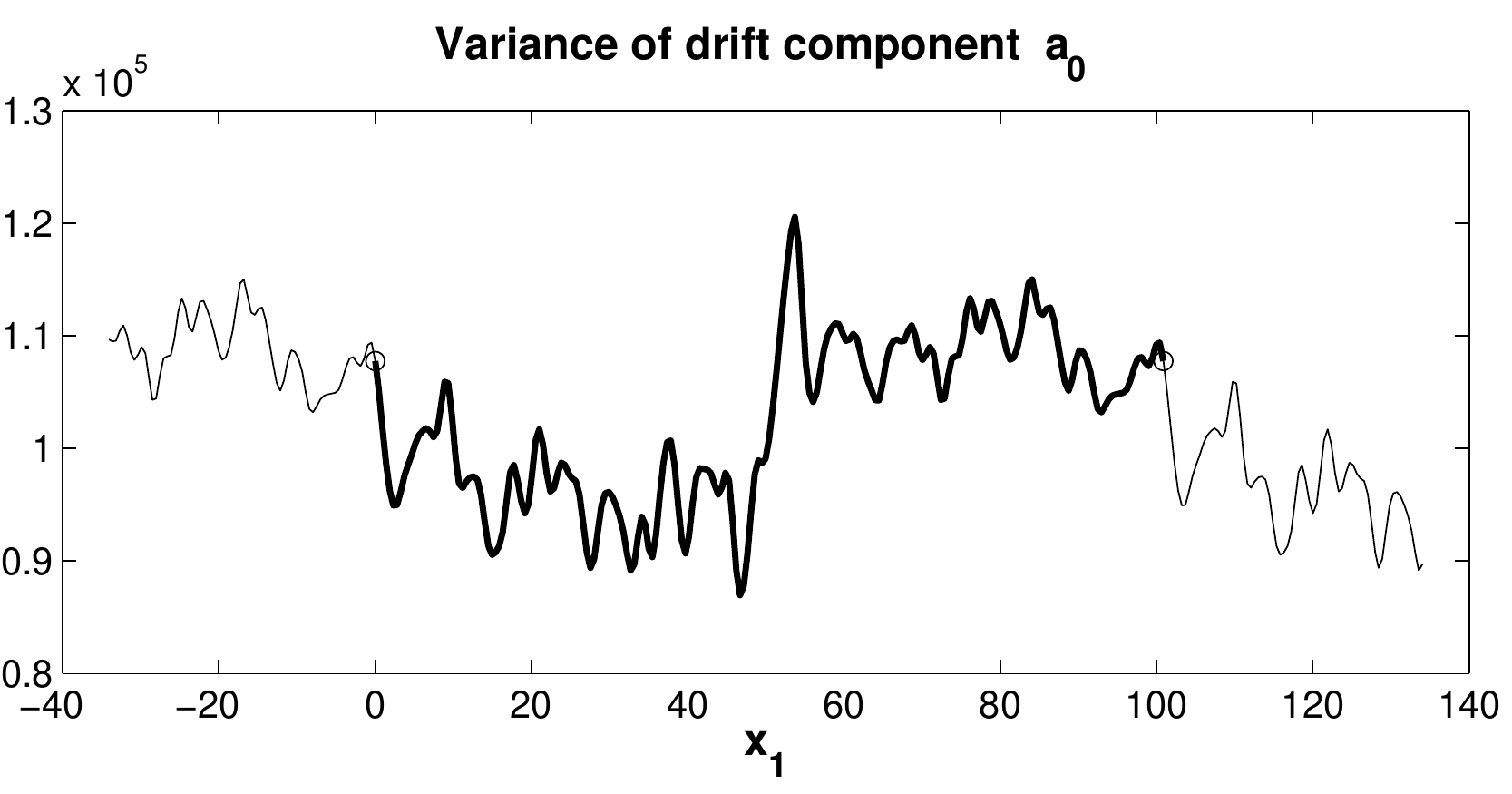}}
  \caption{The term $\average{\driftzero}{\sampleset}$ is the 
    slowest converging average in the drift average; a comparison
    with Figure~\ref{fig:conv_total_drift} shows that this term
    dominates the total drift average. 
    This explicit form of the term, given in~\eqref{eq:driftzero}
    is a sum over all particles of terms that are second order in
    the particle forces and a term containing the divergence of
    the particle force; in the molecular dynamics simulation,
    these terms are large and so is the function \driftzero, when
    computed from a single configuration, as in~(e).
    Eventually the average must decrease to order 1 through
    cancellation, but for the number of configurations available
    here fluctuations dominate the computed averages
    $\average{\driftzero}{\sampleset}$.}
  \label{fig:conv_drift_zero}
\end{figure}

A closer look on the terms of the drift, reveals that the
different terms are of different orders of magnitude.
The term \average{\driftzero(x_1;\cdot)}{\sampleset}, 
with \driftzero\ defined in~\eqref{eq:driftzero}, contains 
both second order differentials of the potential with respect to
the particle positions and second powers of first order
differentials. These terms, as illustrated in 
Figure~\ref{fig:exp6_with_derivatives}, attain much larger 
values than the potential itself and cancellation is required to
reduce \average{\driftzero(x_1;\cdot)}{\sampleset}\ to a
size comparable with the two other terms in the drift. 
Figure~\ref{subfig:instant_no_diff} shows an individual
$\driftzero(x_1;\cdot)$ computed from one configuration; in
the length of the computational cell, the values range from
approximately -500 to +500, whereas the phase-field,
$\pfen(x_1;\cdot)$, is of the order 1, and $\driftone(x_1;\cdot)$
is of intermediate magnitude.
A comparison between the computed averages
$\average{\driftmd(x_1;\cdot)}{\sampleset}$ in
Figure~\ref{fig:conv_total_drift} and 
\average{\driftzero(x_1;\cdot)}{\sampleset} in
Figure~\ref{fig:conv_drift_zero} shows that 
\average{\driftzero(x_1;\cdot)}{\sampleset}\ is the dominates the
other two terms completely here. 

\begin{figure}[hbp]
  \centering
  \subfigure[Mean based on 111 configurations, $\tend=0.8875$]
  {\label{subfig:one_diff_111}
    \includegraphics[width=6.5cm]
    {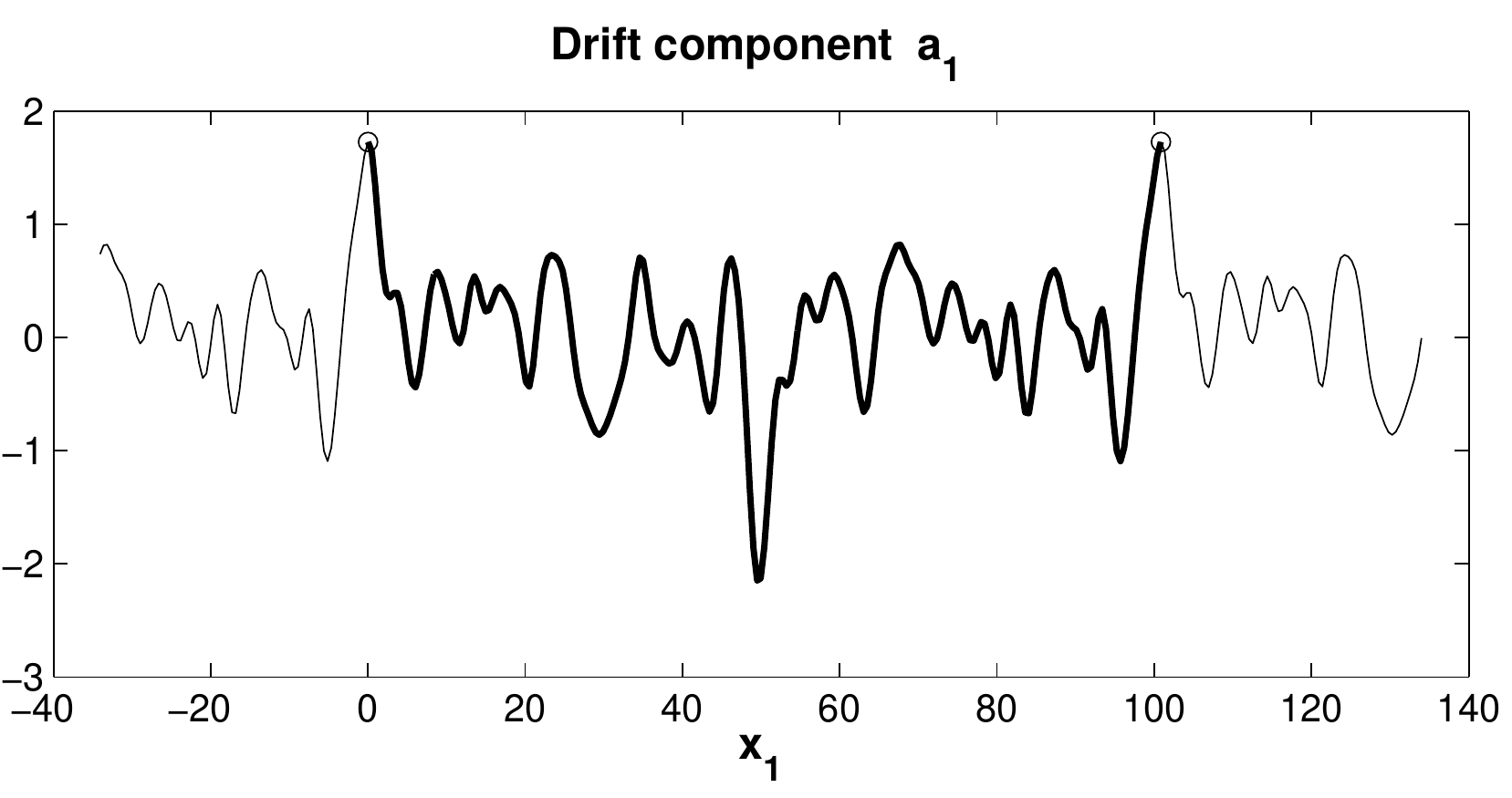}}
  \subfigure[Mean based on 444 configurations, $\tend=0.8875$]
  {\label{subfig:one_diff_444}
    \includegraphics[width=6.5cm]
    {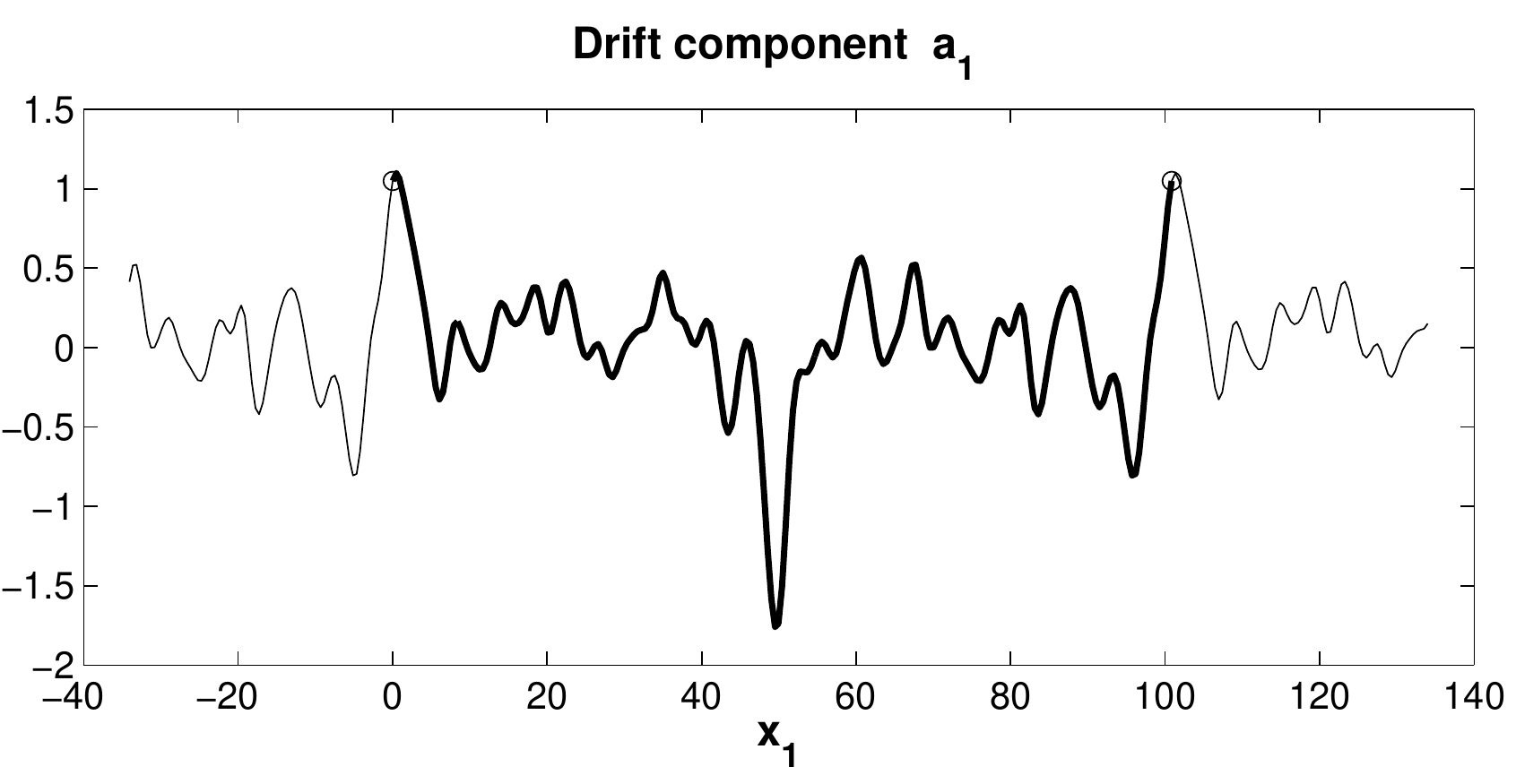}}
  \subfigure[Mean based on 444 configurations, $\tend=0.2220$]
  {\label{subfig:one_diff_444_dense}
    \includegraphics[width=6.5cm]
    {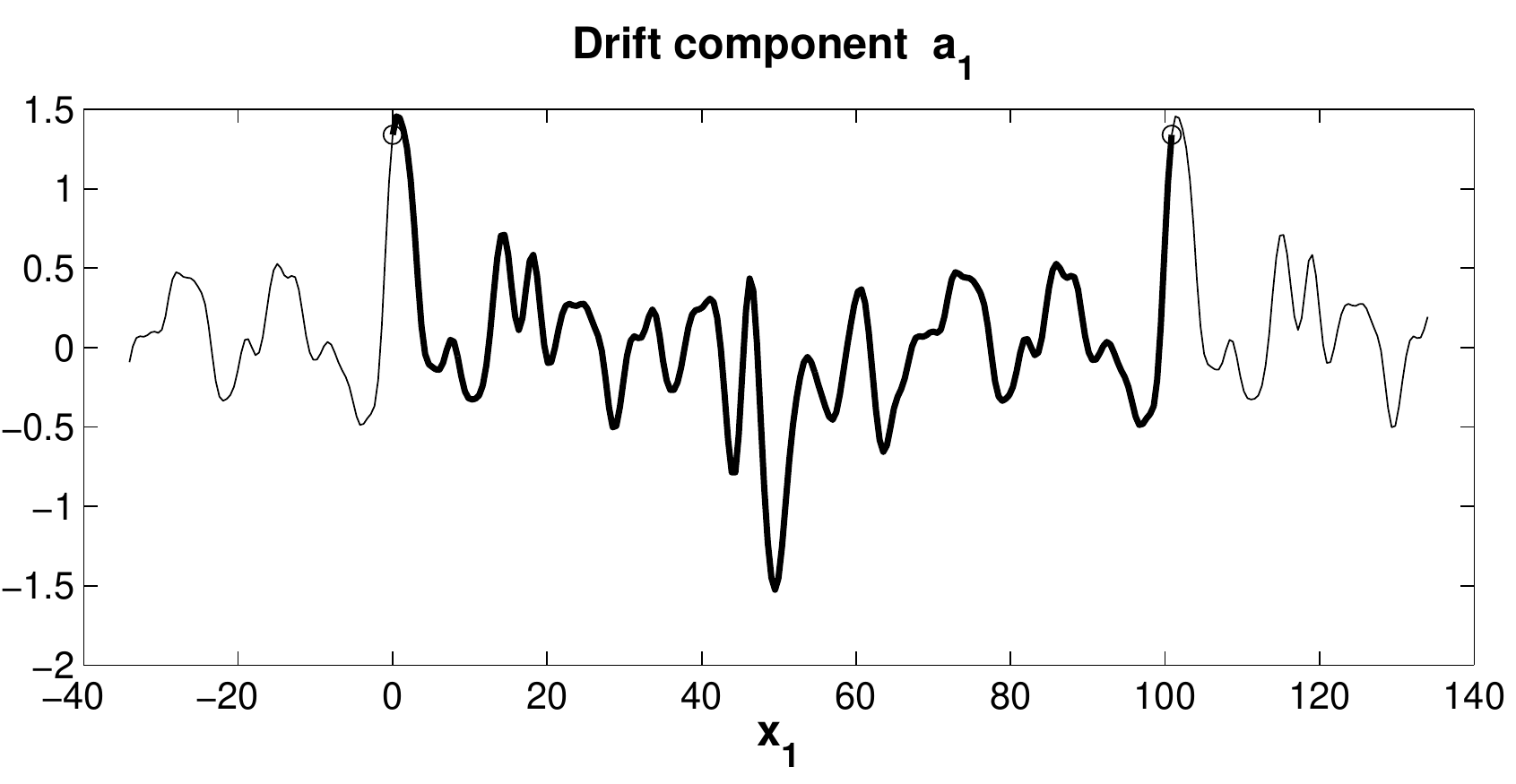}}
  \subfigure[Mean based on 1775 configurations, $\tend=0.8875$]
  {\label{subfig:one_diff_1775}
    \includegraphics[width=6.5cm]
    {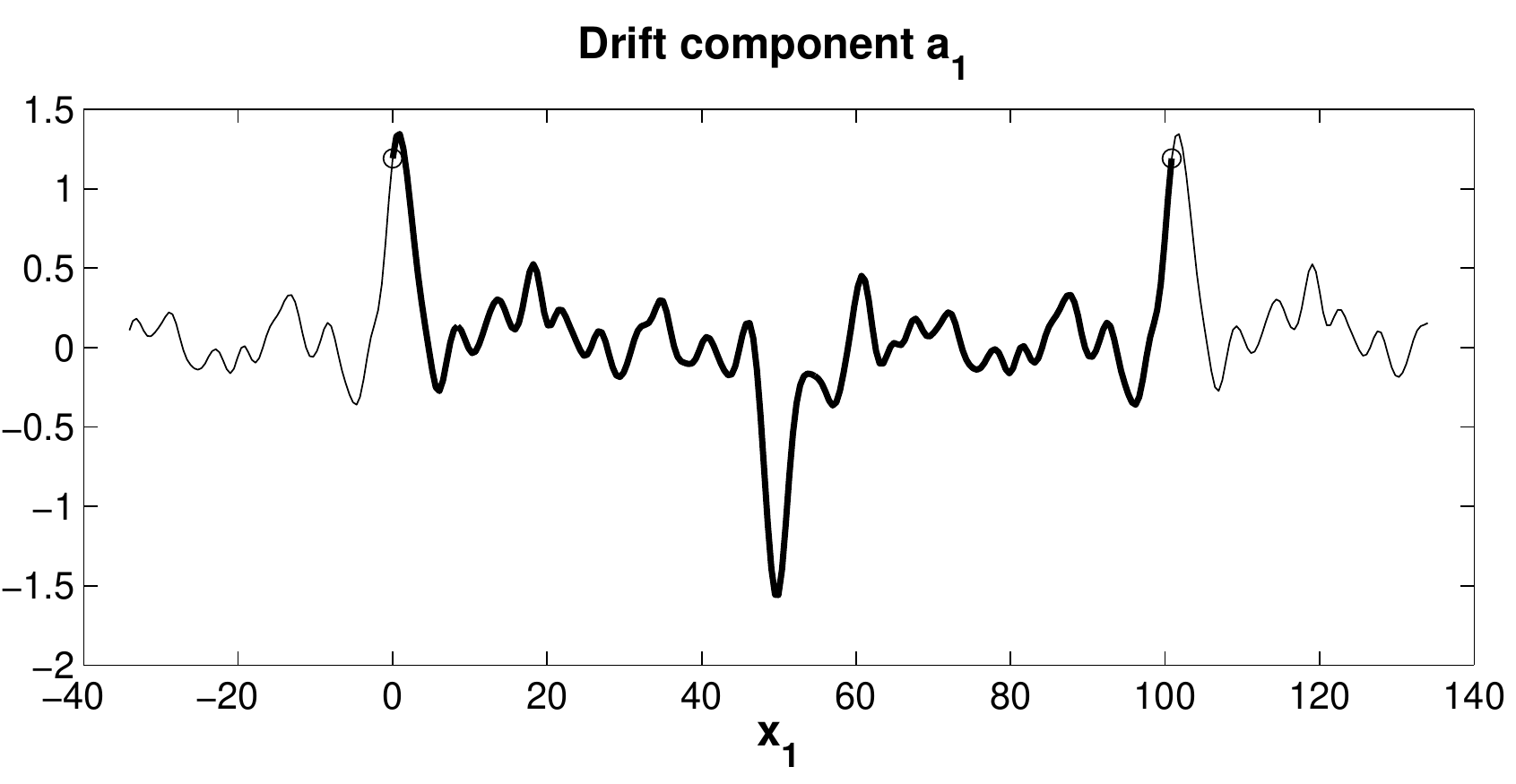}}
  \caption{The term $\average{\driftone}{\sampleset}$ is supposed
    to approach zero as the number of configurations, and \tend,
    increases, provided that the interfaces are stationary. 
    Though the fluctuations are large here, they are much smaller
    than in Figure~\ref{fig:conv_drift_zero}.
    When the fluctuations decrease a pattern appears with peaks
    at the two interfaces. This supports the observation, from
    the computed $\average{\pfen}{\sampleset}$ in
    Figure~\ref{fig:phasefield_not_eq}, that the two phase system
    is not in equilibrium yet and the interfaces are not really
    stationary on the time scale of the average.
  }
  \label{fig:conv_drift_one}
\end{figure}

\begin{figure}[hbp]
  \centering
  \subfigure[Mean based on 111 configurations, $\tend=0.8875$]
  {\label{subfig:d2m_111}
    \includegraphics[width=6.5cm]
    {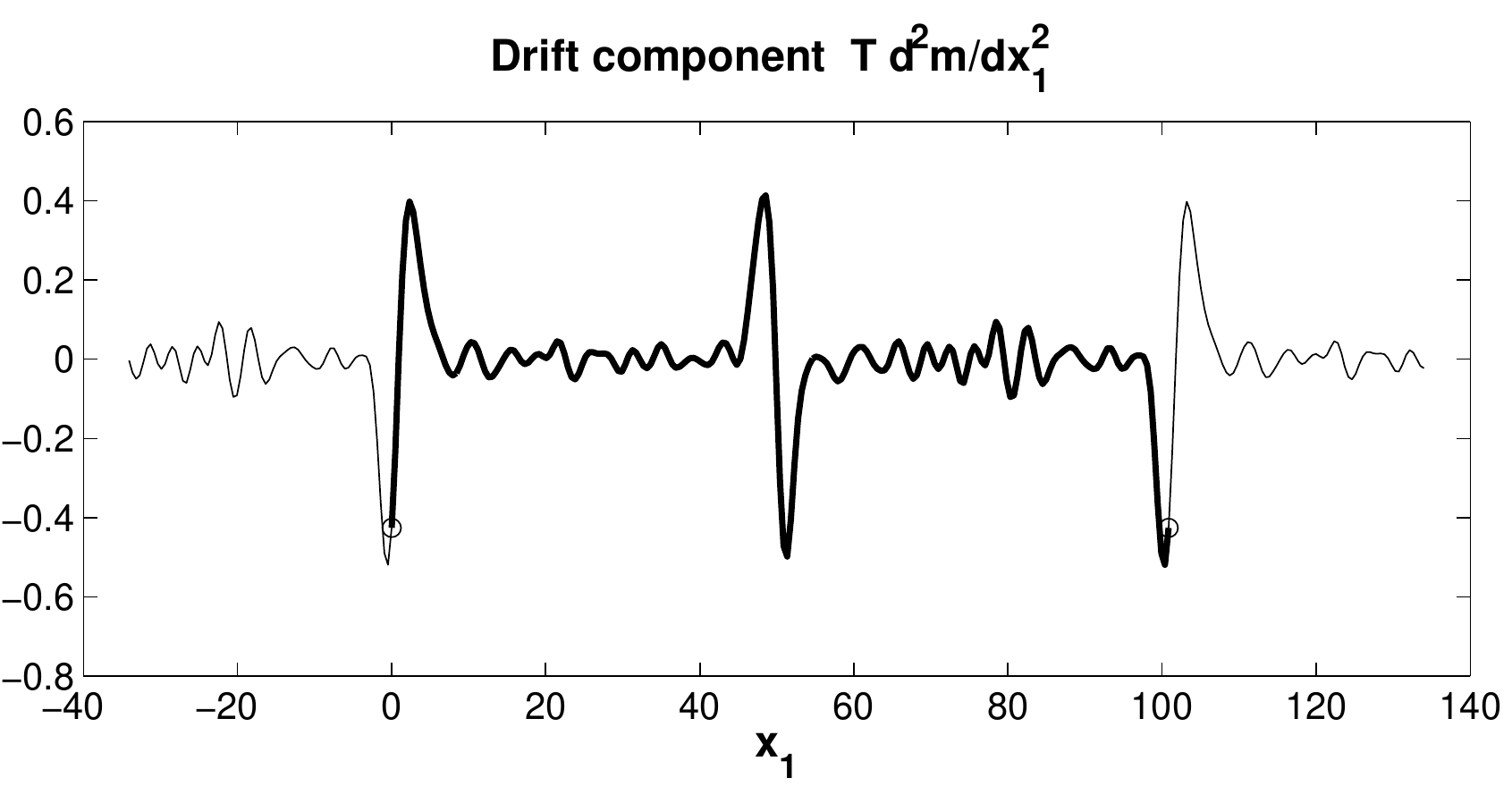}}
  \subfigure[Mean based on 444 configurations, $\tend=0.8875$]
  {\label{subfig:d2m_444}
    \includegraphics[width=6.5cm]
    {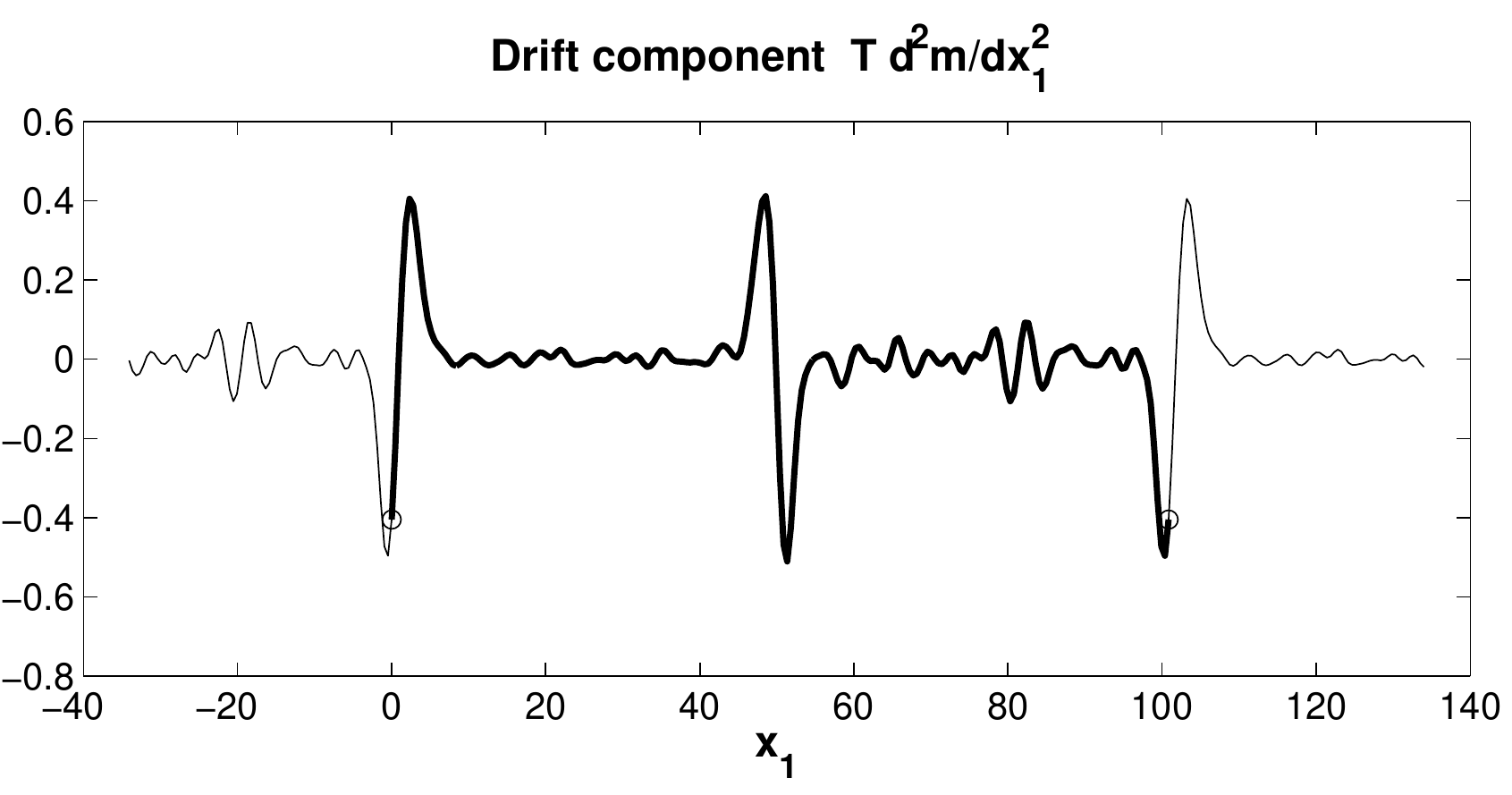}}
  \subfigure[Mean based on 444 configurations, $\tend=0.2220$]
  {\label{subfig:d2m_444_dense}
    \includegraphics[width=6.5cm]
    {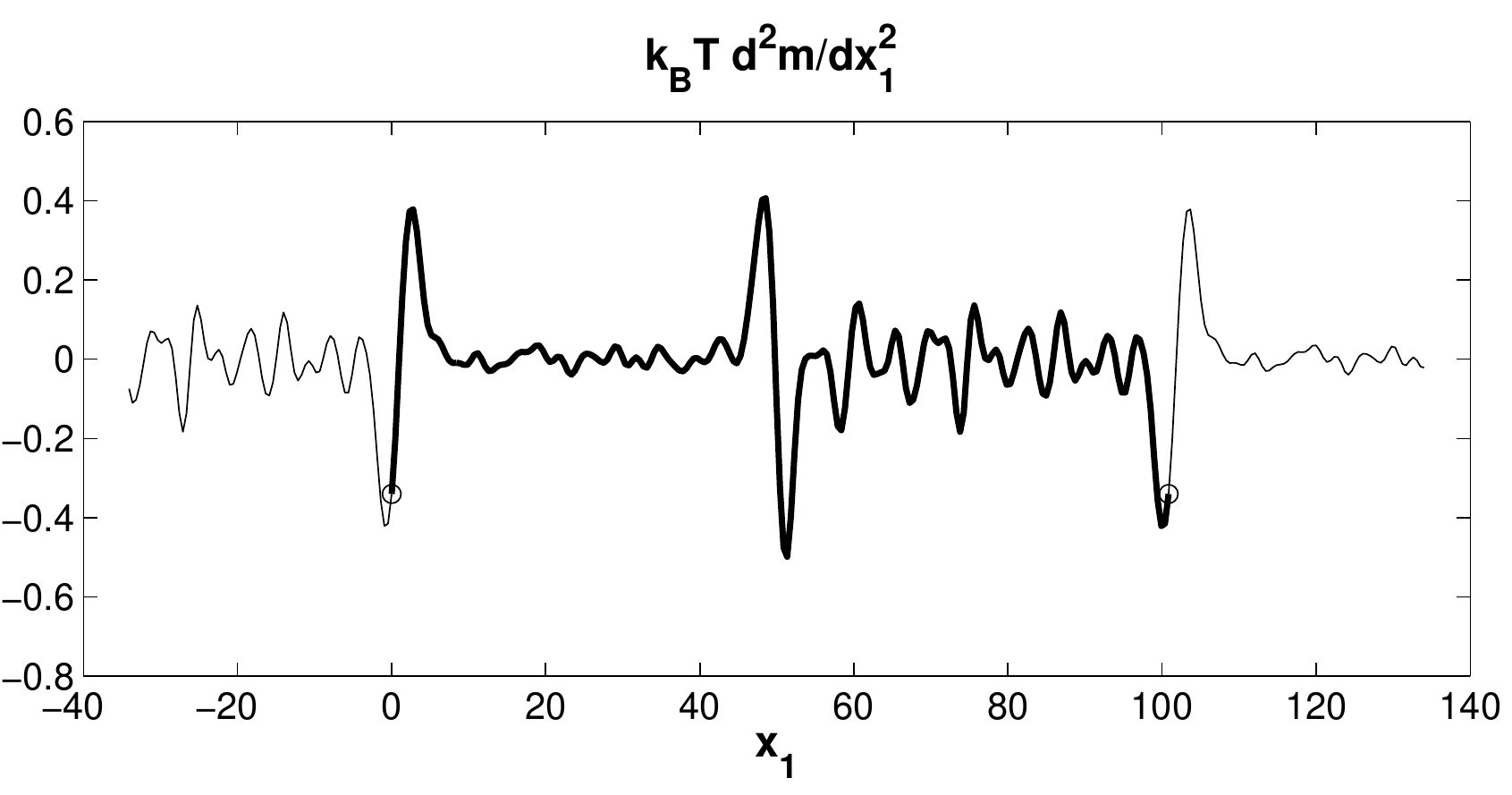}}
  \subfigure[Mean based on 1775 configurations, $\tend=0.8875$]
  {\label{subfig:d2m_1775}
    \includegraphics[width=6.5cm]
    {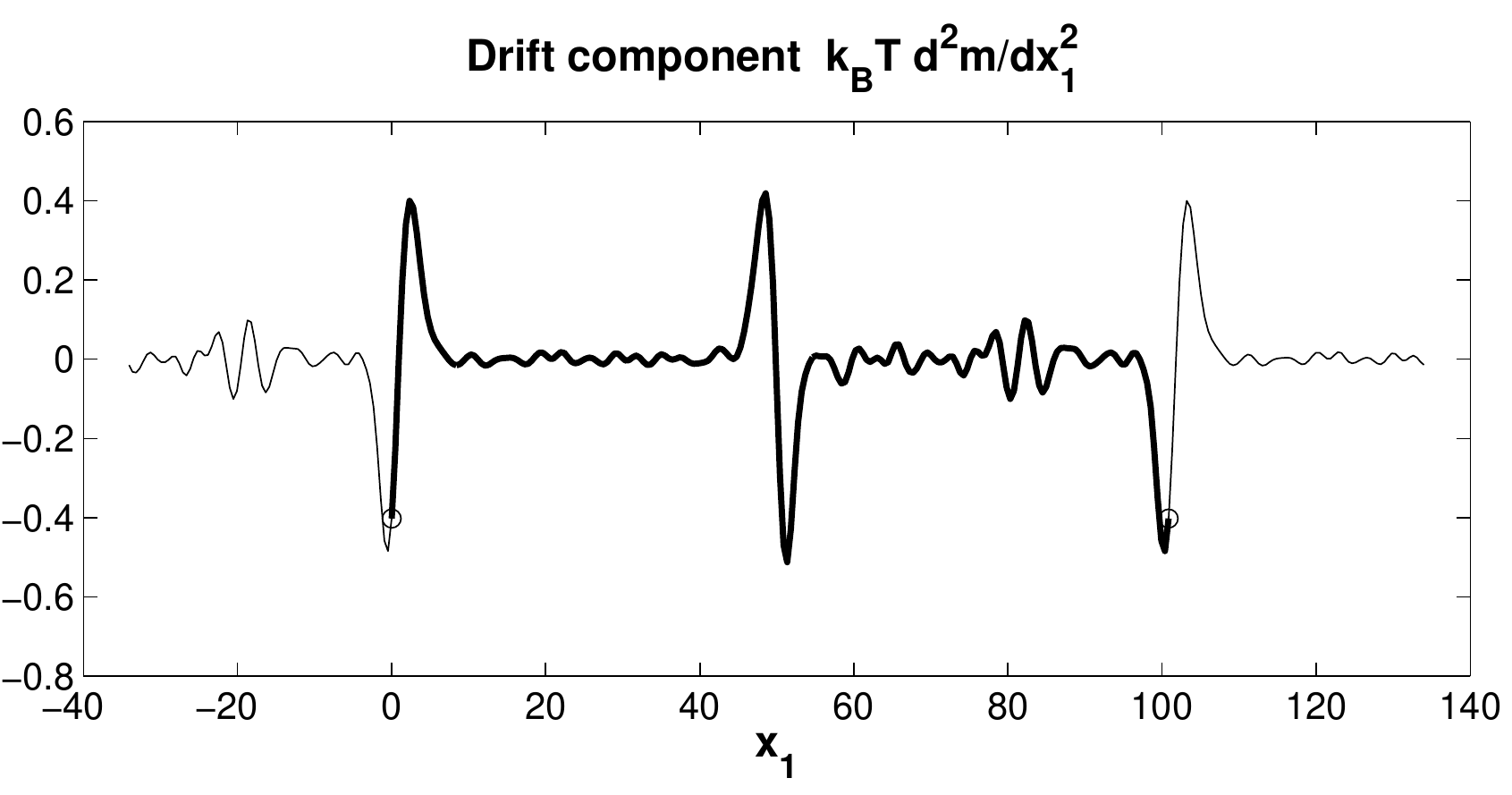}}
  \caption{The average
    $\kb\temp\ddxett\average{\driftone}{\sampleset}$ converges
    faster than the two other terms in
    $\average{\driftmd}{\sampleset}$. 
    The fluctuations are larger in subfigure~(c) than in~(b),
    which indicates that the error is dominated by the length of
    the averaging time interval rather than the number of
    configurations sampled within the time interval.
  }
  \label{fig:conv_drift_d2m}
\end{figure}

The average \average{\driftone(x_1;\cdot)}{\sampleset}, contains
first order differentials of the potential, but only to the
first power. The convergence of is faster than that of
\average{\driftzero(x_1;\cdot)}{\sampleset}, but the computed
averages in Figure~\ref{fig:conv_drift_one} still show
significant fluctuations.
The final term in $\average{\driftmd(x_1;\cdot)}{\sampleset}$ is 
$\kb\temp \ddxett \average{\pfen(x_1;\cdot)}{\sampleset}$, which
only depends on the potential and not its derivatives. This
average converges faster than the other two and, even after two
differentiations with respect to $x_1$, the fluctuations are
small compared to the distinct structures at the interfaces; see
Figure~\ref{fig:conv_drift_d2m}.

\subsubsection{Obtaining the phase-field double-well potential
  from the drift}
\label{sec:double-well}

When defining a phase-field variable in terms the potential
energy in the microscale model in Section~\ref{sec:Intro}, the
goal was to compute a reaction--diffusion equation, like the
Allen-Cahn equation~\eqref{eq:AllenCahn}, for the coarse-grained
phase-field. 
In a one dimensional problem, with $\temp\equiv\Tmelt$ and $k_1$
constant, the Allen-Cahn equation reduces to 
\begin{align}
  \label{eq:Reaction-Diffusion}
  \frac{\partial \pfgen}{\partial t}
  & = 
  k_1\ddxett\pfgen - k_2\dwell'(\pfgen) + noise,
\end{align}
where the derivative of the double-well
potential \dwell\ gives the reaction part in this
reaction--diffusion equation.
Now, the coarse-grained equation
\begin{align*}
  d\pfcgt(x_1) & = 
  \left(
  \kb\temp\frac{\partial^2}{\partial x_1^2}\pfcgt(x_1)
  + \frac{\partial}{\partial x_1} \driftcgx_1(x_1)
  + \driftcgx_0(x_1)\right)dt + 
  \sumall[\nrw]{j} \diffucgx_j(x)\;d\indepw_j^t,  
\end{align*}
where 
\begin{align*}
  \driftcgx_1(x_1) & = 
  \average{\driftone}{\sampleset}(x_1),
  &
  \driftcgx_0(x_1) & = 
  \average{\driftzero}{\sampleset}(x_1)
  &&, \text{for $x_1\in D_K$,} 
\end{align*}
and the diffusion coefficient vectors, $\diffucgx_j$, are
obtained from the factorisation~\eqref{eq:def_diffucgxmat}, is a
stochastic convection--reaction--diffusion equation.
As the described above the time averaged drift is zero in a
stationary situation, but in the computations presented here the
fluctuations are still too large. In the ideal situation for a
stationary interface, when all three components in the drift
average have converged, the convection should vanish, that is
\begin{align*}
  \frac{\partial}{\partial x_1} \driftcgx_1
  & \equiv 0,
\end{align*}
and the reaction and diffusion parts should cancel each other, so
that
\begin{align}
  \label{eq:stat_react_diffu}
  0 & = 
  \kb\temp\frac{\partial^2}{\partial x_1^2}\pfcgt(x_1)
  + \driftcgx_0(x_1).
\end{align}
The second best thing, when some of the computed averages contain
too large errors, is to extract information from the most
accurate part, that is 
$\kb\temp\frac{\partial^2}{\partial x_1^2}
\pfen_\mathrm{av}(x_1)$.
Assuming that this computed average already is close to what it
would be in the ideal situation, an approximation of
the reaction term can be obtained
from~\eqref{eq:stat_react_diffu}. 

\begin{figure}[hbp]
  \centering
  \includegraphics[width=7.5cm]
  {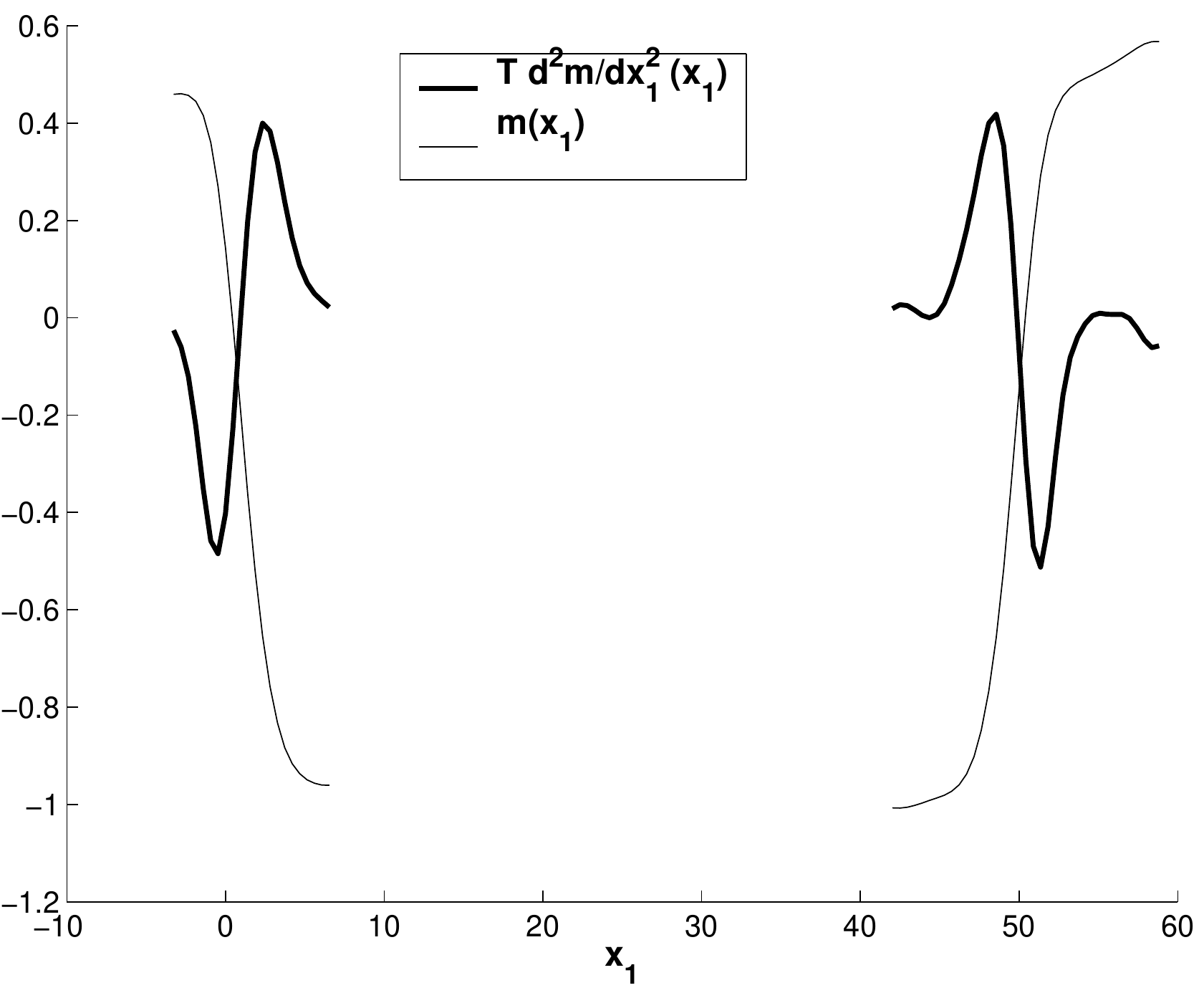}  
  \caption{The computed $\pfen_\mathrm{av}(x_1)$ in its monotone
    intervals in the interfaces together with the corresponding
    diffusion part of the drift
    $\kb\temp\ddxett\average{\driftone}{\sampleset}$. 
    The curves shown are part of the those in 
    Figure~\ref{subfig:rho_pf_O2}
    and Figure~\ref{subfig:d2m_1775}.}
  \label{fig:monoton_interface}
\end{figure}

\begin{figure}[hbp]
  \centering
  \subfigure[Double well potential, $\molliscale=1.0$]{
    \label{subfig:double_well_eps_1}
    \includegraphics[width=6cm]
    {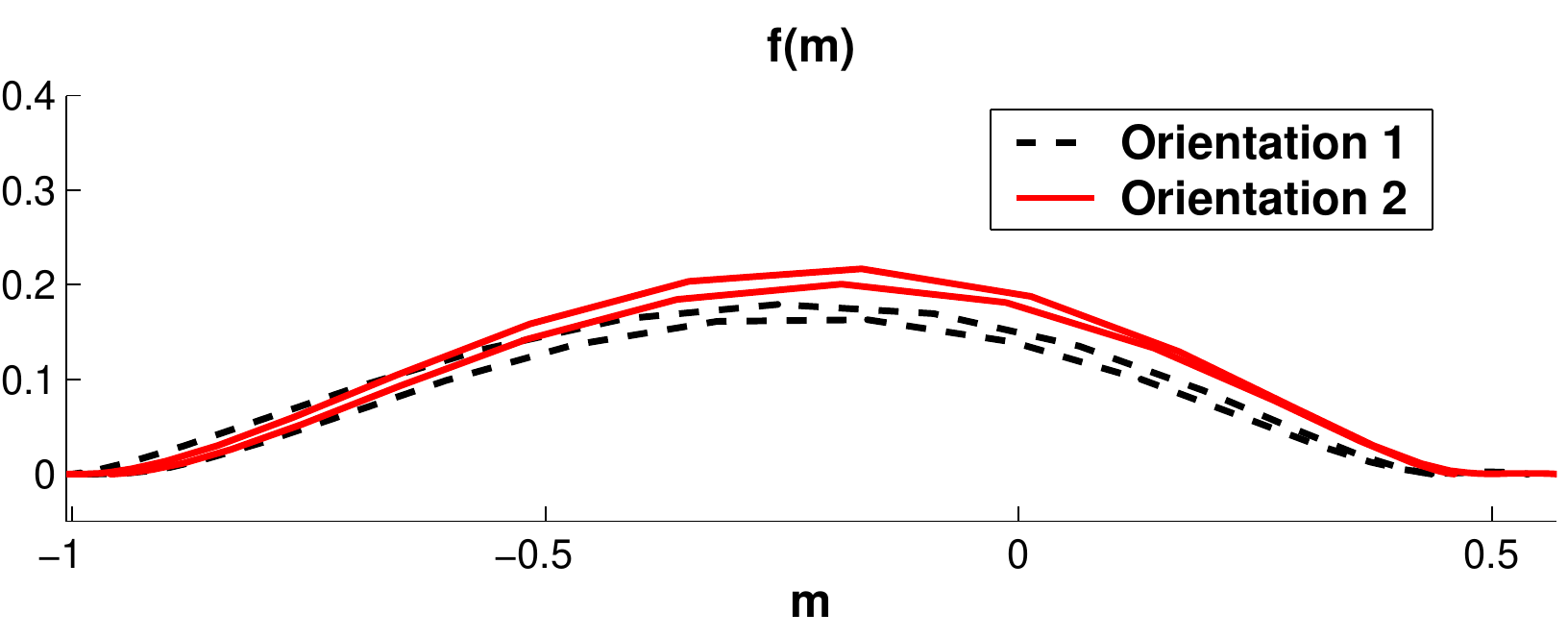}
  }
  \subfigure[Orientation 2, different \molliscale]{
    \label{subfig:double_well_comp_eps}
    \includegraphics[width=6cm]
    {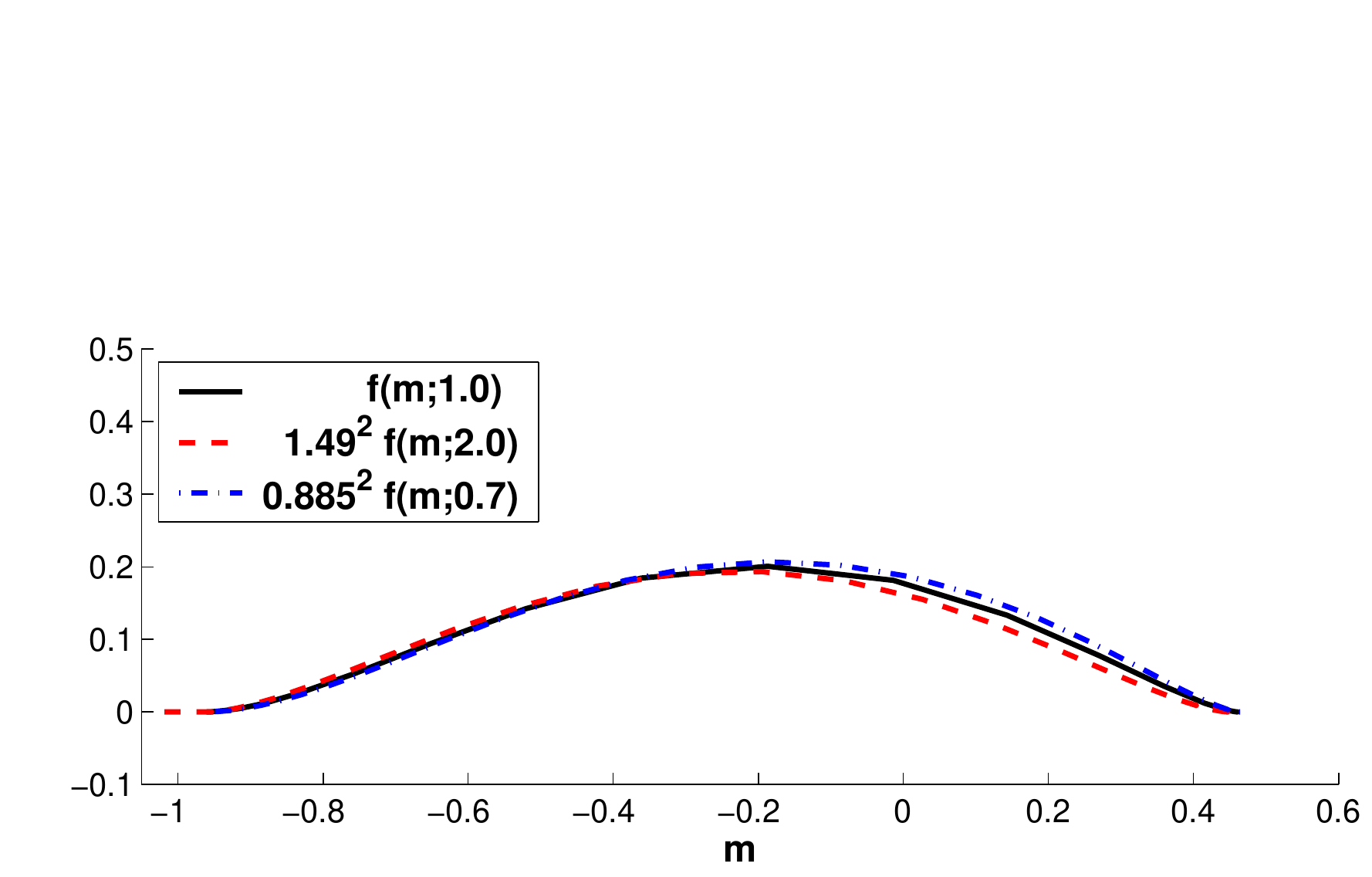}
  }
  \caption{\newline
    (a) The computed double well potentials from both
    simulation O1 and O2 using $\pfen_\mathrm{av}$ shown in 
    Figure~\ref{fig:phasefield_not_eq} and the corresponding 
    $\kb\temp\frac{\partial^2}{\partial x_1^2}
    \pfen_\mathrm{av}(x_1)$.
    \newline
    (b) The computed double well potentials from one of the
    interfaces in O2, using three different values of the
    smoothing parameter \molliscale\ in the mollifier. Since the
    interface width varies with \molliscale\ the height of the
    potential barriers vary with \molliscale. Here double-wells
    have been rescaled with factors obtained in the analysis of
    the \molliscale-dependence in
    Figure~\ref{fig:interface_scaling} to compare the shape of
    the curves. 
  }
  \label{fig:double_well_computed}
\end{figure}

The expression of the drift in the coarse-grained
equation~\eqref{eq:sde_pfcg} as a function of the coarse-grained
phase-field \pfcg\ in the interface regions, instead of the space
variable $x_1$, assumes monotonicity of the phase-field near the
interfaces to allow the inversion in~\eqref{eq:def_coeff_cg}.
Figure~\ref{fig:monoton_interface} shows $\pfen_\mathrm{av}(x_1)$
and 
$\kb\temp\frac{\partial^2}{\partial x_1^2}\pfen_\mathrm{av}(x_1)$ 
in the interval of monotonicity for $\pfen_\mathrm{av}(x_1)$ in
the simulation~O2.
Using the computed 
$\kb\temp\frac{\partial^2}{\partial x_1^2}\pfen_\mathrm{av}(x_1)$ 
in~\eqref{eq:stat_react_diffu}, gives
\begin{align*}
  \driftcgx_0(x_1) & = 
  -\kb\temp\frac{\partial^2}{\partial x_1^2}
  \pfen_\mathrm{av}(x_1).
\end{align*}
Inverting the computed function $\pfen_\mathrm{av}(x_1)$ in the
interface intervals, the derivative of the double-well potential
\dwell\ can be identified as 
\begin{align*}
  \dwell'(\pfcg) & = 
  \driftcgx_0(\pfen_\mathrm{av}^{-1}(\pfcg)).
\end{align*}
Integration with respect to $\pfcg$ in the interval between
$\pfcg_\mathrm{solid}$ and $\pfcg_\mathrm{liquid}$ gives the
double-well potentials shown in
Figure~\ref{subfig:double_well_eps_1}. 
As expected the potentials obtained from the two different
simulations O1 and O2 are slightly different. However, the
potentials obtained from the two different interfaces in one
molecular dynamics simulation cell also differ slightly and it is
not possible to say that difference between simulations O1 and O2
depend on the orientation of the interfaces with respect to the 
crystal lattice.
The computed double wells seem to be qualitatively right. 

\subsection{The averaged diffusion matrix \diffumat\ and the
  coarse-grained diffusion coefficients $\diffucgx_j$.}
\label{sec:computed_diffu}

\begin{figure}[hbp]
  \centering
  \subfigure[\diffumat]{
    \label{subfig:diffumat}
    \includegraphics[width=6.5cm]
    {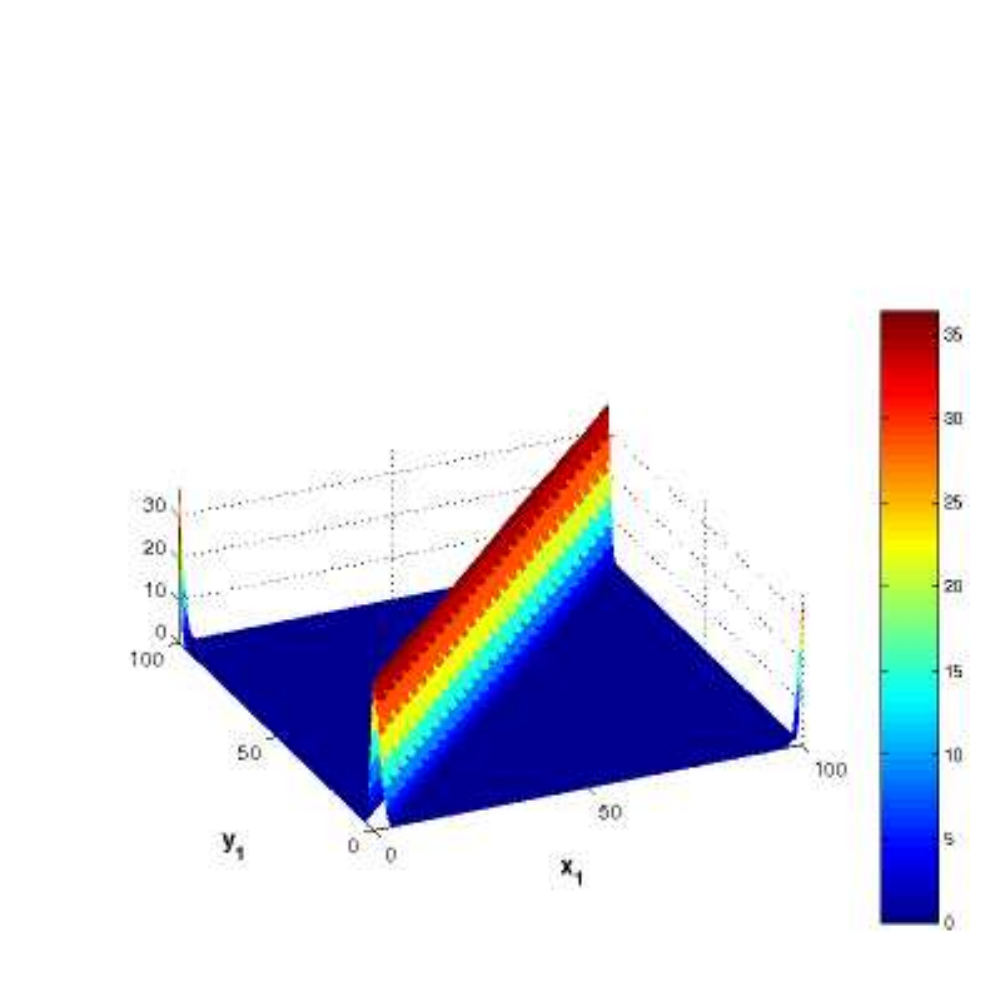}
  }
  \subfigure[$\diffucgxmat$]{
    \label{subfig:diffucgxmat}
    \includegraphics[width=6.5cm]
    {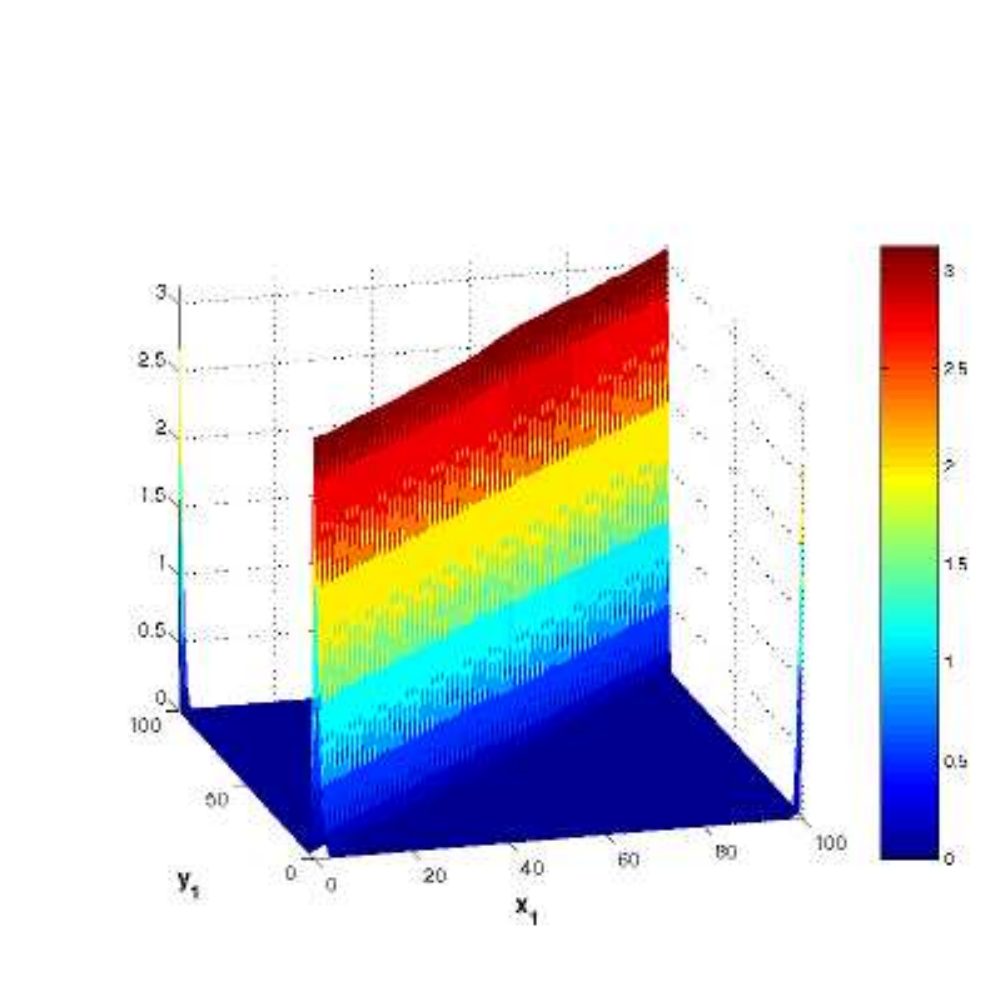}
  }
  \subfigure[Some $\diffucgx_j$:s]{
    \label{subfig:nine_bjs}
    \includegraphics[width=6.5cm]{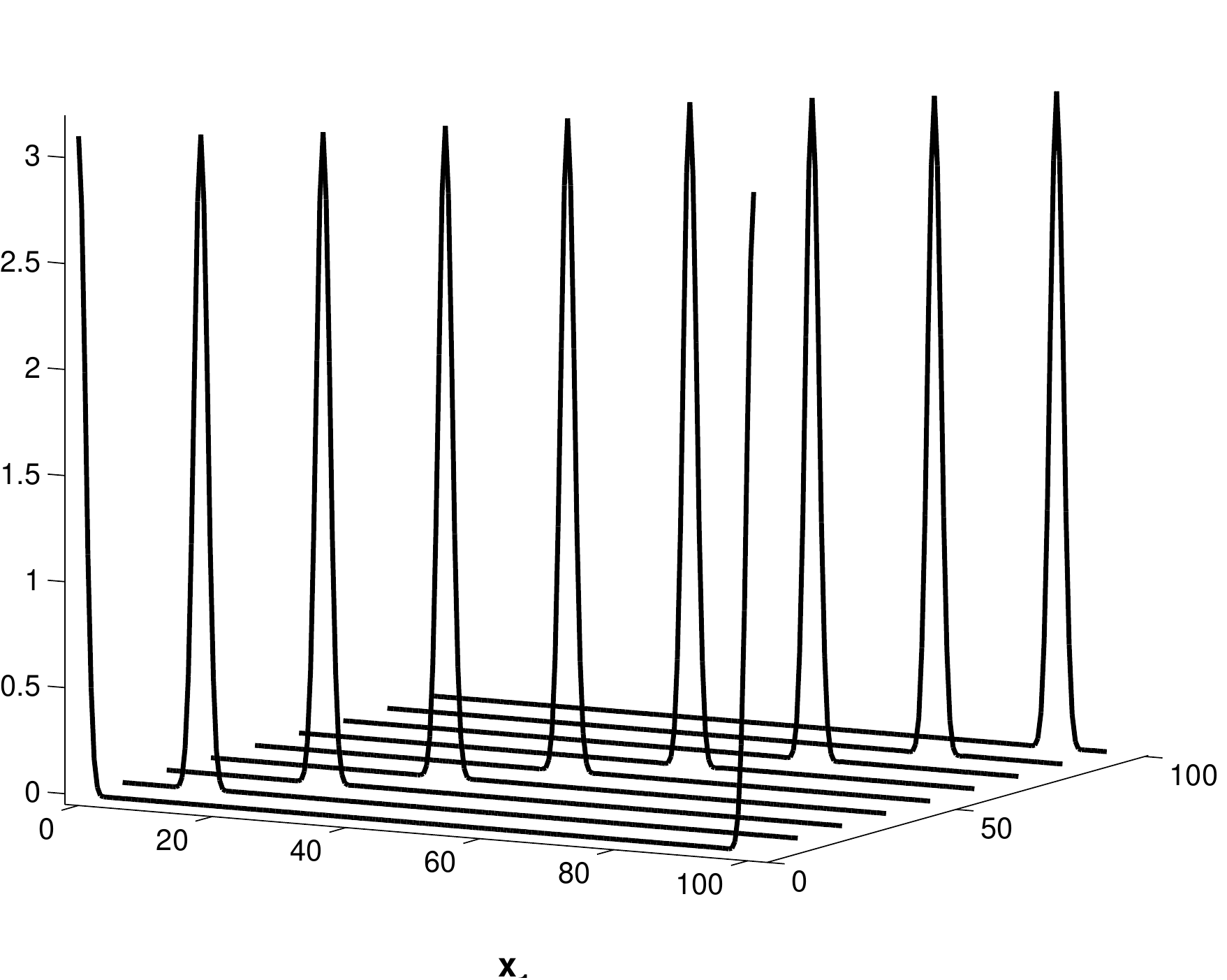}
  }
  \caption{The computed average diffusion matrix \diffumat, for
    $\molliscale=1.0$, using the same configurations from
    simulation O2 as in Figure~\ref{subfig:rho_pf_O2} and
    Figure~\ref{subfig:d2m_1775}, is shown in~(a).
    The square root \diffucgxmat\ of \diffumat, as defined
    in~\eqref{eq:def_diffucgxmat} is shown in~(b). 
    The individual columns in \diffucgxmat\ are the diffusion
    coefficient functions, $\diffucgx_j$, in the stochastic
    differential equation for the coarse-grained phase-field
    \pfent. Some of these column vectors have been plotted as
    functions of the space variable $x_1$ in~(c).
    The support of each $\diffucgx_j$ is centred around the grid
    point $x_1^j$.}
  \label{fig:diffusions_eps_1}
\end{figure}

\begin{figure}[hbp]
  \centering
  \begin{picture}(50,180)(0,0)
    \put(0,0){\makebox(0,0)}
  \end{picture}
  \subfigure[FCC, \molliscale=2.0]{
    \label{subfig:mean_bj_FCC_2_0}
    \includegraphics[width=2.5cm]{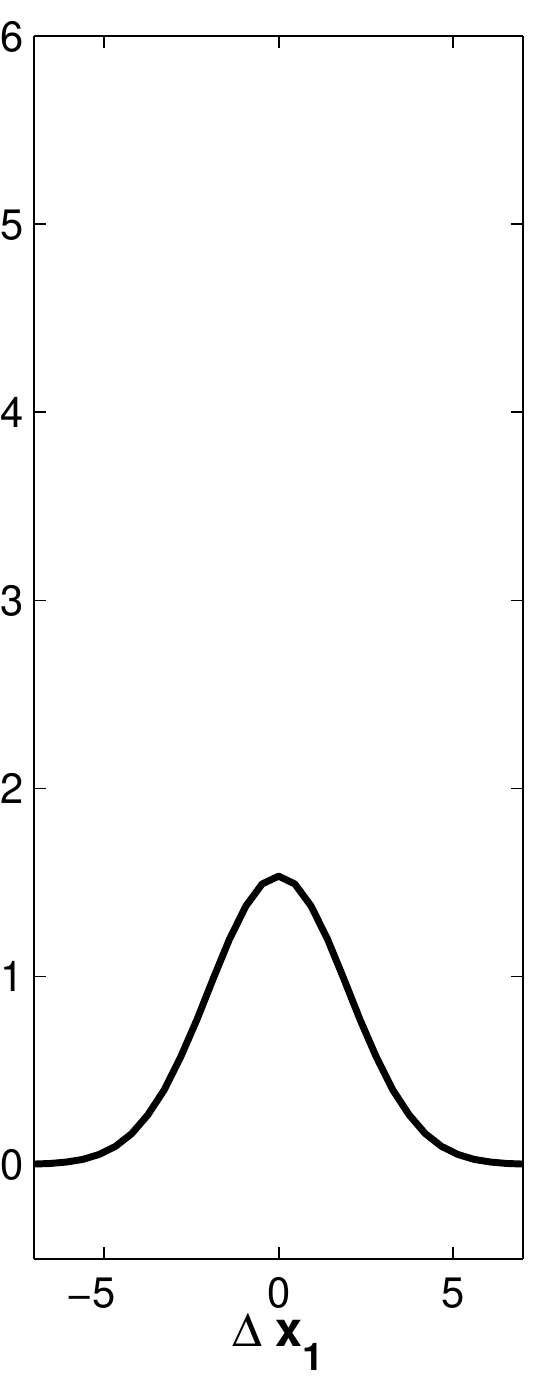}
  }
  \subfigure[FCC, \molliscale=1.0]{
    \label{subfig:mean_bj_FCC_1_0}
    \includegraphics[width=2.5cm]{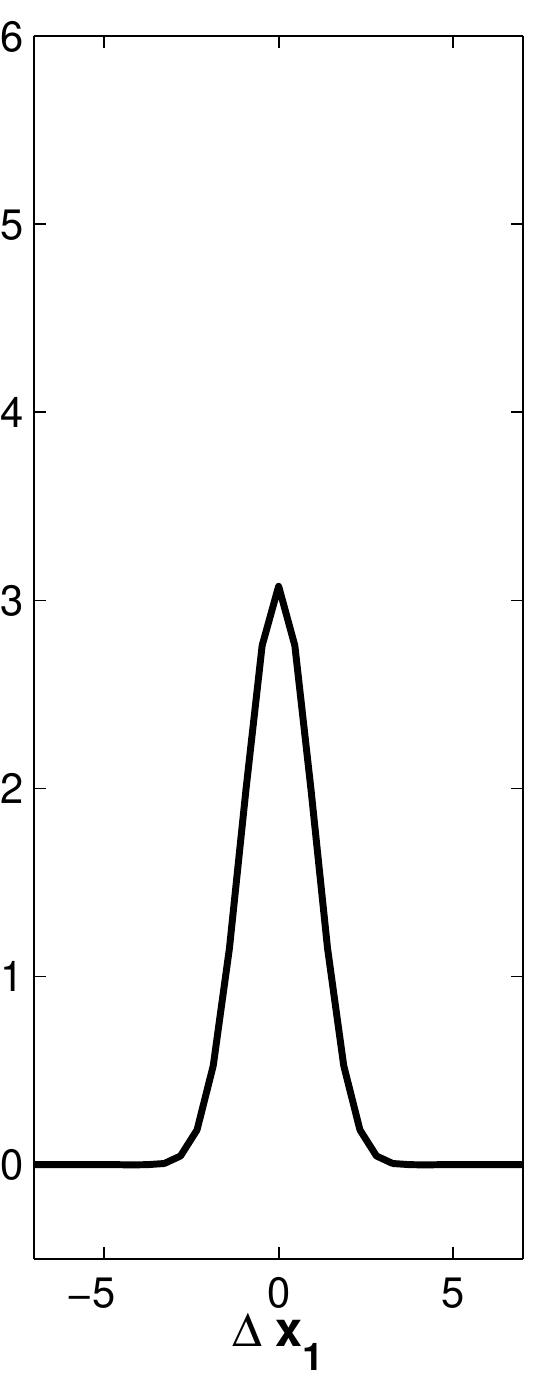}
  }
  \subfigure[FCC, \molliscale=0.55]{
    \label{subfig:mean_bj_FCC_0_55}
    \includegraphics[width=2.5cm]{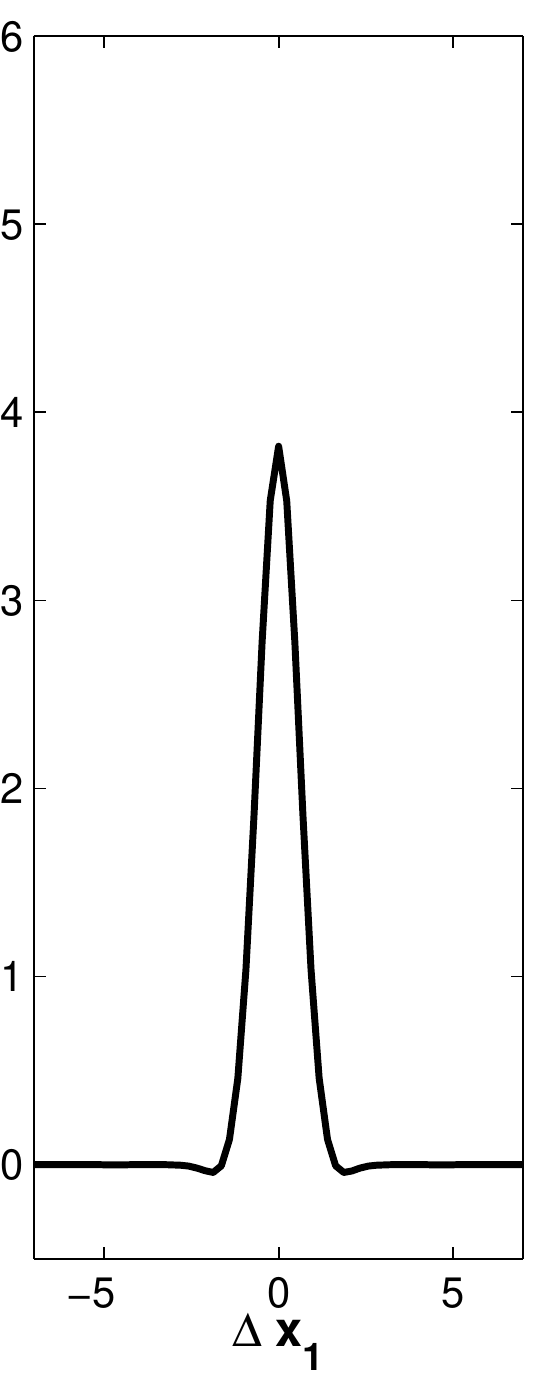}
  }
  \begin{picture}(50,180)(0,0)
    \put(0,0){\makebox(0,0)}
  \end{picture}
  \begin{picture}(50,180)(0,0)
    \put(0,0){\makebox(0,0)}
  \end{picture}
  \subfigure[Liquid, \molliscale=2.0]{
    \label{subfig:mean_bj_Liquid_2_0}
    \includegraphics[width=2.5cm]{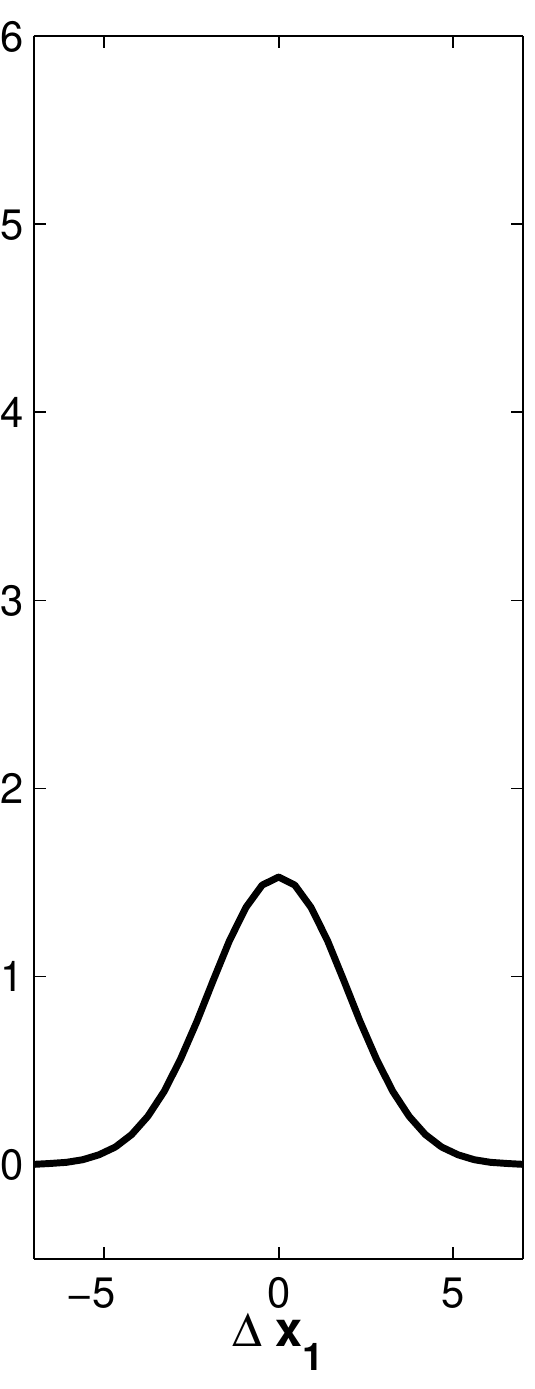}
  }
  \subfigure[Liquid, \molliscale=1.0]{
    \label{subfig:mean_bj_Liquid_1_0}
    \includegraphics[width=2.5cm]{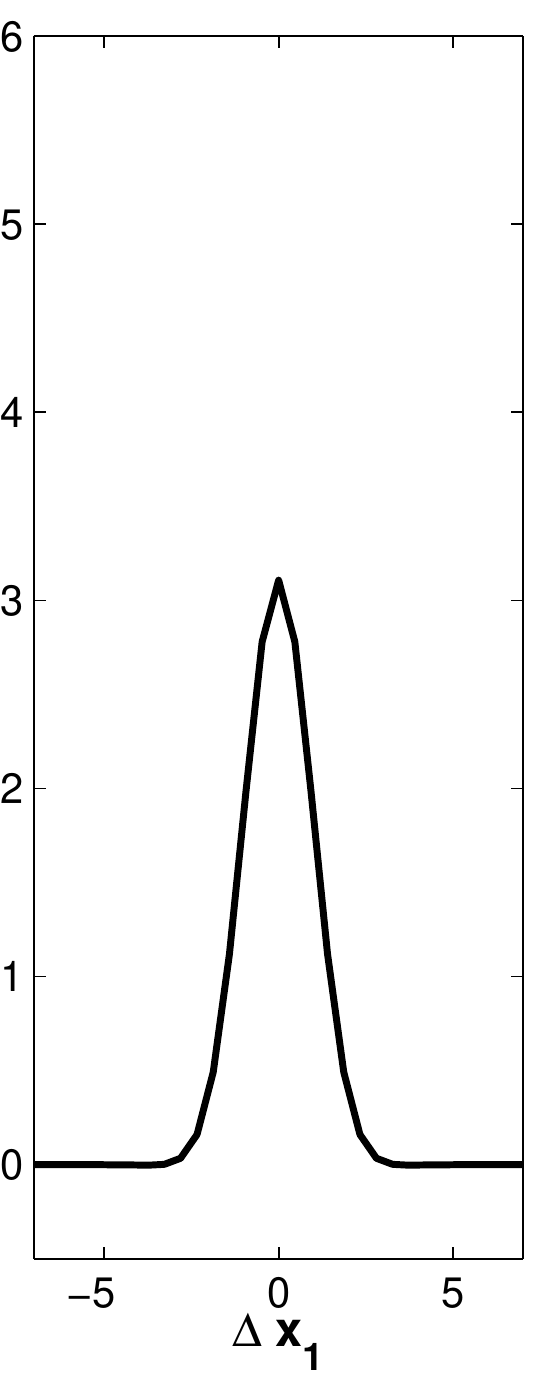}
  }
  \subfigure[Liquid, \molliscale=0.55]{
    \label{subfig:mean_bj_Liquid_0_55}
    \includegraphics[width=2.5cm]{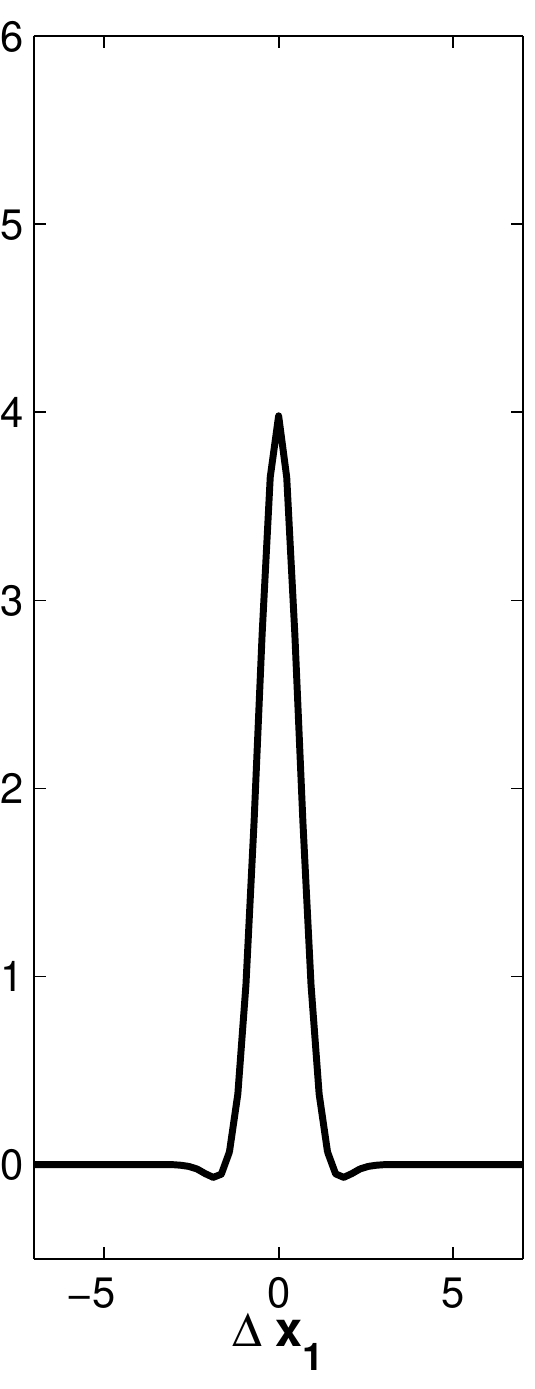}
  }
  \begin{picture}(50,180)(0,0)
    \put(0,0){\makebox(0,0)}
  \end{picture}
  \caption{The average diffusion coefficient functions,
    $\tilde{\diffucgx}(\Delta x_1) =
    \mathrm{mean}\left\{\diffucgx_j(x_1^j+\Delta x_1)\right\}$
    have been computed for different values of \molliscale, with
    the mean taken over points $x_1^j$ in the interior of the
    solid and the liquid domains, respectively. 
    The configurations used are the same as in
    Figure~\ref{fig:diffusions_eps_1}. 
    The difference between the solid and liquid parts is small.}
  \label{fig:mean_bjs}
\end{figure}

The final component to extract in the coarse-grained model is the
diffusion in the stochastic differential equation for \pfcgt.
Using $\molliscale=1.0$ and the same 1775 configurations that
were used in the computation of the averaged phase-field and
drift for simulation O2, the averaged diffusion matrix \diffumat,
has been computed, with the result shown in
Figure~\ref{subfig:diffumat}.
As described in Section~\ref{sec:coarse}, 
the square root of \diffumat is computed by an eigenvector
decomposition where all negative eigenvalues are set to zero;
the result is shown in Figure~\ref{subfig:diffucgxmat}.
The negative eigenvalues are very small in absolute value,
compared to the dominating positive ones, so the error made by
neglecting them is insignificant when
$\diffucgxmat\transpose{\diffucgxmat}$ is compared to \diffumat.
By choosing the diffusion coefficients $\diffucgx_j$ in the
coarse-grained stochastic differential equation as the columns
of \diffucgxmat, they become localised in space; see
Figure~\ref{subfig:nine_bjs}.
With $\molliscale=1.0$ the observed difference between the
diffusion in the solid part and the liquid part is small, as
shown in Figure~\ref{fig:mean_bjs}.

\subsection{Dependence on the smoothing parameter}
\label{sec:eps_dep}

The mollifier \molli\ includes a parameter, \molliscale, determining
the scale on which the local average is taken. This is in itself an
ad~hoc variable in the micro model and it is important to analyse
its effects on the computed quantities. 

\begin{figure}[hbp]
  \centering
  \subfigure[Orientation 1]{
    \label{fig:atom_plane_O1}
    \includegraphics[height=5cm]{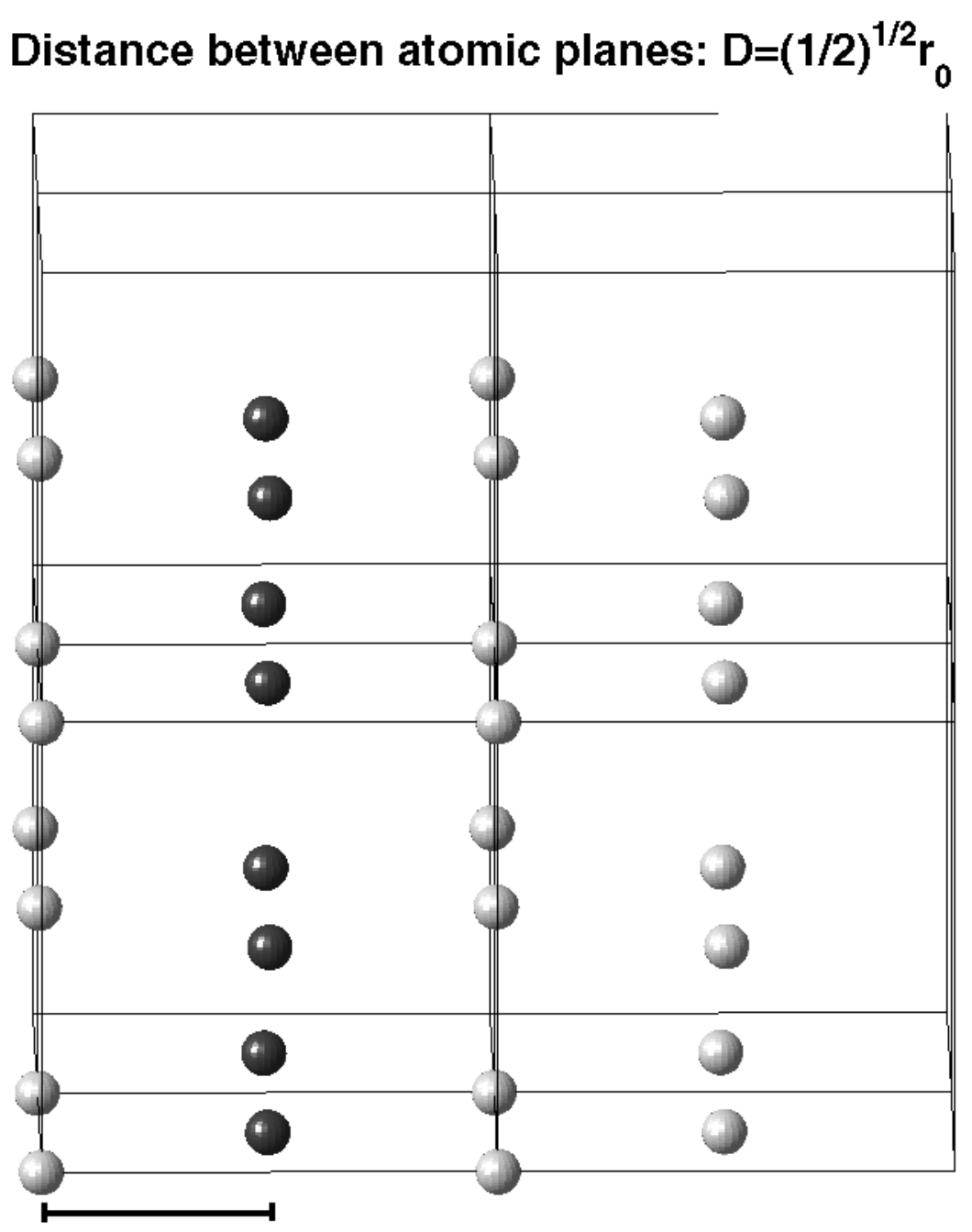}
  }
  \hspace{2cm}
  \subfigure[Orientation 2]{
    \label{fig:atom_plane_O2}
    \includegraphics[height=5cm]{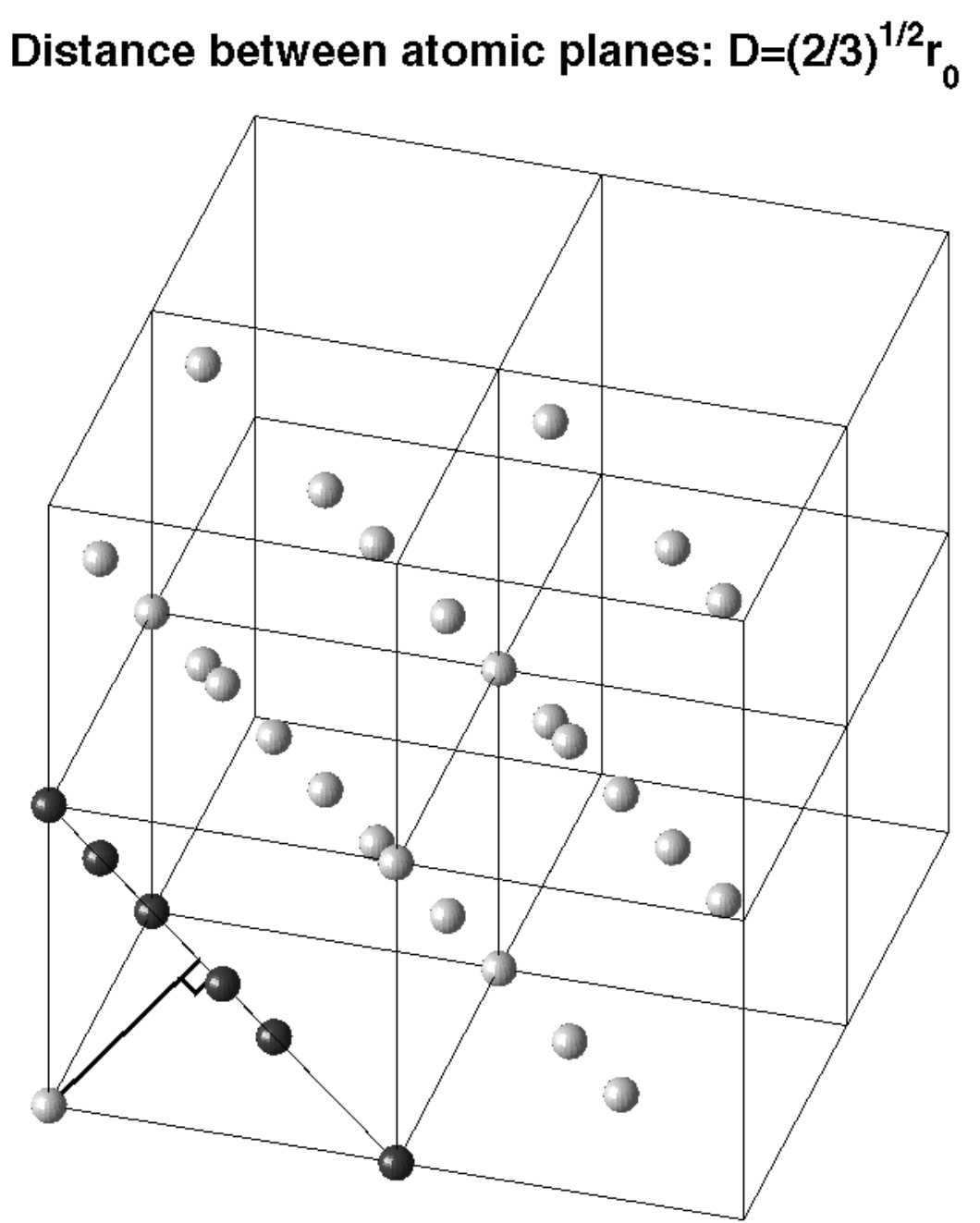}
    }
  \caption{The distance between two adjacent atom layers in a
    perfect FCC 
    lattice is~$\sqrt{1/2}\,r_0$ in orientation 1
    and $\sqrt{2/3}\,r_0$ in orientation 2, where
    $r_0$ is the nearest neighbour distance.}
  \label{fig:atom_planes}
\end{figure}

\begin{figure}[hbp]
  \centering
  \subfigure[Orientation 1]{
    \label{fig:pf_rho_O1_epsi_0_45}
    \includegraphics[width=6.5cm]
    {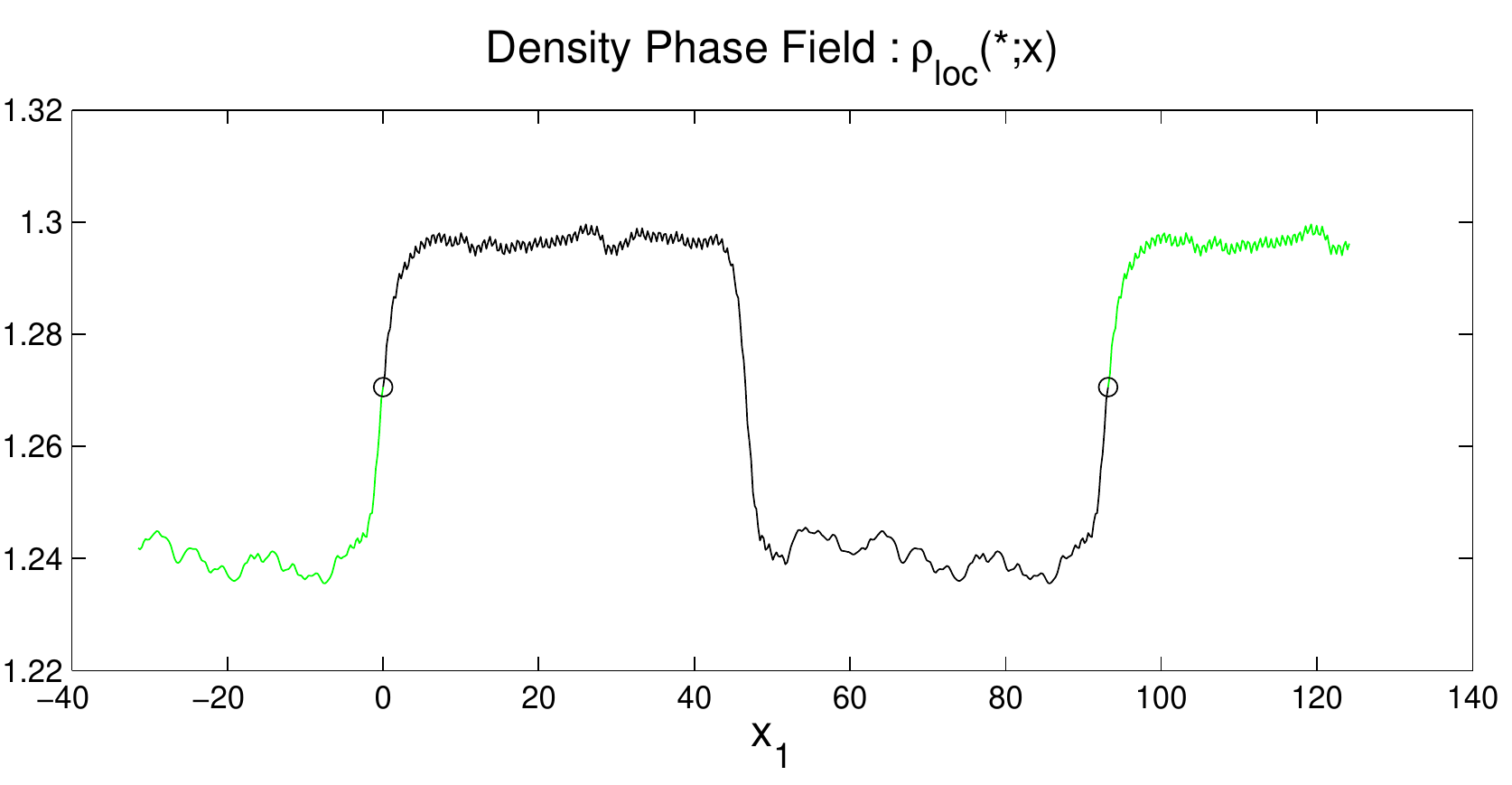}
  }
  \subfigure[Orientation 2]{
    \label{fig:pf_rho_O2_epsi_0_45}
    \includegraphics[width=6.5cm]
    {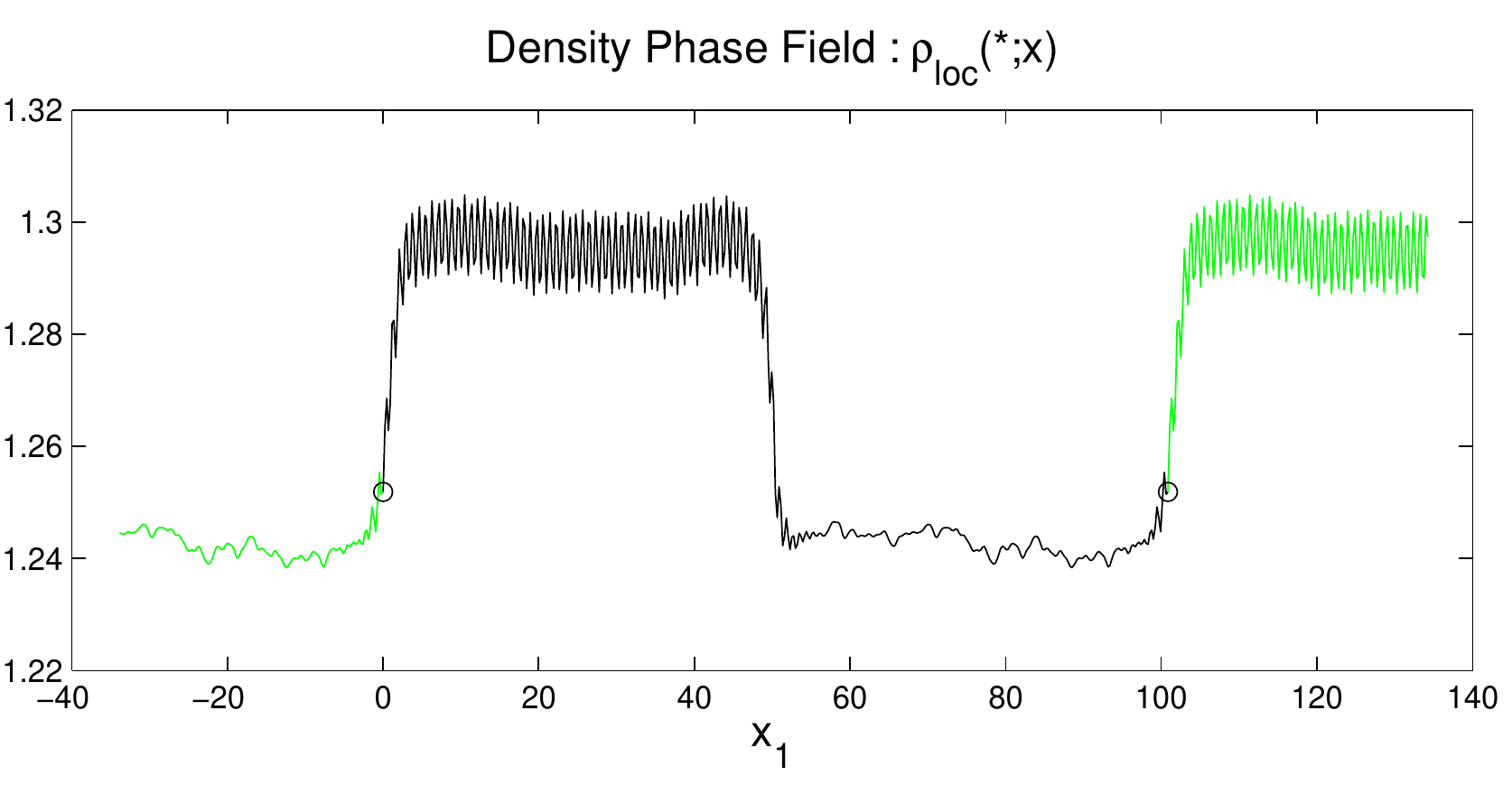}
  }
  \subfigure[Orientation 1]{
    \label{fig:pf_pot_O1_epsi_0_45}
    \includegraphics[width=6.5cm]
    {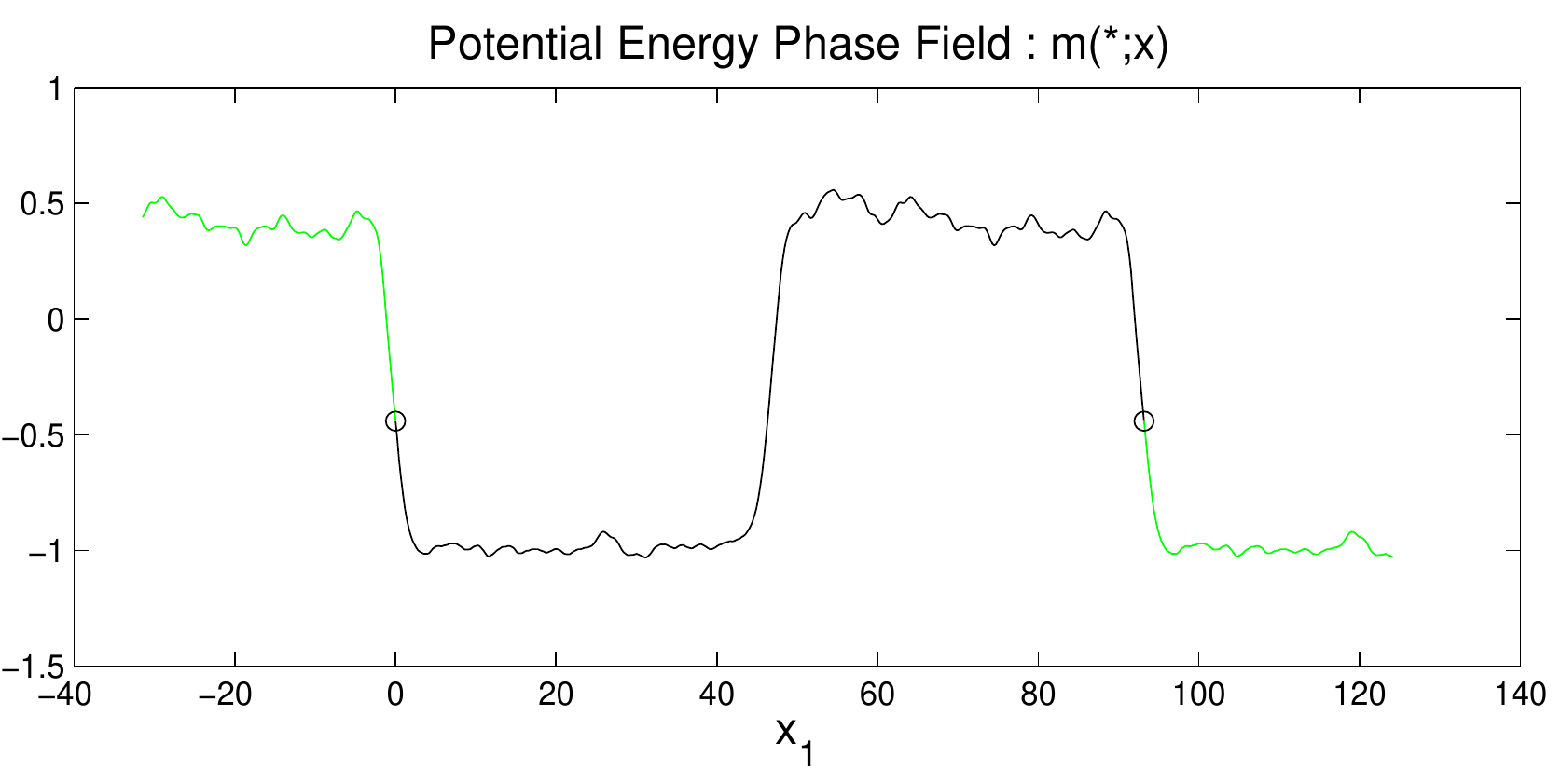}
  }
  \subfigure[Orientation 2]{
    \label{fig:pf_pot_O2_epsi_0_45}
    \includegraphics[width=6.5cm]
    {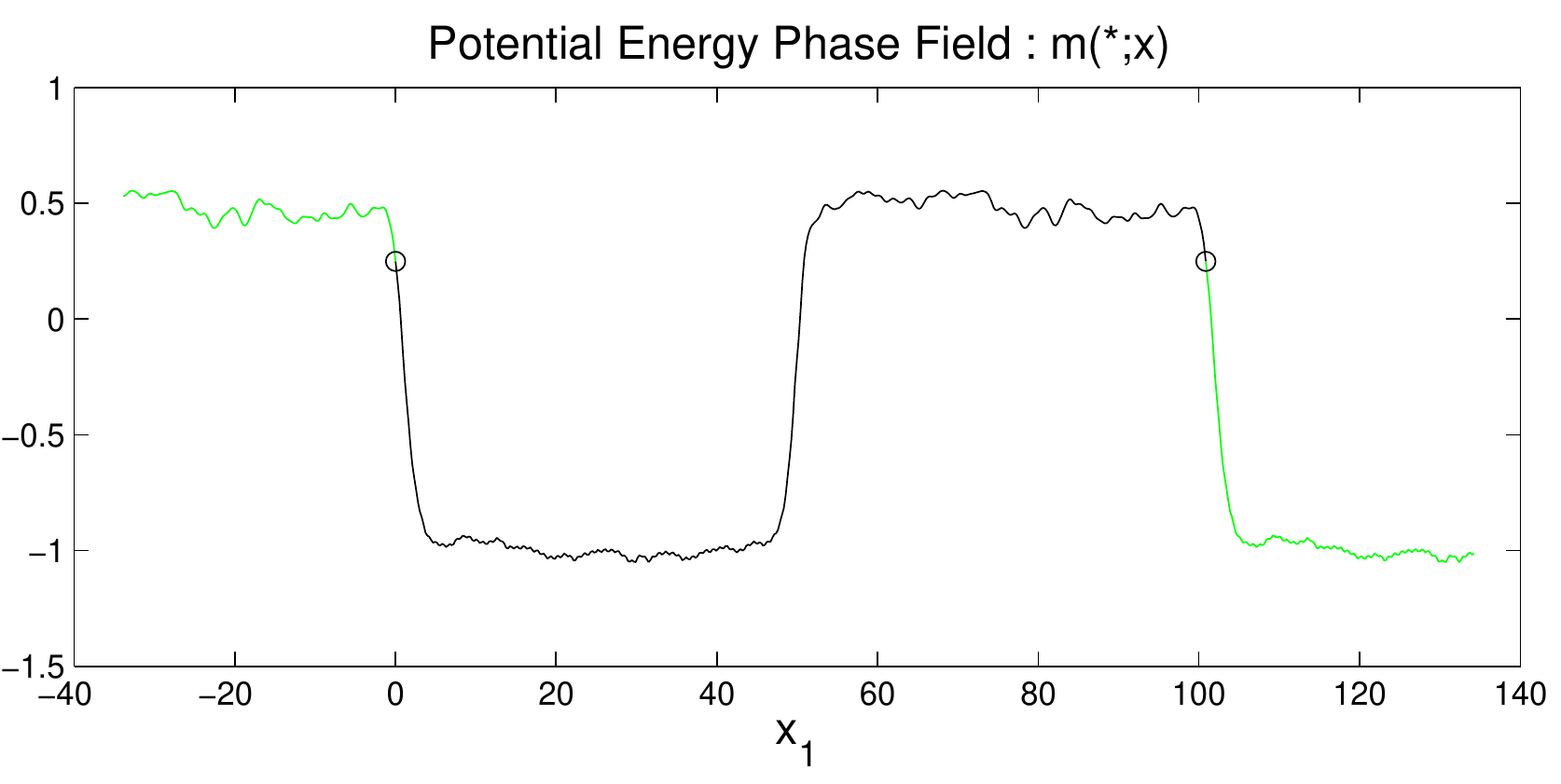}
  }
  \caption{Computed density, $\density_{loc}$, and potential
    energy phase fields for simulations O1 and O2 using
    $\molliscale=0.45$.}
  \label{fig:pf_epsi_0_45}
\end{figure}

\begin{figure}[hbp]
  \centering
  \subfigure[Orientation 1]{
    \label{fig:pf_rho_O1_epsi_0_70}
    \includegraphics[width=6.5cm]
    {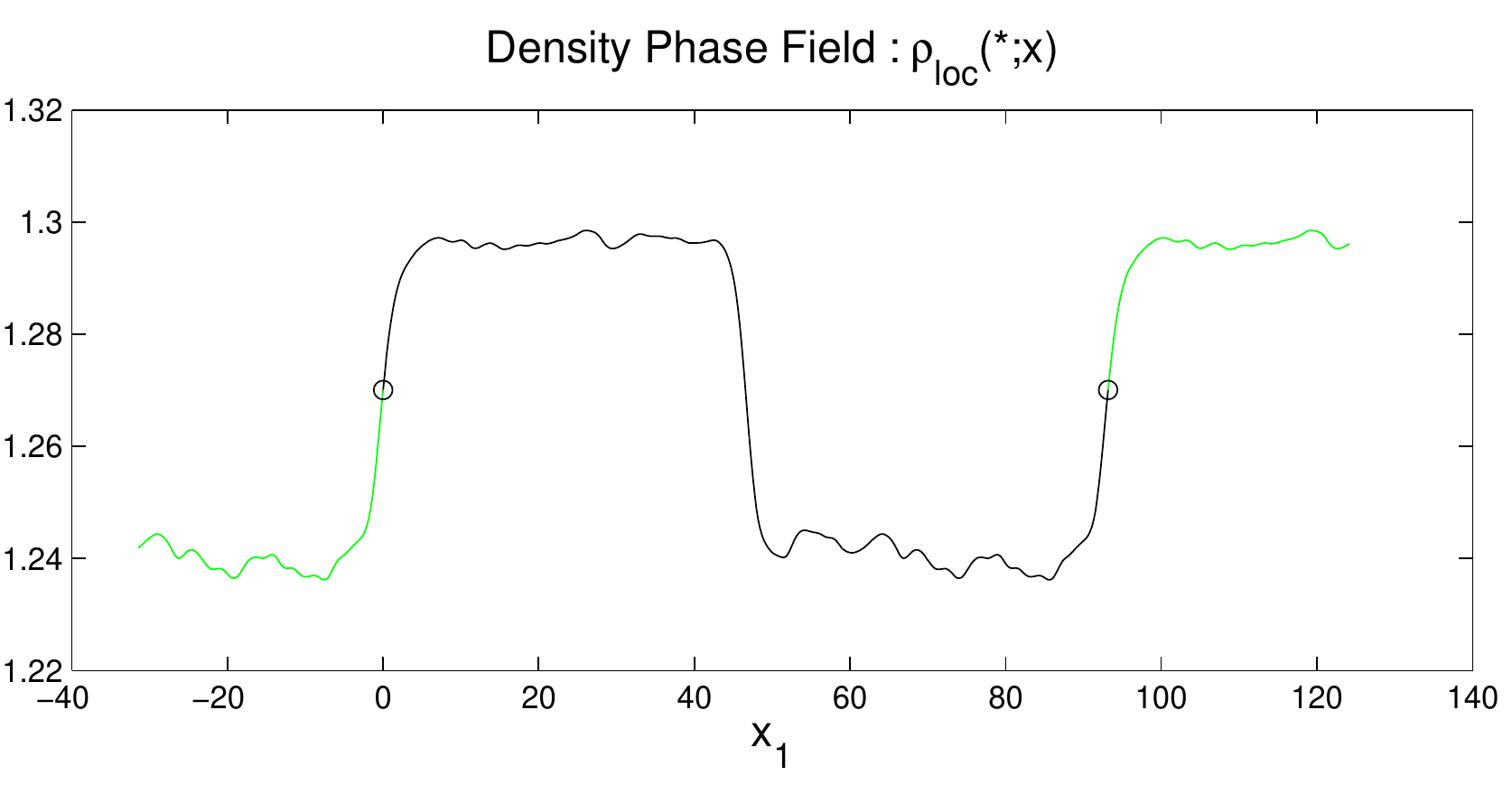}
  }
  \subfigure[Orientation 2]{
    \label{fig:pf_rho_O2_epsi_0_70}
    \includegraphics[width=6.5cm]
    {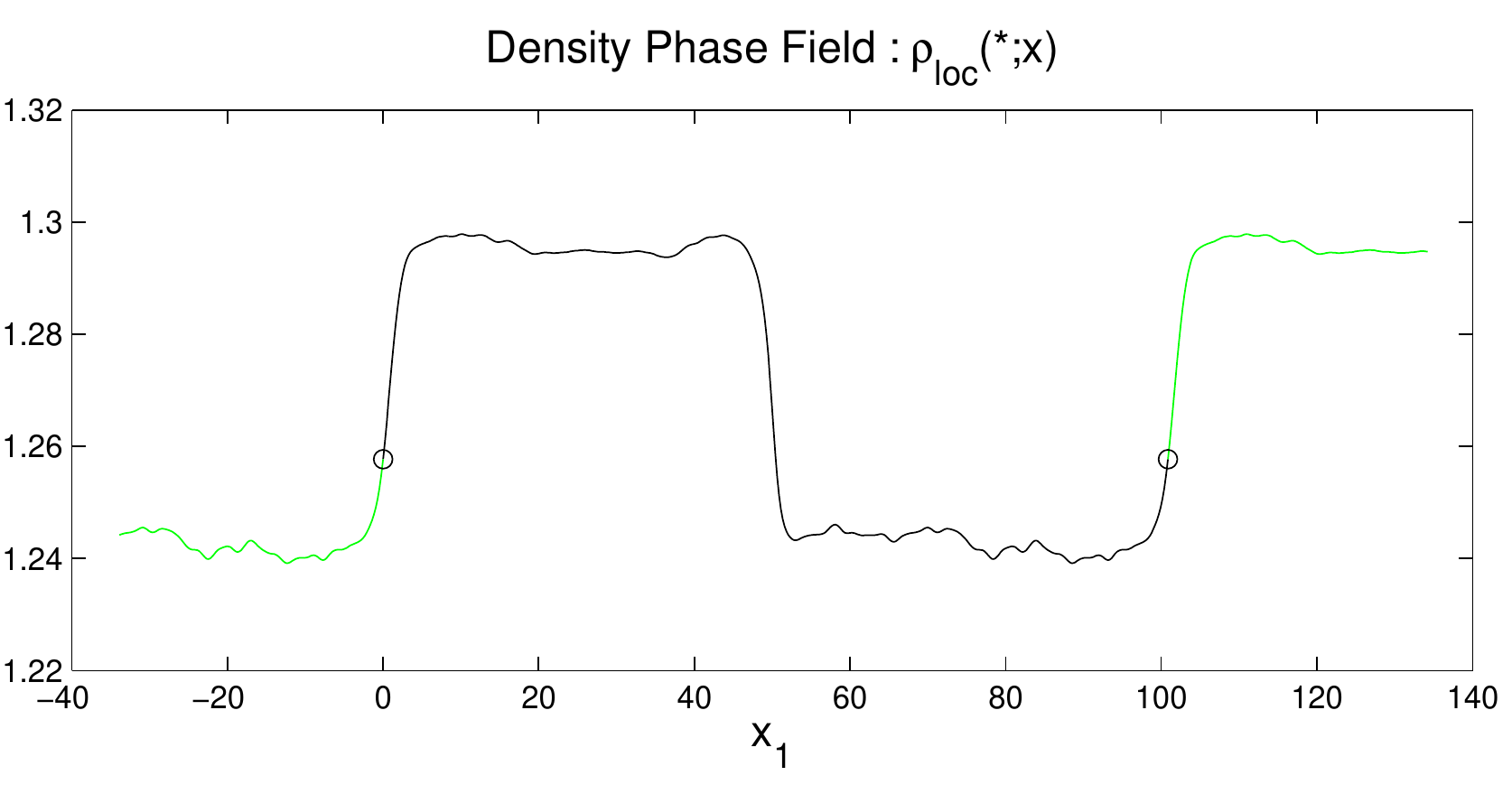}
  }
  \subfigure[Orientation 1]{
    \label{fig:pf_pot_O1_epsi_0_70}
    \includegraphics[width=6.5cm]
    {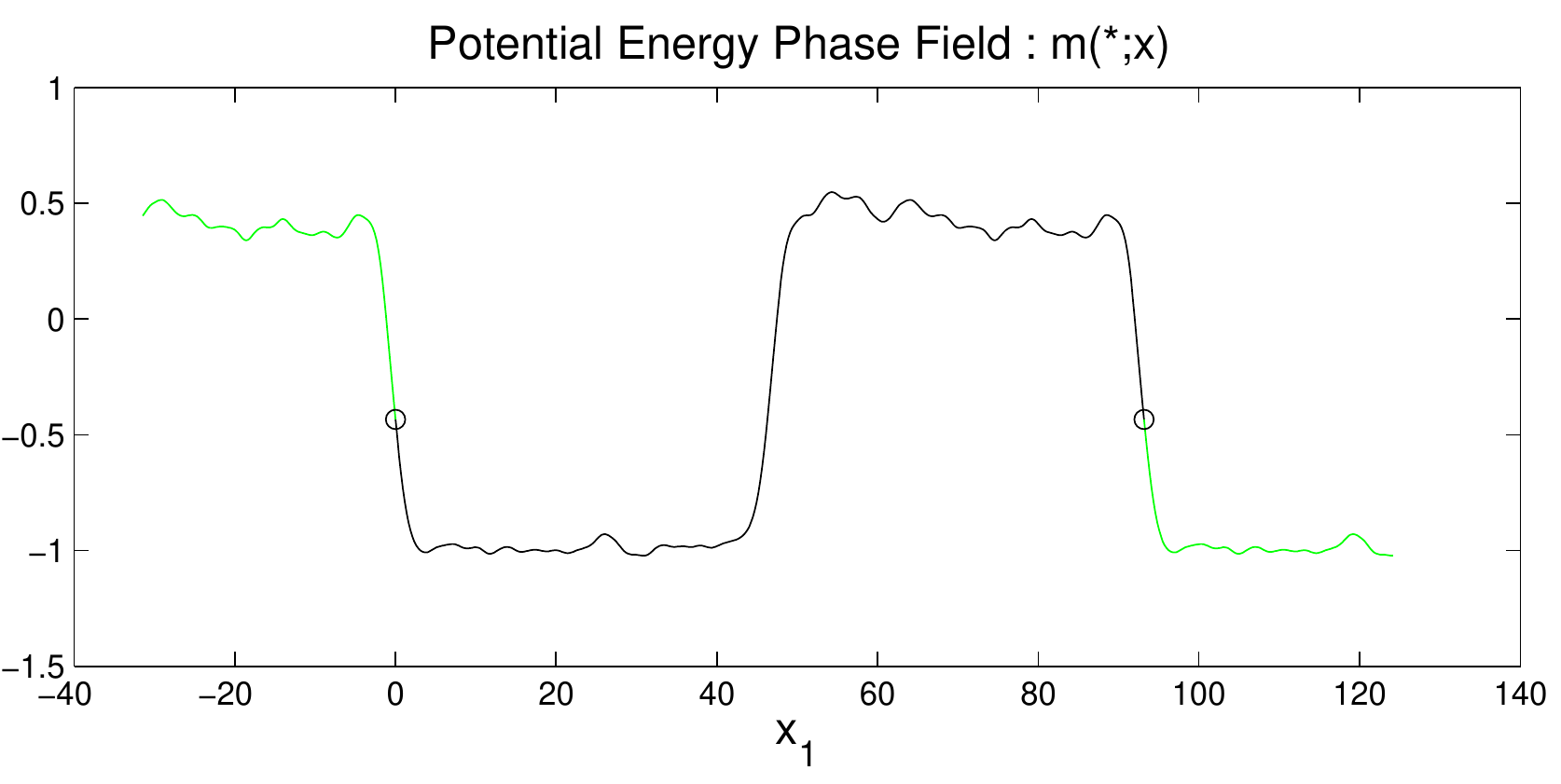}
  }
  \subfigure[Orientation 2]{
    \label{fig:pf_pot_O2_epsi_0_70}
    \includegraphics[width=6.5cm]
    {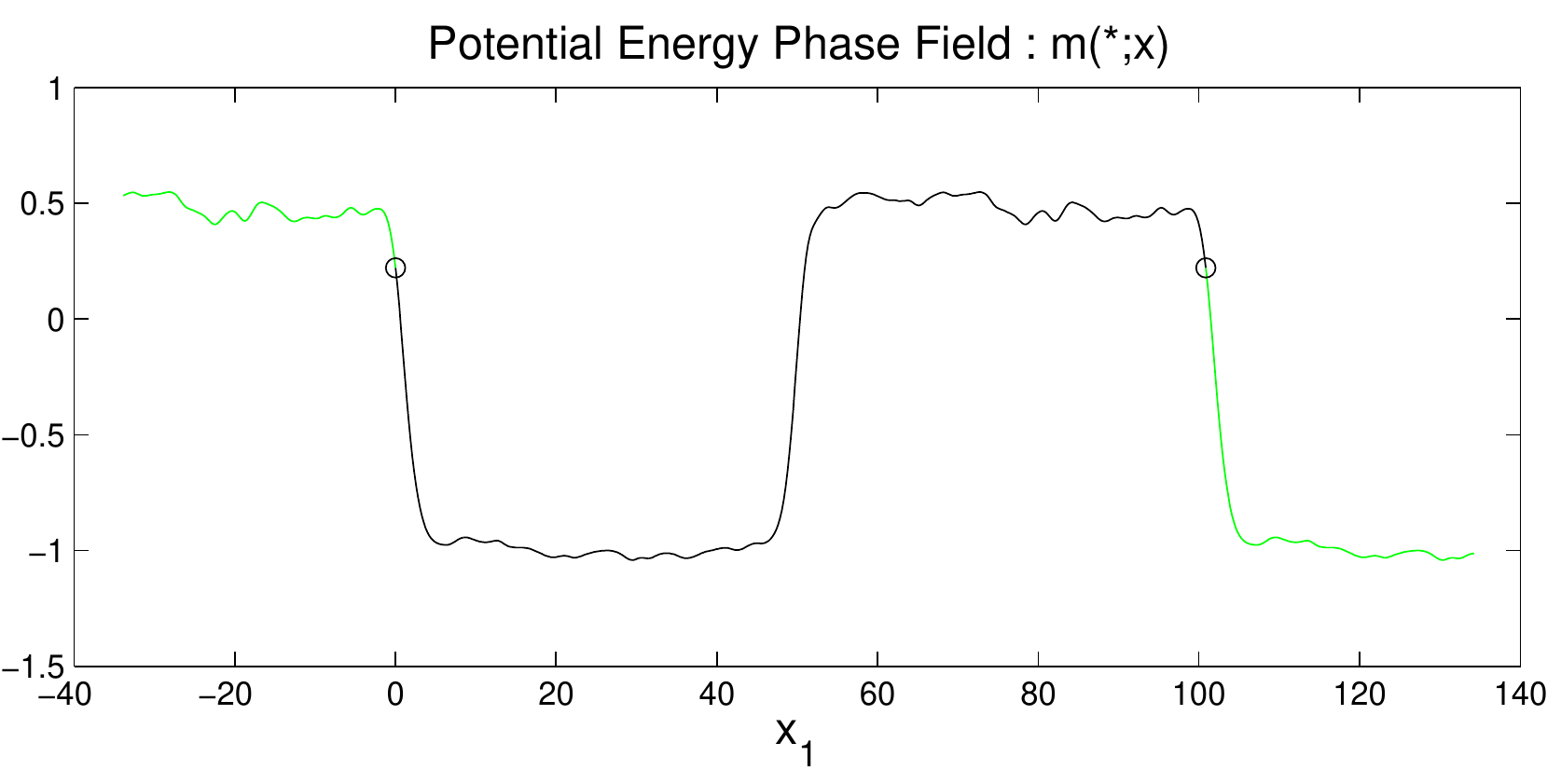}
  }
  \caption{Computed density, $\density_{loc}$, and potential
    energy phase fields for simulations O1 and O2 using
    $\molliscale=0.70$.}
  \label{fig:pf_epsi_0_70}
\end{figure}

A lower limit on
\molliscale\ is set by the demand that the phase-field be
approximately constant in the solid in spite of the periodic
structure. 
If the solid structure is aligned with the 
computational domain in such a way that the global spatial
averages are taken parallel to atomic layers, then the parameter 
\molliscale\ controlling the width of the average in the
orthogonal direction must be large enough to smooth the gaps
between the atomic layers.
In the numerical simulations the orientations of the FCC lattice
with respect to the solid--liquid interface, and hence the planes
of averaging, are precisely such that averages are computed
parallel to atomic planes, as illustrated in
Figure~\ref{fig:atom_planes}. 
In the present case the distance to the nearest neighbours in the
FCC-lattice is around 1.02; with \molli\ on the
form~\eqref{eq:mollifier} the parameter \molliscale\ must be
taken greater than $0.43$ to ensure that \molli\ decreases with
at most a factor $1/2$ in half the distance to the nearest
neighbour, which seems a reasonable demand.
Figure~\ref{fig:pf_epsi_0_45}, presenting computed phase-fields
based on local averages of the density and the potential energy
using $\molliscale=0.45$, shows that the smoothing parameter has
to be larger than this to avoid oscillations in the solid part.
The phase-fields based on $\molliscale=0.70$ in
Figure~\ref{fig:pf_epsi_0_70} do not show these oscillations on
the length scale smaller than the distance between atom layers.


\begin{figure}[hbp]
  \centering
  \begin{picture}(100,120)(0,0)
    \put(0,0){\makebox(0,0)}
  \end{picture}
  \hspace{0.2cm}
  \subfigure[Mollifier, \molli]{
    \label{fig:mollifier}
    \includegraphics[width=3.1cm]{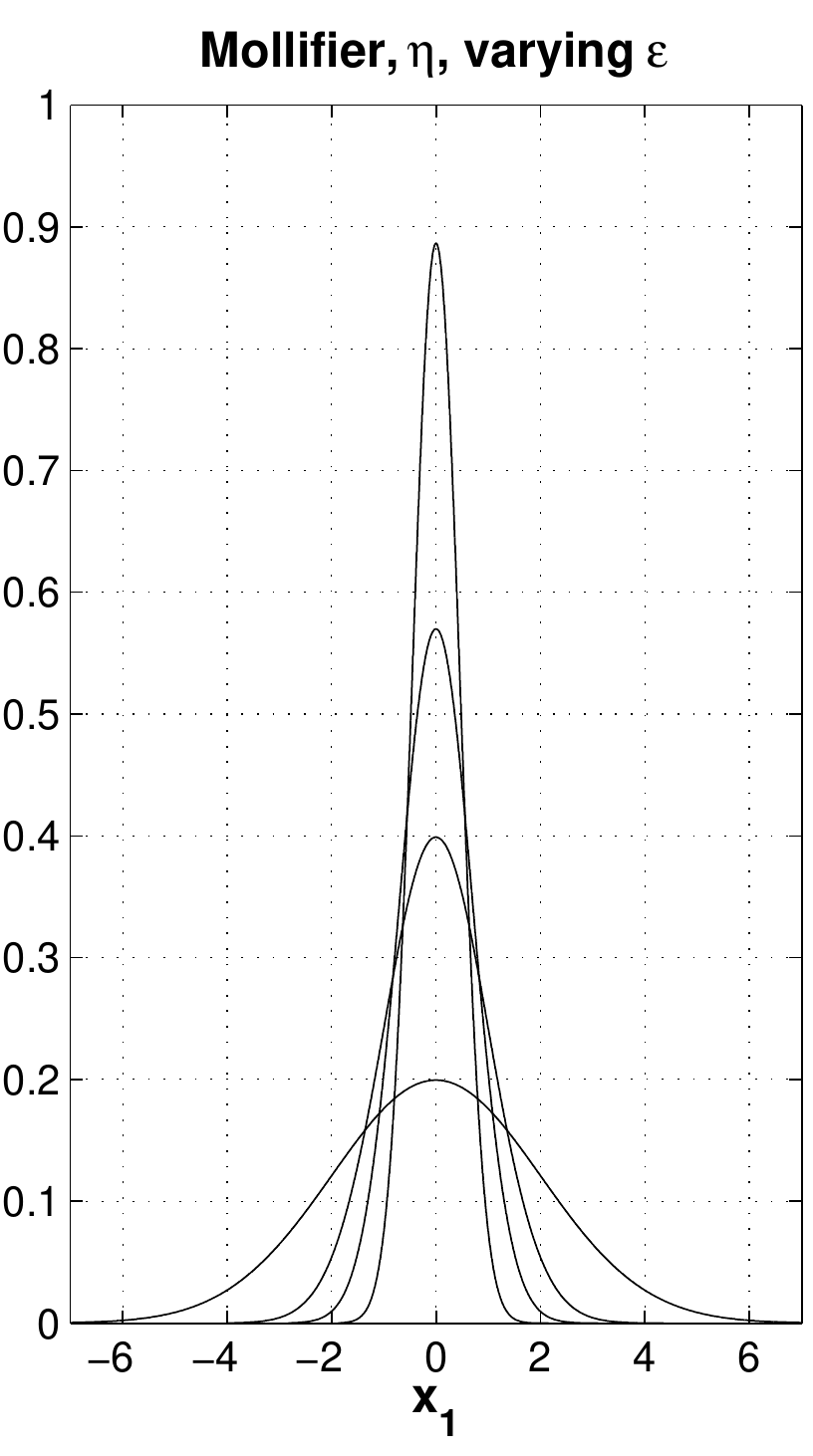}
  }
  \hspace{0.2cm}
  \begin{picture}(100,120)(0,0)
    \put(0,0){\makebox(0,0)}
  \end{picture}
  \subfigure[Orientation 1]{
    \label{fig:if_1_1}
    \includegraphics[width=3.1cm]{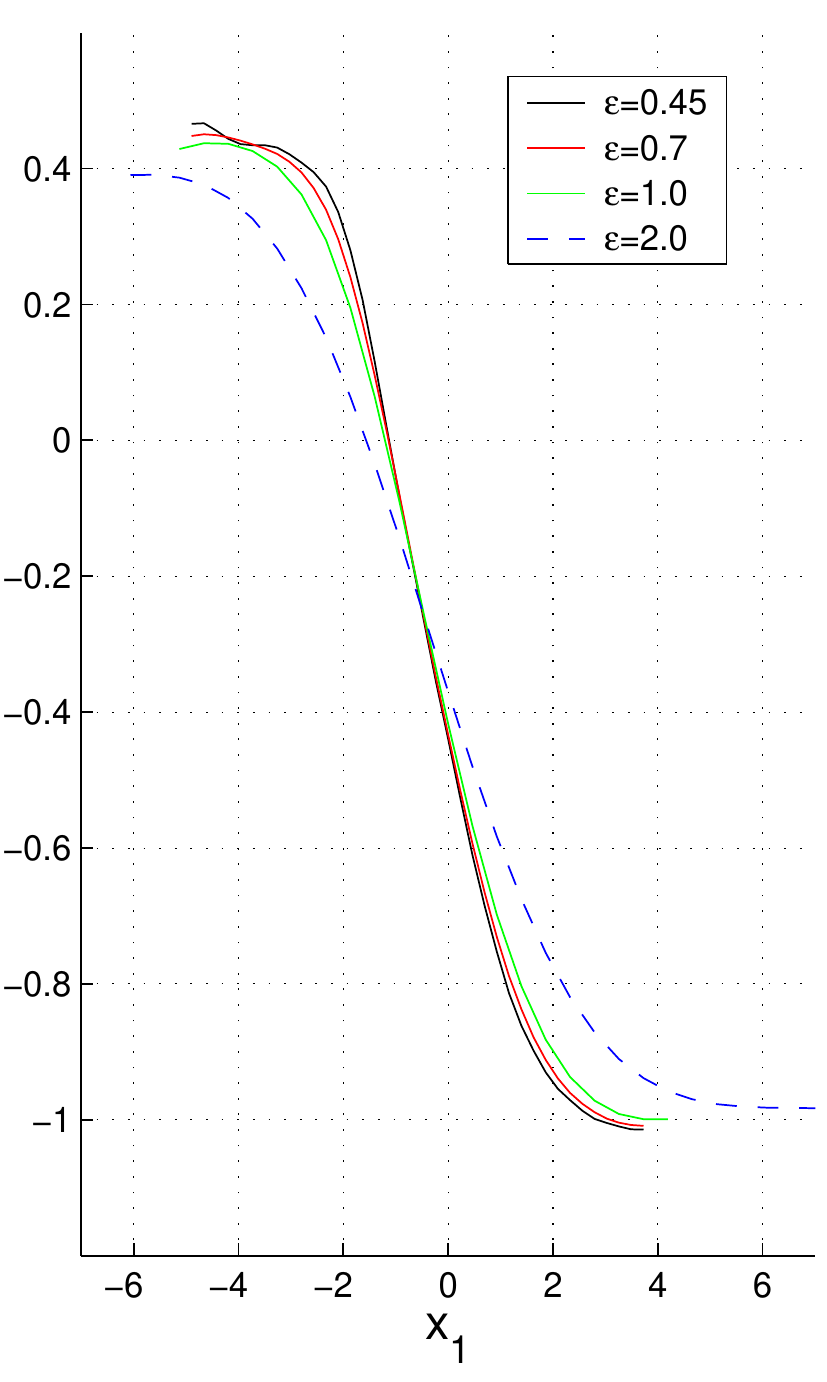}
  }
  \subfigure[Orientation 1]{
    \label{fig:if_2_1}
    \includegraphics[width=3.1cm]{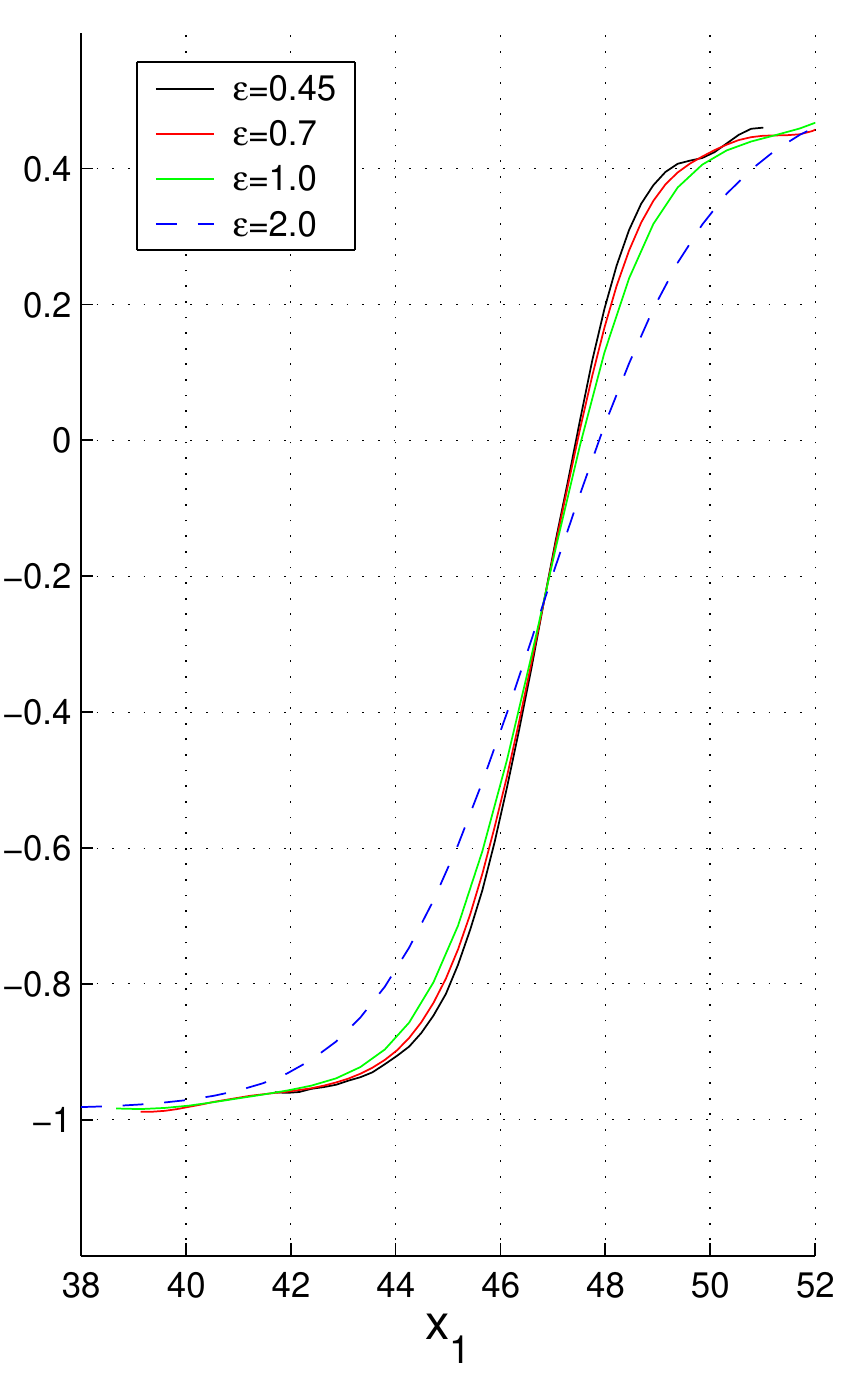}
  }
  \subfigure[Orientation 2]{
    \label{fig:if_1_2}
    \includegraphics[width=3.1cm]{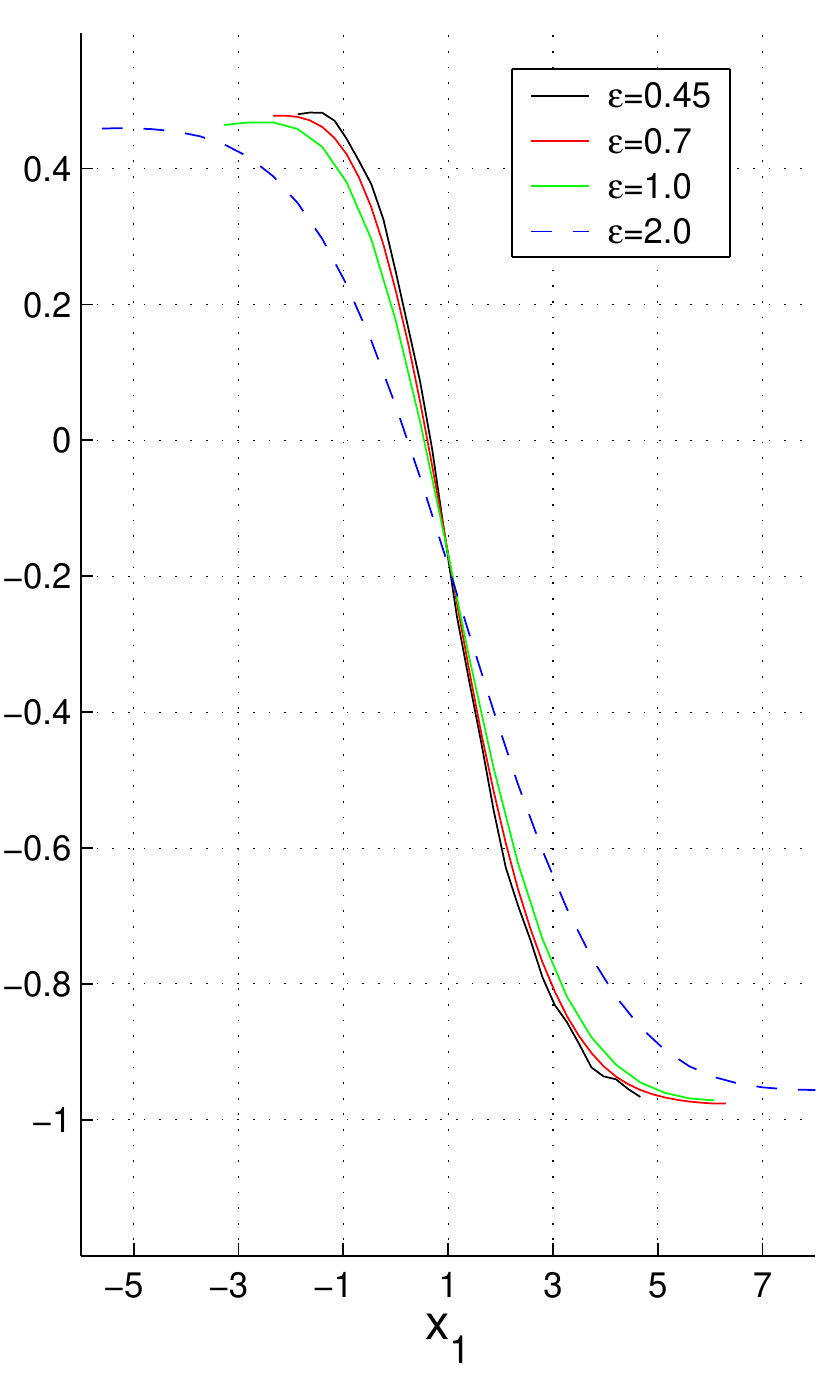}
  }
  \subfigure[Orientation 2]{
    \label{fig:if_2_2}
    \includegraphics[width=3.1cm]{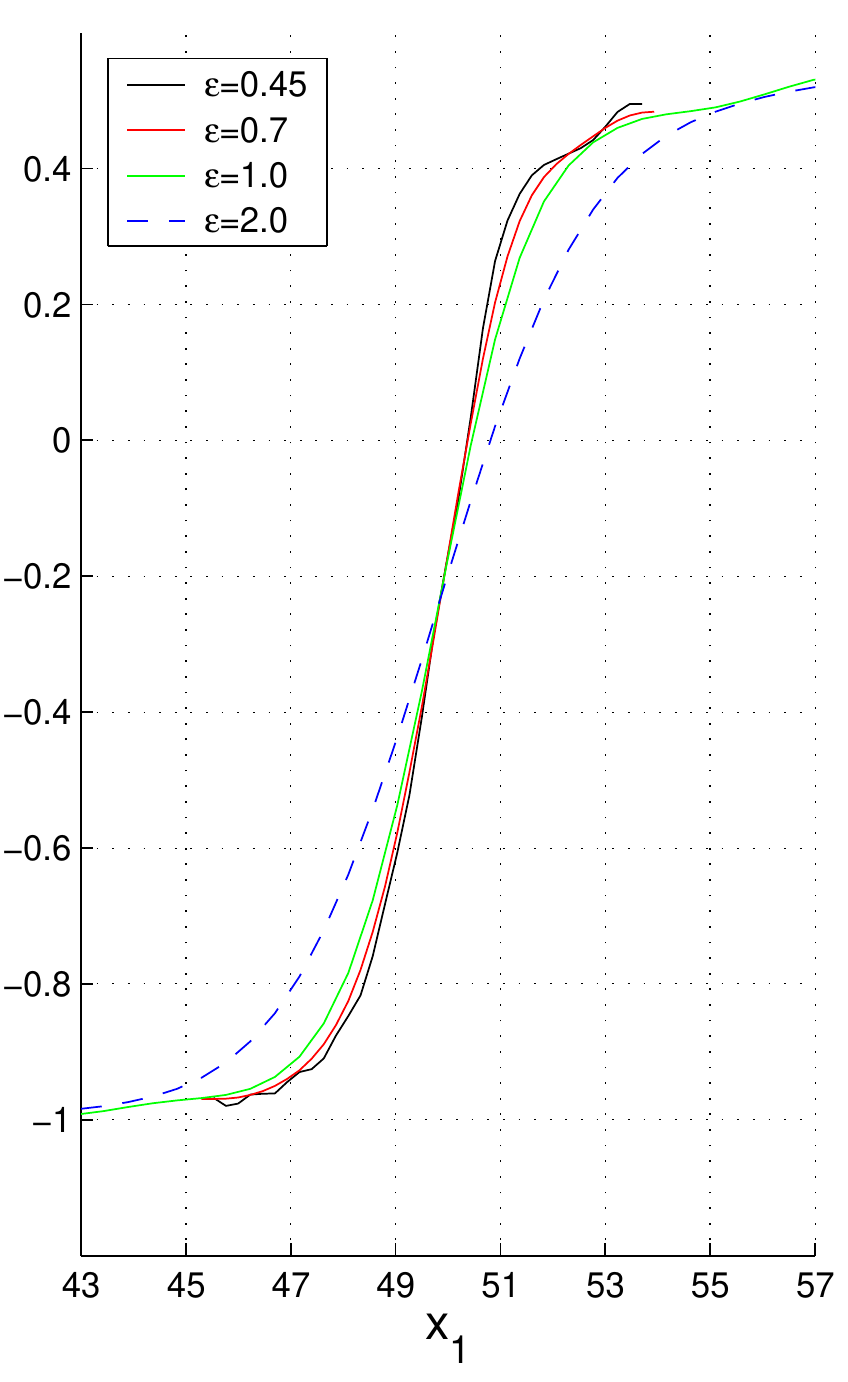}
  }  
  \caption{The mollifier, $\molli$, in the definition of the
    phase field, \pfen, depends on the model parameter
    \molliscale. The width of the averaging is proportional to
    \molliscale, as illustrated in~(a) which shows \molli\ for
    $\molliscale=0.45,\,0.7,\,1.0,\,2.0$. 
    \newline
    The phase field, \pfen, in the interface regions has been
    computed from 174 configurations with the four
    \molliscale-values listed above. 
    In (b) and (c) the configurations are taken from 
    simulation O1, and in (d) and (e) from simulation O2.  
    In each case the time interval between two successive 
    configurations is $2.5\cdot10^{-3}$, corresponding to
    $5\cdot10^3$ time steps. 
    Though the interface width in the computed phase-fields varies
    with \molliscale, it is not proportional to \molliscale\ in
    this range.
  } 
  \label{fig:interface_width}
\end{figure}

\begin{figure}[hbp]
  \centering
  \subfigure[]{
    \label{fig:comp_0_45}
    \includegraphics[width=6cm]
    {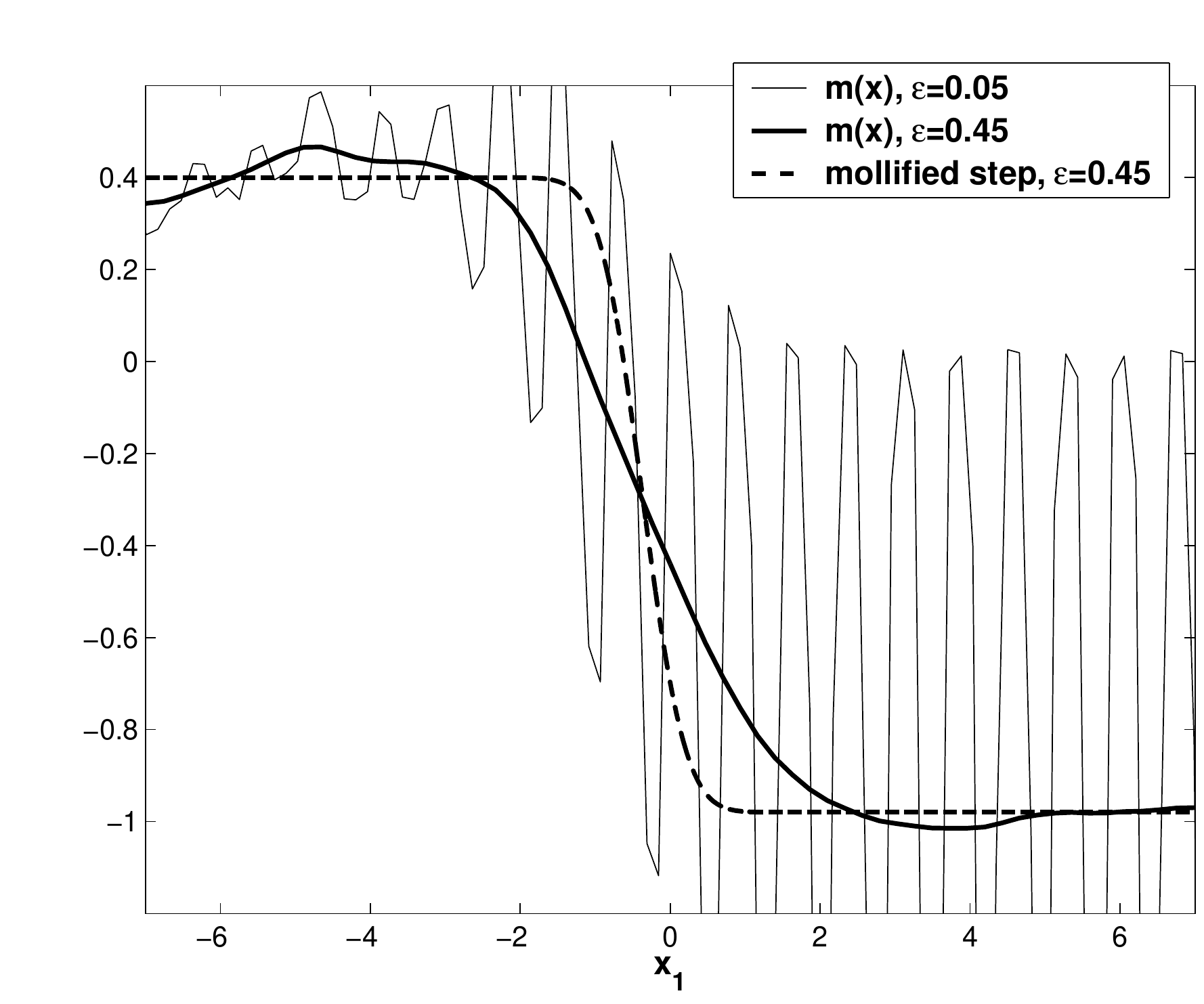}
  }
  \subfigure[]{
    \label{fig:comp_0_70}
    \includegraphics[width=6cm]
    {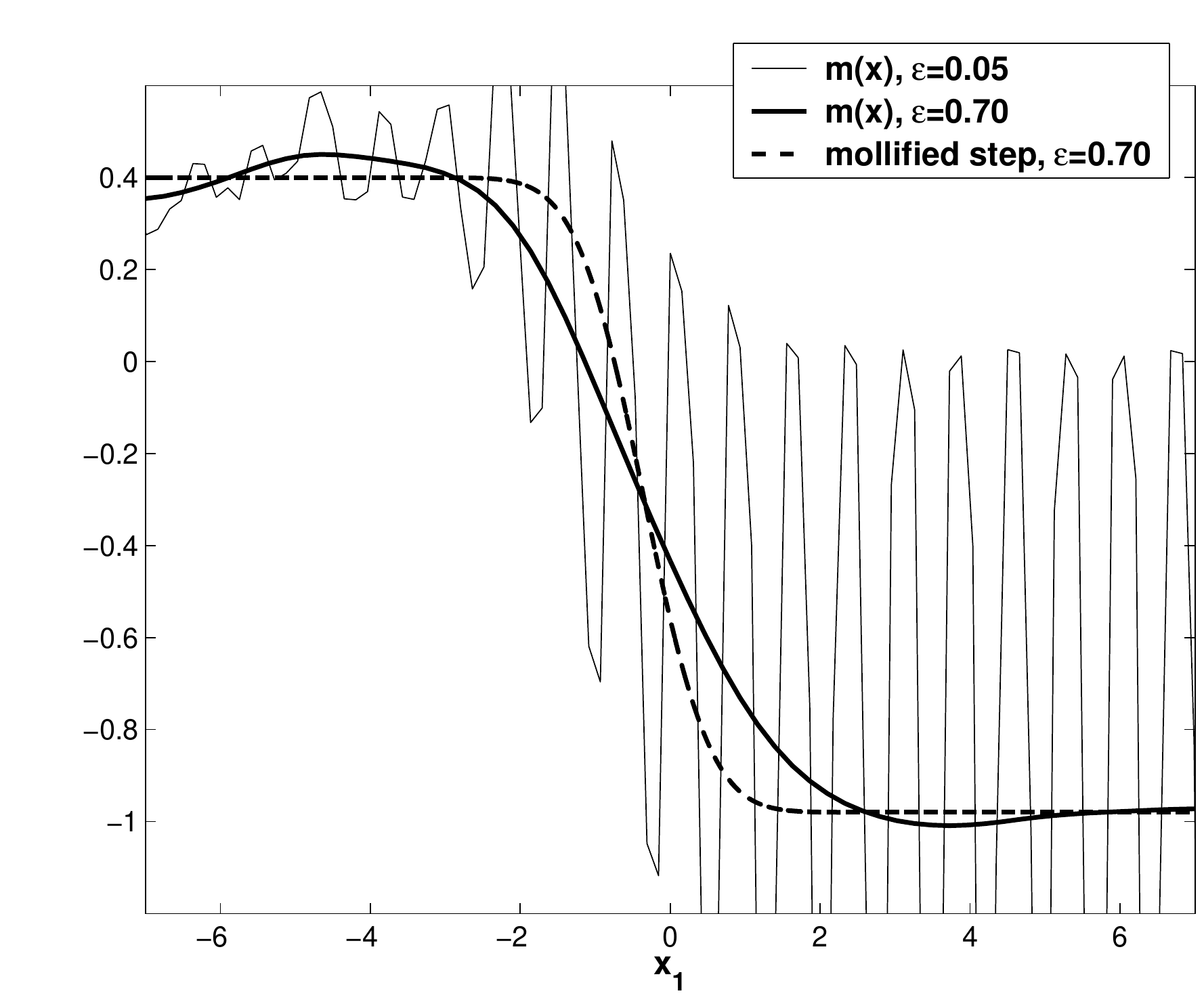}
  }
  \subfigure[]{
    \label{fig:comp_1_0}
    \includegraphics[width=6cm]
    {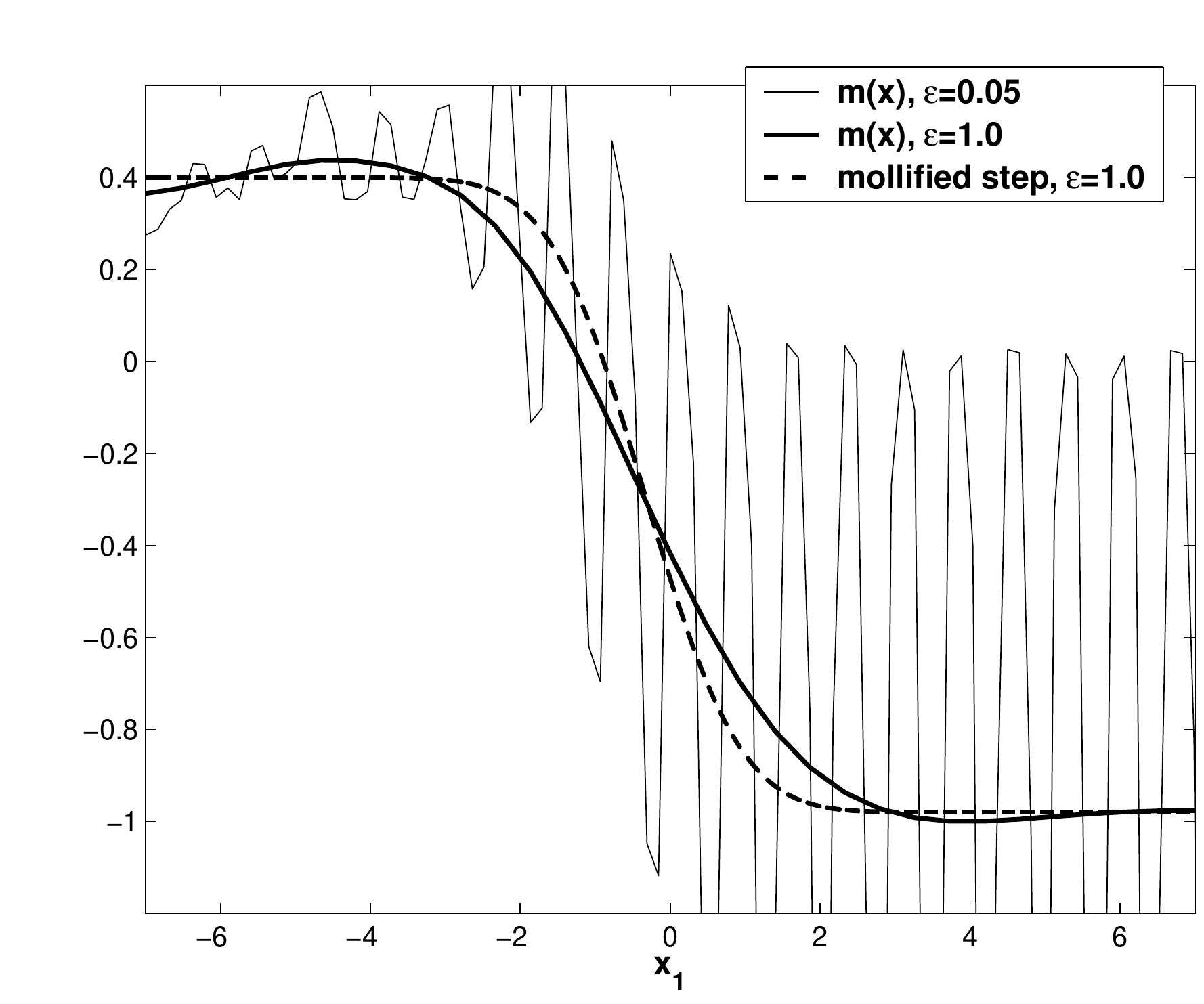}
  }
  \subfigure[]{
    \label{fig:comp_2_0}
    \includegraphics[width=6cm]
    {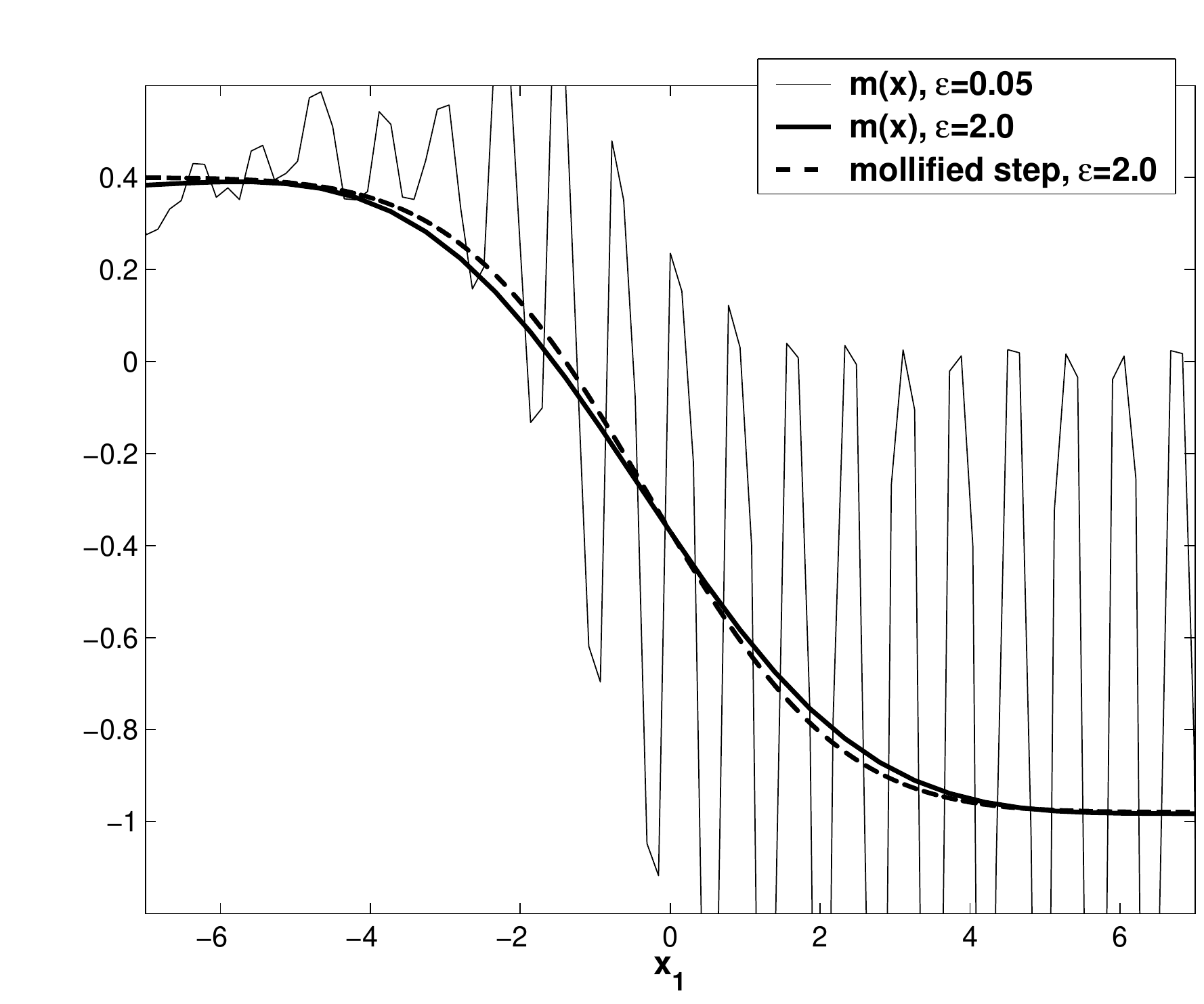}
  }
  \caption{For the phase field based on local contributions to the
    potential energy the transition from solid to liquid occurs
    on a length scale of at least several nearest neighbour
    distances for any choice of the smoothing
    parameter~\molliscale. 
    \newline
    The four subfigures are based on the same configurations from
    simulation O1 as were Figure~\ref{fig:if_1_1}--\ref{fig:if_2_1}.
    The oscillating curve present in all subfigures is the
    computed phase-field, \pfen,\ using $\molliscale=0.05$ with a
    cutoff of \molli\ at $0.3$. The nearest neighbour distance is
    approximately 1 and, for the present orientation of the FCC
    structure with respect to the $x_1$-axis, the $x_1$-distance
    between the atomic layers becomes approximately $1/\sqrt{2}$.
    Since the cutoff is less than half the distance between the atomic
    layers the phase-field would be exactly zero at the middle
    distance if the crystal were perfect and it is very close to
    zero here. The transition from the stable oscillation pattern
    in the solid to diminishing oscillations around the mean in
    the liquid is extended over a distance corresponding to at
    least four or five atomic layers in the solid.
    \newline
    The phase-field, \pfen, for
    $\molliscale=0.45,\,0.70,\,1.0,\,\mathrm{and}\,2.0$ is shown
    as the heavy solid curve in subfigures~(a)--(d). For
    reference the convolutions
    $\int_{-\infty}^{\infty}f(y)\molli(x-y)\;dy$ 
    of a sharp interface, given by the step function
    $f(y)=m_{\mathrm{liq}}\I{{\R^-}}(y)
    -m_{\mathrm{FCC}}\I{{\R^+}}(y)$, 
    and the mollifier using the respective
    \molliscale-value is included as the heavy dashed curve.
    For the smaller \molliscale-values the mollified step
    function is significantly sharper than the corresponding
    phase-field.
  }
  \label{fig:comp_to_step}
\end{figure}

For the method to be reasonable, the lower bound on \molliscale\
must not hide an interface width in the phase-field that is sharp
even on the atomic scale.
In addition to the computations with $\molliscale=1.0$, 
the phase field has been computed for
$\molliscale=0.45,\,0.70,\,\mbox{and}\,2.0$. 
The computed phase-fields in the regions around the interfaces,
for both orientation~1 and~2, are shown in 
Figure~\ref{fig:interface_width}.
The comparison shows that the interface width varies with the
smoothing parameter. It would not, however, become infinitely
sharp in the limit when \molliscale\ goes to zero, even if the
lower bound on \molliscale\ were disregarded. This is clear from
the results presented in Figure~\ref{fig:comp_to_step} where, in
addition to the values of \molliscale\ above, a phase-field
obtained with $\molliscale=0.05$, violating the lower bound, is
shown around one of the interfaces in O1. 
This value of the smoothing parameter, and the corresponding
mollifier cutoff, $R_c=6\cdot0.05=0.3$, is so small that the 
contribution to the phase-field of an individual atom in the FCC
lattice is restricted to an interval extending less than half way
to the next atom layer in either direction. Still the change in
the phase-field, from strong oscillations in the solid to decaying 
oscillations around the average in the pure liquid, occurs gradually
on a length scale corresponding to at least several atom layers
and thus several times the artificial smoothing introduced by
\molliscale. 
Figure~\ref{fig:comp_to_step} also shows that the interface
region of the phase-field obtained with
$\molliscale=0.45,\,0.70,\,\mbox{and}\,1.0$ is wider than the
transition region of a step function, representing an infinitely
sharp interface, smoothed by a convolution with the mollifier
using the corresponding \molliscale. 
For $\molliscale=2.0$ the interface is very close to that of a
mollified step function in both width and profile. 
The interface width of the smoothed step
function is proportional to \molliscale\ and it is expected that 
the same will hold for the phase-field, \pfen,
if the smoothing parameter is increased beyond the present range. 

\begin{figure}[hbp]
  \centering
  \subfigure[]{
    \label{fig:if_scale_1}
    \includegraphics[width= 6.5cm]
    {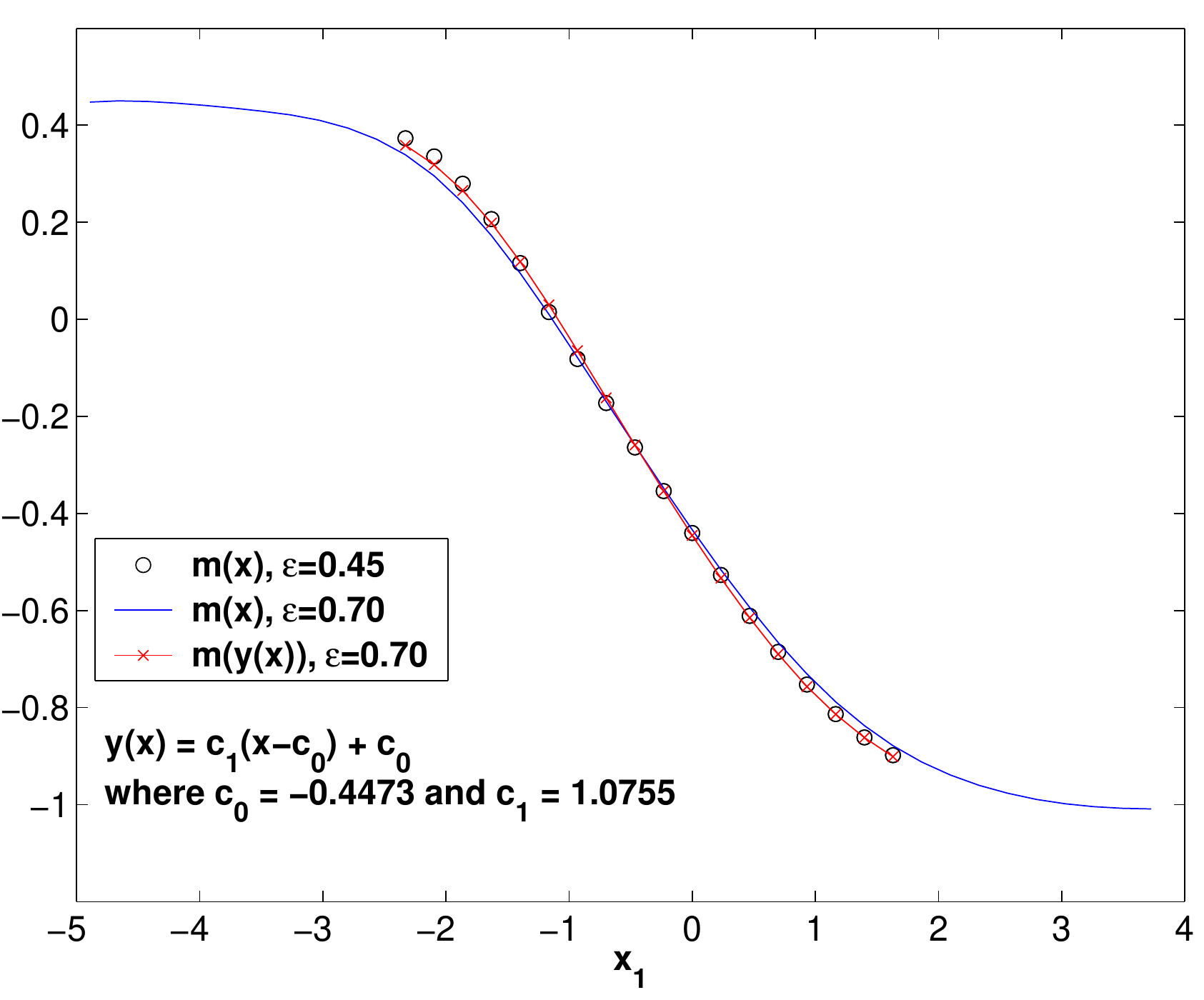}
  }
  \subfigure[]{
    \label{fig:if_scale_2}
    \includegraphics[width= 6.5cm]
    {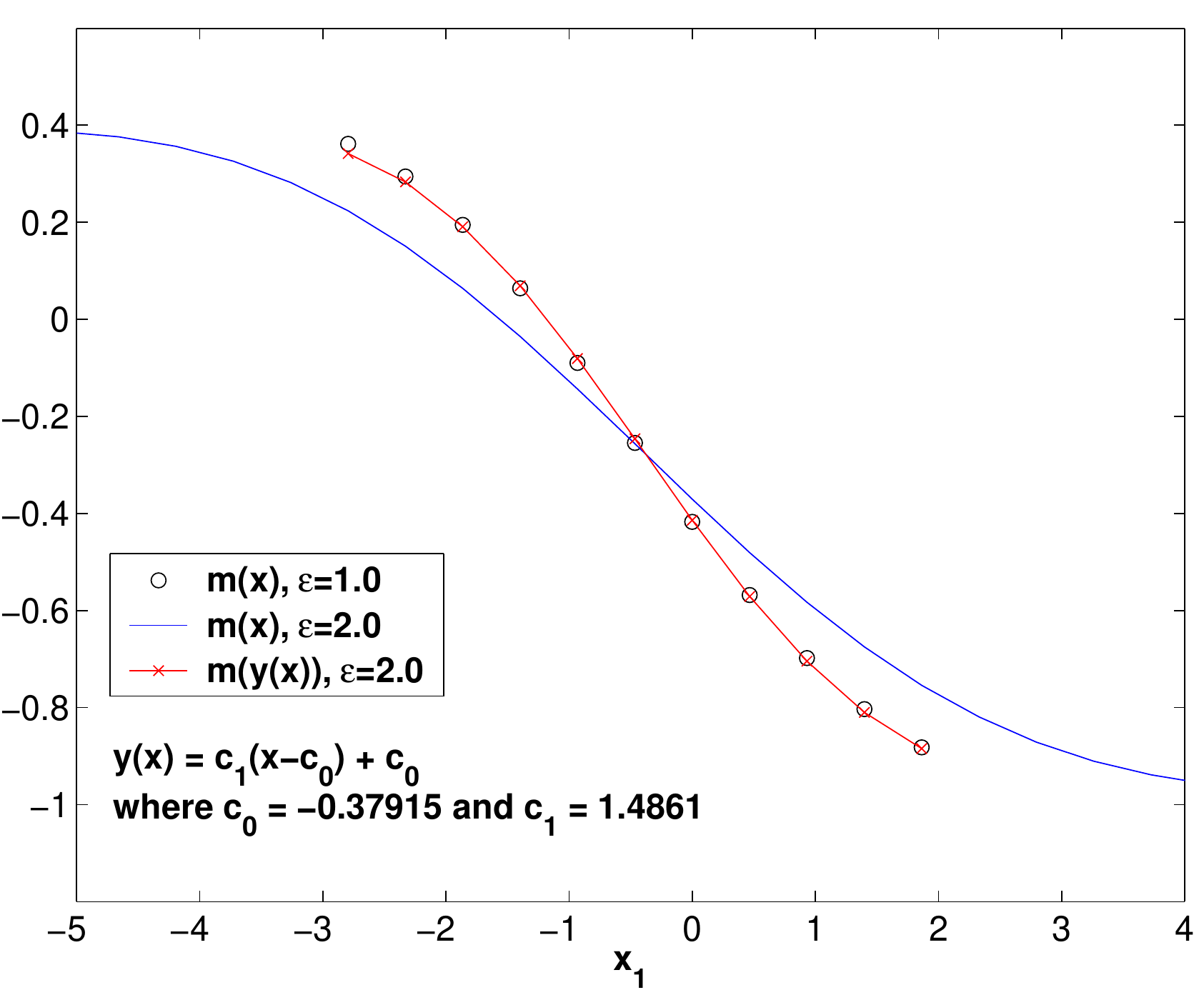}
  }
  \subfigure[Rescaling factors -- the accuracy is approximately
    $\pm 0.05$.]{
    \begin{tabular}{|l|l|c|c|c|}
      \hline
      & & \multicolumn{3}{|c|}{approximating \molliscale} \\
      \hline
      & & 0.70 & 1.0 & 2.0 \\
      \hline
      & 0.45 & 1.09 (1.56) & 1.23 (2.22) & 1.82 (4.44) \\
      reference \molliscale & 0.70 & & 1.13 (1.43) & 1.68 (2.86) \\
      & 1.0 & & & 1.49 (2.00) \\
      \hline
    \end{tabular}
  }
  \caption{For an interface given by the convolution of a sharp
    step function and the mollifier, as in
    Figure~\ref{fig:comp_to_step}, the interface width
    is directly proportional to \molliscale, since interface
    profiles corresponding to different \molliscale\ are
    identical up to affine coordinate transformations around the
    interface, $x_{\mathrm{if}}$, that is:  
    $ \pfgen_{\molliscale_2}
    (\frac{\molliscale_2}{\molliscale_1}(x-x_{\mathrm{if}})
    + x_{\mathrm{if}}) = 
    \pfgen_{\molliscale_1}(x) $.
    On a sufficiently large scale the same scaling of the
    interface width can be expected from the phase-field~\pfen\ 
    obtained from MD simulations. This is not the case when
    \molliscale\ is of the order of the nearest neighbour
    distance; then the interface width grows more slowly than
    \molliscale. 
    One way to quantify this statement is to consider the
    tabulated phase-field, $\pfen_{\molliscale_1}$, using the
    parameter value~$\molliscale_1$, 
    as given data to be approximated by the phase-field,
    $\pfen_{\molliscale_2}$,  
    based on the parameter value~$\molliscale_2$; the allowed
    approximations use 
    affine coordinate transformations
    $y(x) = c_1(x-c_0) + c_0$
    of the independent coordinate.
    The data points $((x_k),\pfen_{\molliscale_1}(x_k))$ are
    taken from the interior of an interface, 
    $\pfen_{\mathrm{solid}} < m_0 \leq 
    \pfen_{\molliscale_1}(x_k) \leq m_1 <
    \pfen_{\mathrm{liquid}}$, and the function
    $\pfen_{\molliscale_2}$ is defined by linear interpolation
    between tabulated values.
    A least squares approximation of the overdetermined system
    $\pfen_{\molliscale_2}(y(x_k))=\pfen_{\molliscale_1}(x_k)$
    for $c_0$ and $c_1$ 
    gives a value of the scaling factor $c_1$ to be compared to 
    $\molliscale_2/\molliscale_1$.
    \newline
    Subfigures (a) and (b) show two examples for the
    interface in Figure~\ref{fig:if_1_1}. The circles, $\circ$,
    denote the reference data points, the solid line shows the
    linear interpolation of the tabulated values for the
    approximating phase-field, and the line marked with crosses,
    $\times$, is the least square approximation. 
    \newline
    The table (c) shows the scaling constants obtained after
    averaging over all four interfaces in
    Figure~\ref{fig:if_1_1}--\ref{fig:if_2_2}. 
    The corresponding
    quotients~$\frac{\molliscale_2}{\molliscale_1}$ are included
    in parenthesis for reference.
  }
  \label{fig:interface_scaling}
\end{figure}

\appendix

\include{calculations}

\bibliographystyle{nada-en}
\bibliography{bibfil}

\end{document}

%% file: calculations.tex
\section{Explicit Calculation of Drift and Diffusion Functions}
\label{sec:calculations}

Let the total potential energy be
\begin{align}
  \totpot(\mdpos{\no}) & = \sumall{i} \potenp{i}(\mdposgen), 
  \nonumber
  \intertext{where }
  \potenp{i}(\mdposgen) & = 
  \frac{1}{2}\sumneq{k}{i} 
  \pairpot(||\mdpos[\no]{i}-\mdpos[\no]{k}||). 
  \nonumber
  \intertext{For the phase-field}
  \pfen(x;\mdposgen) 
  & = \sumall{i} 
  \potenp{i}(\mdposgen) \molli(x-\mdpos[\no]{i}), 
  \nonumber
  \intertext{where the particle positions $\mdposgen\in\R^{3\nrp}$ solve the
    \Ito\ stochastic differential equation }
  d\mdpos{\no} & = -\nabla_{\mdposgen} \totpot(\mdpos{\no})\;dt 
  + \sqrt{2\kb \temp}\;dW^t,
  \nonumber
  \intertext{\Ito's formula gives}
  d\pfen(x;\mdpos{\no}) 
  & = 
  \sumall{j}\driftmd_j(x;\mdpos{\no})\,dt +
  \sumall{j}\sum_{k=1}^3\diffumd_{j,k}(x;\mdpos{\no})\,dW_{j,k}^t,
  \nonumber
  \intertext{with}
  \label{eq_app:driftmd_j}
  \driftmd_j(x;\mdposgen) 
  & =
  -\dX\pfen(x;\mdposgen)\cdot\dX \totpot(\mdposgen) 
  + \kb\temp \dX\!\cdot\!\dX \pfen(x;\mdposgen)
  \intertext{and}
  \label{eq_app:diffumd_j}
  \diffumd_{j,\cdot}(x;\mdposgen) 
  & =
  \sqrt{2\kb\temp}\dX\pfen(x;\mdposgen).
\intertext{Introducing the total force, $\force_j$, acting on particle
  $j$, and the contributions from individual pairs, $\pairforce_{ij}$, }
  \force_j(\mdposgen)   & = 
  - \dX \totpot(\mdposgen)
  = \sumneq{i}{j} \pairforce_{ij}(\mdposgen),
  \nonumber
  \\
  \pairforce_{ij}(\mdposgen) 
  & = \pairpot'(|| \mdpos[\no]{i}-\mdpos[\no]{j} ||) 
  \frac{\mdpos[\no]{i}-\mdpos[\no]{j}}
  {|| \mdpos[\no]{i}-\mdpos[\no]{j} ||},
  \nonumber
\end{align}
the gradient of $\potenp{i}$ with respect to the position of 
particle $j$ is
\begin{align*}
  \dX \potenp{i}(\mdposgen) & = \nonumber
  \frac{1}{2}\sumneq{k}{i} 
  \dX \pairpot(||\mdpos[\no]{i}-\mdpos[\no]{k}||)
  \\ & = \nonumber
  \delta_{ij}\frac{1}{2}\sumneq{k}{j} 
  \dX\pairpot(||\mdpos[\no]{j}-\mdpos[\no]{k}||)
  + (1-\delta_{ij})\frac{1}{2}
  \dX\pairpot(||\mdpos[\no]{i}-\mdpos[\no]{j}||)
  \\ & = -\delta_{ij}\frac{1}{2}\force_j(\mdposgen) 
  - (1-\delta_{ij}) \frac{1}{2} \pairforce_{ij}(\mdposgen),
\end{align*}
where $\delta_{ij}$ is the Kronecker delta: 
$\delta_{ij}=1,\text{if $i=j$}, \delta_{ij}=0,\text{if $i\neq j$}$. 
The gradient of the phase-field variable with respect to the position
of particle $j$ is
\begin{align*}
  \dX \pfen(x;\mdposgen) & = 
  \potenp{j}(\mdposgen) \dX \molli(x-\mdpos[\no]{j}) + 
  \sumall{i}\dX \potenp{i}(\mdposgen) \molli(x-\mdpos[\no]{i})  
  \\ & = 
  - \potenp{j}(\mdposgen) \dx \molli(x-\mdpos[\no]{j}) 
  \\ & \quad 
  - \frac{1}{2} \sumall{i} 
  \delta_{ij}\force_j(\mdposgen) \molli(x-\mdpos[\no]{i})
  - \frac{1}{2} \sumall{i} 
  (1-\delta_{ij}) \pairforce_{ij}(\mdposgen) \molli(x-\mdpos[\no]{i})
  \\ & = 
  - \dx (\potenp{j}(\mdposgen) \molli(x-\mdpos[\no]{j})) 
  \\ & \quad
  - \frac{1}{2}\force_j(\mdposgen)\molli(x-\mdpos[\no]{j})
  -\frac{1}{2} \sumneq{i}{j} 
  \pairforce_{ij}(\mdposgen)\molli(x-\mdpos[\no]{i}).
\end{align*}
Introducing the notation $-\divforce_j$ for the divergence of the
force $\force_j$ with respect to $\mdpos[\no]{j}$ and the notation 
$\divpairforce_{ij}$ for the individual contributions,
\begin{align*}
  \divforce_j(\mdposgen) 
  & = - \dX \cdot \force_j(\mdposgen) 
  = 
  - \sumneq{i}{j} \dX \cdot \pairforce_{ij}(\mdposgen)
  = \sumneq{i}{j} \divpairforce_{ij}(\mdposgen),
  \\
  \divpairforce_{ij}(\mdposgen) 
  & = \pairpot''(||\mdpos[\no]{i}-\mdpos[\no]{j}||) + 
  \pairpot'(||\mdpos[\no]{i}-\mdpos[\no]{j}||)
  \frac{2}{||\mdpos[\no]{i}-\mdpos[\no]{j}||},
\end{align*}
the divergence of gradient of phase field variable with respect to the 
position of particle $j$ becomes 
\begin{align*}
  \dX\!\cdot\!\dX \pfen(x;\mdposgen) & = \nonumber
  - \dX\cdot 
  \left( \potenp{j}(\mdposgen) \dx \molli(x-\mdpos[\no]{j}) \right)
  \\ & \quad 
  - \frac{1}{2}\dX\cdot 
  \left(\force_j(\mdposgen)\molli(x-\mdpos[\no]{j}) \right)
  - \frac{1}{2}\sumneq{i}{j} \dX\cdot 
  \left(\pairforce_{ij}(\mdposgen)\molli(x-\mdpos[\no]{i}) \right)
  \\ & = 
  -\dX \potenp{j}(\mdposgen)\cdot \dx \molli(x-\mdpos[\no]{j}) 
  - \potenp{j}(\mdposgen) \dX \cdot \dx \molli(x-\mdpos[\no]{j}) 
  \\ & \quad 
  - \frac{1}{2}\dX\cdot \force_j(\mdposgen)\molli(x-\mdpos[\no]{j})
  - \frac{1}{2}\force_j(\mdposgen)\cdot \dX \molli(x-\mdpos[\no]{j})
  \\ & \quad 
  -
  \frac{1}{2}\sumneq{i}{j} \dX\cdot 
  \pairforce_{ij}(\mdposgen)\molli(x-\mdpos[\no]{i}) 
  \\ & = 
  \frac{1}{2} \force_j(\mdposgen) \cdot \dx \molli(x-\mdpos[\no]{j})
  + \potenp{j}(\mdposgen) \dx \cdot \dx \molli(x-\mdpos[\no]{j}) 
  \\ & \quad 
  + \frac{1}{2} \divforce_j(\mdposgen) \molli(x-\mdpos[\no]{j})
  + \frac{1}{2} \force_j(\mdposgen) \cdot \dx \molli(x-\mdpos[\no]{j})
  \\ &  \quad 
  + \frac{1}{2}\sumneq{i}{j} 
  \divpairforce_{ij}(\mdposgen)\molli(x-\mdpos[\no]{i})
  \\ & = 
  \dx\!\cdot\!\dx 
  \left( \potenp{j}(\mdposgen) \molli (x-\mdpos[\no]{j}) \right)
  + \dx \cdot \left( \force_j(\mdposgen) \molli(x-\mdpos[\no]{j}) \right)
  \\ & \quad 
  + \frac{1}{2} \divforce_j(\mdposgen) \molli(x-\mdpos[\no]{j})
  + \frac{1}{2}\sumneq{i}{j} 
  \divpairforce_{ij}(\mdposgen)\molli(x-\mdpos[\no]{i}).
\end{align*}

Using the explicit expressions for 
$\dX \pfen(x;\mdposgen)$ and $\dX\!\cdot\!\dX \pfen(x;\mdposgen)$, 
the components~\eqref{eq_app:driftmd_j} of the drift become
\begin{align*}
  \driftmd_j(x;\mdposgen) & = 
  \dx (\potenp{j}(\mdposgen) \molli(x-\mdpos[\no]{j})) 
  \cdot (-\force_j(\mdposgen))
  + \frac{1}{2}\force_j(\mdposgen)\molli(x-\mdpos[\no]{j}) 
  \cdot (-\force_j(\mdposgen))
  \\ & \quad 
  + \frac{1}{2} \sumneq{i}{j} 
  \pairforce_{ij}(\mdposgen)\molli(x-\mdpos[\no]{i})
  \cdot (-\force_j(\mdposgen))
  \\ & \quad 
  + \kb\temp \dX\!\cdot\!\dX \pfen(x;\mdposgen)
  \\ & = 
  - \dx \cdot (\potenp{j}(\mdposgen) \force_j(\mdposgen) \molli(x-\mdpos[\no]{j}))
  - \frac{1}{2} ||\force_j(\mdposgen)||^2\molli(x-\mdpos[\no]{j})
  \\ & \quad 
  - \frac{1}{2} \sumneq{i}{j} \pairforce_{ij}(\mdposgen) \cdot 
  \force_j(\mdposgen)\molli(x-\mdpos[\no]{i})
  \\ & \quad 
  + \kb\temp \dx\!\cdot\!\dx 
  \left( \potenp{j}(\mdposgen) \molli (x-\mdpos[\no]{j}) \right)
  + \kb\temp 
  \dx \cdot \left( \force_j(\mdposgen) \molli(x-\mdpos[\no]{j}) \right)
  \\ & \quad 
  + \kb\temp
  \frac{1}{2} \divforce_j(\mdposgen) \molli(x-\mdpos[\no]{j})
  + \kb\temp
  \frac{1}{2}\sumneq{i}{j} 
  \divpairforce_{ij}(\mdposgen)\molli(x-\mdpos[\no]{i}).
  \\ & = 
  \kb\temp \dx\!\cdot\!\dx
  \bigl( \potenp{j}(\mdposgen) \molli (x-\mdpos[\no]{j}) \bigr)
  \\ & \quad +
  \dx\cdot \biggl( 
    (\kb\temp - \potenp{j}(\mdposgen))\force_j(\mdposgen)
    \molli(x-\mdpos[\no]{j}) 
  \biggr)
  \\ & \quad +
  \frac{1}{2} \bigl( 
    \kb\temp \divforce_j(\mdposgen) - ||\force_j(\mdposgen)||^2
  \bigr) \molli(x-\mdpos[\no]{j})
  \\ & \quad +
  \frac{1}{2} \sumneq{i}{j} \bigl(
    \kb\temp \divpairforce_{ij}(\mdposgen) 
    - \pairforce_{ij}(\mdposgen) \cdot \force_j(\mdposgen)
  \bigr) \molli(x-\mdpos[\no]{i})
\end{align*}
so that, after summing over $j$,
\begin{align*}
  \driftmd(x;\mdposgen) & = 
  \kb\temp \dx\!\cdot\!\dx \pfen(x;\mdposgen)
  + \dx\cdot \widetilde\driftone(x;\mdposgen) 
  + \driftzero(x;\mdposgen)
  \intertext{with}
  \widetilde\driftone(x;\mdposgen) & = 
  \sumall{j}
  (\kb\temp - \potenp{j}(\mdposgen))\force_j(\mdposgen)
  \molli(x-\mdpos[\no]{j})
  \intertext{and}
  \driftzero(x;\mdposgen) & = 
  \sumall{j}
  \left(\kb\temp \divforce_j(\mdposgen) 
    - \frac{1}{2} ||\force_j(\mdposgen)||^2 \right) 
  \molli(x-\mdpos[\no]{j}) 
  \\ & \quad 
  - \frac{1}{2} 
  \sumall{j} \sumneq{i}{j}
  \pairforce_{ij}(\mdposgen) \cdot 
  \force_j(\mdposgen) \molli(x-\mdpos[\no]{i}). 
\end{align*}

Using the one-dimensional mollifier 
\begin{align}
  \label{eq_app:mollifier}
  \molli(x) & = \molli(x_1) = \text{constant}\cdot
  \exp{
    \left(-\frac{1}{2}\left(\frac{x_1}{\epsilon}\right)^2\right)
  },
\end{align}
that only varies in the $x_1$-direction, the expression for the drift
reduces to 
\begin{align*} 
  \driftmd(x;\mdposgen) & = 
  \kb\temp \ddxett \pfen(x;\mdposgen)
  + \dxett \driftone(x;\mdposgen) 
  + \driftzero(x;\mdposgen)
  \intertext{with}
  \driftone(x;\mdposgen) & = 
  \sumall{j}
  (\kb\temp - \potenp{j}(\mdposgen))[\force_j(\mdposgen)]_1
  \molli(x-\mdpos[\no]{j}),
\end{align*}
where $[\force_j(\mdposgen)]_1$ is the $x_1$ component of
$\force_j(\mdposgen)$.

For the purpose of computing an approximation of
\begin{align*}
  \frac{1}{\tend}
  \E\left[\int_1^\tend
    \sumall{j}\sum_{k=1}^3\diffumd_{j,k}\otimes\diffumd_{j,k}\right]
\end{align*}
it is not practical to postpone the differentiation of the mollifier
with respect to the space varible.
Using the choice~\eqref{eq_app:mollifier}, the gradient of the
mollifier can be expressed in terms of the mollifier itself as
\begin{align*}
  \dx\molli(x-\mdpos[\no]{j}) & = 
  \frac{-1}{\epsilon^2}\molli(x-\mdpos[\no]{j})
  \left(\left[x-\mdpos[\no]{j}\right]_1,0,0\right)^T.
\end{align*}
Then the expression for $\dX \pfen(x;\mdposgen)$ becomes
\begin{align*}
  \dX \pfen(x;\mdposgen) & = 
  \left(\frac{\potenp{j}(\mdposgen)}{\epsilon^2}
    \left(\left[x-\mdpos[\no]{j}\right]_1,0,0\right)^T
    - \frac{1}{2}\force_j(\mdposgen) \right)
  \molli(x-\mdpos[\no]{j})
  \\ & \quad
  -\frac{1}{2} \sumneq{i}{j} 
  \pairforce_{ij}(\mdposgen)\molli(x-\mdpos[\no]{i})
\end{align*}
and, using the diffusion component~\eqref{eq_app:diffumd_j},
\begin{align*}
  \sum_{k=1}^3
  \diffumd_{j,k}(x;\mdposgen)\diffumd_{j,k}(y;\mdposgen)
  & = 
  2\kb\temp\biggl(
  p_j(x,y;\mdposgen) + 
  q_j(x,y;\mdposgen)
  \biggr),
\end{align*}
where
\begin{align*}
  p_j(x,y;\mdposgen) & = 
  \left(\frac{\potenp{j}(\mdposgen)}{\epsilon^2}\right)^2
  \left[x-\mdpos[\no]{j}\right]_1 \left[y-\mdpos[\no]{j}\right]_1
  \molli(x-\mdpos[\no]{j})\molli(y-\mdpos[\no]{j})\\
  & \quad - 
  \frac{\potenp{j}(\mdposgen)}{2\epsilon^2} 
  [x-\mdpos[\no]{j}]_1 \molli(x-\mdpos[\no]{j})
  \biggl(
  [\force_j(\mdposgen)]_1 \molli(y-\mdpos[\no]{j})
  + \!\!\sumneq{i}{j}
  [\pairforce_{ij}(\mdposgen)]_1 \molli(y-\mdpos[\no]{i})
  \biggr)
  \\
  & \quad - 
  \frac{\potenp{j}(\mdposgen)}{2\epsilon^2} 
  [y-\mdpos[\no]{j}]_1 \molli(y-\mdpos[\no]{j})
  \biggl(
  [\force_j(\mdposgen)]_1 \molli(x-\mdpos[\no]{j})
  + \!\!\sumneq{i}{j}
  [\pairforce_{ij}(\mdposgen)]_1 \molli(x-\mdpos[\no]{i})
  \biggr)
  \intertext{and}
  q_j(x,y;\mdposgen) & = 
  \frac{1}{4}
  \biggl( 
  \force_j(\mdposgen) \molli(x-\mdpos[\no]{j})
  + \sumneq{i}{j} \pairforce_{ij}(\mdposgen) 
  \molli(x-\mdpos[\no]{i}) 
  \biggr) \\
  & \quad\quad
  \cdot
  \biggl( 
  \force_j(\mdposgen) \molli(y-\mdpos[\no]{j})
  + \sumneq{i}{j} \pairforce_{ij}(\mdposgen) 
  \molli(y-\mdpos[\no]{i}) 
  \biggr).
\end{align*}